\documentclass[12pt,twoside,a4paper]{article}
\usepackage[utf8]{inputenc}
\usepackage{amsmath, amssymb, mathtools, enumerate}
\usepackage{mathrsfs} 
\usepackage{amsfonts, amsthm, stmaryrd, fullpage}
\usepackage{authblk}
\usepackage{tikz-cd}
\usepackage[german, english]{babel}
\usepackage[nottoc,numbib]{tocbibind}

\usepackage{needspace}
\usepackage{color}
\usepackage{hyperref}
\hypersetup{colorlinks=true,
linktoc=all,
citecolor=blue,
filecolor=black,
linkcolor=blue,
urlcolor=blue}

\newtheorem{mthm}{Theorem}[section]

\newtheorem{thm}{Theorem}[subsection]
\newtheorem{lemma}[thm]{Lemma}
\newtheorem{cor}[thm]{Corollary}
\newtheorem{prop}[thm]{Proposition}

\theoremstyle{definition} 
\newtheorem{mydef}[thm]{Definition}
  
\newtheorem{example}[thm]{Example}

\newtheorem{mmydef}[mthm]{Definition}

\theoremstyle{remark}
\newtheorem{rmk}[thm]{Remark}

\newcommand\act{{\rm act}}
\newcommand\ad{{\rm ad}}
\newcommand\Aff{{\rm Aff}}
\newcommand\Affd{{\rm Affd}}

\newcommand\alg{{\rm alg}}
\newcommand\an{{\rm an}}

\newcommand\Ban{{\rm Ban}}

\newcommand\can{{\rm can}}

\newcommand\Frac{{\rm Frac}}

\newcommand\Hom{{\rm Hom}}
\newcommand\id{{\rm id}}

\newcommand\Perf{{\rm Perf}}

\newcommand\Rat{{\rm Rat}}

\newcommand\res{{\rm res}}

\newcommand\Sets{{\rm Sets}}

\newcommand\spc{{\rm sp}}

\newcommand\Spa{{\rm Spa}}

\newcommand\Spec{{\rm Spec}} 
\newcommand\Spf{{\rm Spf}}

\begin{document}

\title{Quasi-coherent sheaves in the geometry of adic spaces}
\author{Dimitri Dine \\ 
\small University of California San Diego \\ La Jolla, United States \\
\text{ddine@ucsd.edu} \\ \small \href{https://orcid.org/0000-0002-8340-5358}{ORCiD: 0000-0002-8340-5358}}
\date{}
\maketitle

\begin{abstract}We develop several aspects of a theory of non-derived quasi-coherent sheaves on analytic adic spaces, with the aim of doing EGA-style geometry with adic spaces. For adic spaces equipped with a global pseudo-uniformizer (which we call Tate adic spaces) our notion of quasi-coherent sheaf extends the notion of a pseudo-coherent sheaf on an analytic adic space studied by Kedlaya and Liu. We construct an analog \begin{equation*}(\mathcal{A}, \mathcal{A}^{+})\mapsto \underline{\Spa}_{S}(\mathcal{A}, \mathcal{A}^{+})\end{equation*}of the relative $\Spec$ functor from algebraic geometry in the context of quasi-coherent algebras on Tate adic spaces and show that this construction characterizes locally affinoid morphisms of Tate adic spaces $f: X\to S$ in the same way as the usual relative $\Spec$ construction characterizes affine morphisms of schemes. We also use our relative $\Spa$ construction to give an explicit description of Huber's universal compactifications in the case of locally affinoid morphisms. Finally, we relate our results on locally affinoid morphisms to affine morphisms of formal schemes via the theory of formal models. \end{abstract}

\tableofcontents

\section{Introduction}

\subsection{Motivation and first theorem}

It is a widely known fact in analytic geometry in general, and in the geometry of adic spaces in particular, that a well-behaved theory of quasi-coherent analytic sheaves must be of a derived nature. This philosophy is dictated by the pathological behaviour of the non-derived complete tensor product of Banach modules and algebras, which precludes us from setting up a general formalism of non-derived quasi-coherent sheaves on analytic spaces (e.g., adic spaces) with kernels, cokernels, direct images and pullbacks under more-or-less arbitrary morphisms, etc. 

However, there exist natural restricted classes of quasi-coherent sheaves which are reasonably well-behaved such as \textit{coherent sheaves} on locally Noetherian adic spaces, \textit{vector bundles} on arbitrary analytic adic spaces and, more generally, \textit{pseudocoherent sheaves} on analytic adic spaces, which were introduced by Kedlaya and Liu in \cite{Kedlaya-Liu2}. Another such class of examples are the \textit{Banach bundles} introduced by Kedlaya, which can be regarded as infinite-dimensional generalizations of vector bundles. 

In this paper, our aim is to place all these examples under one umbrella notion of non-derived quasi-coherent sheaf and, more importantly, to consider yet another class of examples where the notion of (non-derived) quasi-coherence can be studied: Direct images of all the above examples under locally affinoid morphisms, a class of morphisms which can be considered as a natural analog of affine morphisms from algebraic geometry (see Definition~\ref{Locally affinoid, introduction} below).  
\begin{mthm}[Example~\ref{Banach bundles}, Example~\ref{Pseudocoherent}, Proposition~\ref{Direct images}]\label{Metatheorem}Consider adic spaces $X$ equipped with a global pseudo-uniformizer. There exists a notion of (non-derived) quasi-coherent sheaves of topological Banach modules $\mathcal{F}$ (and algebras) on such adic spaces which generalizes pseudocoherent sheaves and Banach bundles and such that for every locally affinoid morphism \begin{equation*}f: X\to S\end{equation*}and every pseudocoherent sheaf or Banach bundle $\mathcal{F}$, the direct image $f_{\ast}\mathcal{F}$ is again a quasi-coherent sheaf.  \end{mthm}
In the future, we hope to relate the notion of non-derived quasi-coherence and the examples studied in this paper to the derived categories of quasi-coherent sheaves based on the formalism of condensed mathematics of Clausen and Scholze \cite{Condensed}, \cite{ScholzeAnalytic}. We hope that the quasi-coherent sheaves studied in this paper could become useful as examples of $0$-truncations of quasi-coherent complexes in the general theory. 

\subsection{Three classes of morphisms}

In the meantime, our goal is to use quasi-coherent sheaves to give characterizations, in the spirit of EGA, of several classes of morphisms of analytic adic spaces: The locally affinoid morphisms, the locally affinoid proper morphisms (not necessarily of finite type) and the finite morphisms. The class of locally affinoid morphisms, which was already referred to in Theorem~\ref{Metatheorem}, generalizes finite (or, more generally, pro-finite) morphisms of not necessarily Noetherian adic spaces as well as all morphisms between affinoid adic spaces. It can be seen as a natural (even though, not so faithful, see below) analog of affine morphisms of schemes. 
\begin{mmydef}[Locally affinoid morphisms]\label{Locally affinoid, introduction}A morphism \begin{equation*}f: X\to S\end{equation*}of adic spaces is called locally affinoid if there exists an affinoid open cover $(U_{i})_{i}$ of $S$ such that the pre-image $f^{-1}(U_{i})$ is affinoid for all $i$. \end{mmydef}
We caution the reader that, unlike affine morphisms of schemes, locally affinoid morphisms need not be ``globally affinoid'' (i.e., not every affinoid open subspace of the target has affinoid pre-image), as is illustrated by the example of the open immersion \begin{equation*}V=\{\, \vert x\vert=1\,\}\cup\{\, \vert y\vert=1\,\}\hookrightarrow \mathbb{B}_{K}^{2}=\Spa(K\langle x, y\rangle, K^{\circ}\langle x, y\rangle)\end{equation*}considered by Scholze and Weinstein \cite{Berkeley}, Example~9.1.2 (see Remark~\ref{A non-adapted affinoid adic space} and Remark~\ref{Locally affinoid morphisms and kernels} for detailed discussions of this example in the context of quasi-coherent sheaves). However, the failure of the pre-image of an affinoid open to be affinoid can be characterized in terms of properties of the quasi-coherent sheaf $f_{\ast}\mathcal{O}_{X}$ (Corollary~\ref{Locally affinoid 3}) or, more precisely, of the pair \begin{equation*}(f_{\ast}\mathcal{O}_{X}, f_{\ast}\mathcal{O}_{X}^{+}).\end{equation*}This pair is an example of what we call a \textit{quasi-coherent sheaf of adic Banach algebras}, which we define to be a pair $(\mathcal{A}, \mathcal{A}^{+})$ consisting of a quasi-coherent sheaf of (Banach) algebras $\mathcal{A}$ and a subsheaf \begin{equation*}\mathcal{A}^{+}\subset \mathcal{A}\end{equation*}``of integral elements'' which plays the same role in this definition as the ring of integral elements in the definition of a Huber pair. Our main result on locally affinoid morphisms can then be formulated as follows. 
\begin{mthm}[Theorem~\ref{The affinoid equivalence}]\label{Main theorem 1}\label{Relative Spa}For every adic space $S$ equipped with a global pseudo-uniformizer, there is an anti-equivalence of categories \begin{equation*}(\mathcal{A}, \mathcal{A}^{+})\mapsto \underline{\Spa}_{S}(\mathcal{A}, \mathcal{A}^{+})\end{equation*}between the category of quasi-coherent sheaves of adic Banach algebras and the category of adic spaces $f: X\to S$ over $S$ which are locally affinoid over $S$. 

The quasi-inverse is the functor which sends a locally affinoid morphism $f: X\to S$ to the quasi-coherent sheaf of adic Banach algebras $(f_{\ast}\mathcal{O}_{X}, f_{\ast}\mathcal{O}_{X}^{+})$. \end{mthm}
For locally affinoid morphisms which are moreover proper (where we define properness and partial properness by the valuative criterion, with no finite-type assumptions), keeping track of the subsheaf $\mathcal{A}^{+}$ turns out to be unnecessary since there is a canonical such subsheaf $\overline{\mathcal{A}^{+}}$. 
\begin{mthm}[Theorem~\ref{The affinoid equivalence 2}, Theorem~\ref{The affinoid equivalence 3}]\label{Main theorem 2}Let $S$ be an adic space equipped with a global pseudo-uniformizer. The functors \begin{equation*}\mathcal{A}\mapsto \underline{\Spa}_{S}(\mathcal{A}, \overline{\mathcal{A}^{+}})\end{equation*}and \begin{equation*}(f: X\to S)\mapsto f_{\ast}\mathcal{O}_{X}\end{equation*}are quasi-inverse anti-equivalences of categories between the category of quasi-coherent sheaves of sheafy Banach algebras on $S$ and the category of adic spaces $f: X\to S$ over $S$ which are locally affinoid and proper over $S$. 

Moreover, the functors restrict to quasi-inverse anti-equivalences between the category of quasi-coherent sheaves of algebras of finite type and the category of adic spaces $f: X\to S$ finite over $S$. \end{mthm}
The last part of the theorem generalizes a result of Huber, \cite{Huber0}, Korollar~3.12.12, in the case when the rings of local sections of the adic spaces in question have Noetherian rings of definition. We note that locally affinoid and proper morphisms in our sense need not be finite since we define (partially) proper morphisms by a valuative criterion and do not assume them to be (locally) of finite type. 

\subsection{Application to Huber's compactifications}

We relate Theorem~\ref{Main theorem 2} to Huber's notion of universal compactifications and get an explicit description of the universal compactification in the case of locally affinoid morphisms. 
\begin{mthm}[Theorem~\ref{Compactifications}]\label{Main theorem 3}Let $S$ be an adic space equipped with a global pseudo-uniformizer and let $f: X\to S$ be a locally affinoid morphism of adic spaces. Then the equivalences from Theorem~\ref{Main theorem 1} and Theorem~\ref{Main theorem 2} give rise to a commutative triangle of analytic adic spaces\begin{center}\begin{tikzcd}X\arrow{r}{j} \arrow{d}{f} & \overline{X}=\underline{\Spa}_{S}(f_{\ast}\mathcal{O}_{X}, \overline{(f_{\ast}\mathcal{O}_{X})^{+}}) \arrow{dl}{\overline{f}} \\ S, \end{tikzcd}\end{center}in which the morphism $\overline{f}$ is a locally affinoid proper morphism, $j$ is a locally affinoid injective morphism such that the induced map of complete topological rings $\mathcal{O}_{\overline{X}}\to j_{\ast}\mathcal{O}_{X}$ is an isomorphism and every point of $\overline{X}$ is a specialization of a point of $j(X)$. Moreover, the triangle $(\overline{X}, \overline{f}, j)$ enjoys the following universal property: For any commutative triangle of morphisms of analytic adic spaces \begin{center}\begin{tikzcd}X \arrow{r}{h} \arrow{d}{f} & Y \arrow{dl}{g} \\ S \end{tikzcd}\end{center}with $g$ locally affinoid and proper, there exists a unique morphism $i: \overline{X}\to Y$ making the diagram \begin{center}\begin{tikzcd}X\arrow{dd}{f} \arrow[hook]{dr}{j} \arrow{drr}{h} & & \\ & \overline{X} \arrow{r}{i} \arrow{dl}{\overline{f}} & Y \arrow{dll}{g} \\ S \end{tikzcd}\end{center}commute. The triangle $(\overline{X}, \overline{f}, j)$ is a universal compactification of $f$ in the sense of Huber if and only if $j$ is an open immersion if and only if the morphism $f$ is rationalizing (a condition which generalizes Huber's morphisms of ${}^{+}$weakly finite type, see Definition~\ref{Rationalizing morphism}).\end{mthm}
We note that Huber's theorem on the existence of universal compactifications for separated, taut morphisms of ${}^{+}$weakly finite type (which is proved in \cite{Huber3} in the setting of locally Noetherian analytic adic spaces) generalizes to non-Noetherian situations, and such a generalization (which applies even to some adic spaces that are not analytic) was recently proved by Solodov~\cite{Solodov25}. The main content of Theorem~\ref{Main theorem 3} is thus the explicit description of the compactification in terms of the direct image sheaf $f_{\ast}\mathcal{O}_{X}$ and its subsheaf of integral elements, as well as the general universal property of the triangle $(\overline{X}, \overline{f}, j)$, which applies beyond the case when $f$ is of ${}^{+}$weakly finite type or rationalizing. We note that these latter results are new even for morphisms of weakly finite type between locally Noetherian analytic adic spaces, i.e., even in the setting of Huber's original result. 

\subsection{Relation to other work}

We discuss the relation of the concepts appearing in this paper to similar concepts which appeared in the literature on analytic geometry. 

\begin{enumerate}[(1)]\item First of all, there is a general notion of derived categories of quasi-coherent sheaves on analytic spaces in the formalism of condensed mathematics of Clausen and Scholze~\cite{Condensed}, \cite{ScholzeAnalytic}, which has led to fully-fledged $6$-functor formalisms on analytic adic spaces and small $v$-stacks due to Mann~\cite{Mann22} and Anschütz, LeBras and Mann~\cite{ALBM24}. As we already mentioned in this introduction, we hope to study the relationship between this paper and those general theories in our future work. 

\item Another framework for studying derived categories of quasi-coherent sheaves in analytic geometry, based on categories of ind-Banach modules and complete bornological modules, was proposed by Ben-Bassat, Kelly and Kremnizer~\cite{BBKK24}. This framework builds upon ideas from previous work of the same authors and Bambozzi, for example, \cite{Bambozzi-Ben-Bassat16}, \cite{Ben-Bassat-Kremnizer17}, \cite{BBBK19} and \cite{Bambozzi-Kremnizer20}. The systematic use of categories of (ind-)Banach modules in these papers provided part of the inspiration for our notion of quasi-coherent sheaves of (topological) Banach modules. 

\item Various special cases of our notion of quasi-coherent sheaves have appeared across the literature in both $p$-adic and complex-analytic geometry. We already mentioned that our work was motivated by the notion of pseudocoherent sheaves on analytic adic spaces studied by Kedlaya-Liu~\cite{Kedlaya-Liu2} as well as by the notion of Banach bundles (``infinite-dimensional vector bundles'') introduced by Kedlaya. We remark on some other concrete examples of analytic quasi-coherent sheaves appearing in the literature and their relation to our work. 
\begin{itemize}\item On classical rigid-analytic spaces, a notion of quasi-coherent sheaves has been considered by Conrad~\cite{Conrad06} who defined quasi-coherent sheaves on a rigid-analytic space $X$ as sheaves of $\mathcal{O}_{X}$-modules which, locally for the G-topology on $X$, can be written as inductive limits of coherent sheaves. The main difference between this notion and the notion considered in this paper is that the quasi-coherent sheaves of \cite{Conrad06} are sheaves of abstract modules (not necessarily equipped with any complete norm or topology), while the local sections of our quasi-coherent sheaves are (topological) Banach modules. For example, Conrad's formalism includes the polynomial sheaf $\mathcal{O}_{S}[T_1,\dots, T_n]$ on a rigid space $S$ while our formalism includes the convergent power series sheaf $\mathcal{O}_{S}\langle T_1,\dots, T_n\rangle$. In a similar vein, Conrad constructs a relative $\Spec$ functor, denoted $\Spec^{\an}$, on a subcategory of quasi-coherent sheaves $\mathcal{A}$ (namely, those quasi-coherent sheaves which arise as quotients of the polynomial sheaves $\mathcal{O}_{S}[T_1,\dots, T_n]$) satisfying a universal property analogous to our Theorem~\ref{Universal property of the relative spectrum}, and the main difference between this functor and our relative $\Spa$ is that the rigid space $\Spec^{\an}(\mathcal{A})$ is ``locally affine'' over the base $S$ while our relative $\Spa$ of a quasi-coherent adic Banach algebra is locally affinoid over $S$. In particular, the most basic example of the $\Spec^{\an}$ of a quasi-coherent sheaf in the sense of \cite{Conrad06} is the relative (rigid-analytic) affine space \begin{equation*}\mathbb{A}_{S}^{n,\an}=\Spec^{\an}(\mathcal{O}_{S}[T_1,\dots, T_n])\end{equation*}while the most basic example of our relative $\Spa$ is the relative $n$-dimensional closed unit ball \begin{equation*}\mathbb{B}_{S}^{n}=\underline{\Spa}_{S}(\mathcal{O}_{S}\langle T_1,\dots, T_n\rangle, \mathcal{O}_{S}^{+}\langle T_1,\dots, T_n\rangle)\end{equation*}over $S$. However, the two constructions coincide for coherent sheaves of algebras over a quasi-separated rigid space $S$, i.e., they offer two different characterizations of finite morphisms of quasi-separated rigid-analytic spaces. 
\item For adic spaces $S$ whose Huber rings of local sections $\mathcal{O}_{S}(U)$ (for $U\subseteq S$ affinoid open) admit Noetherian rings of definition and coherent (i.e., finite-type) sheaves of algebras $\mathcal{A}$ over $S$, the special case of our relative adic spectrum $\underline{\Spa}_{S}(\mathcal{A}, \overline{\mathcal{A}^{+}})$ was constructed by Huber in \cite{Huber0}, Prop.~3.12.11. For adic spaces satisfying this very strong Noetherian assumption, Huber also proved the last assertion of Theorem \ref{Main theorem 2} (about coherent sheaves of algebras and finite morphisms), see \cite{Huber0}, Korollar~3.12.12. 
\item For general locally Noetherian analytic adic spaces, the relative $\Spa$ construction was discussed by Ludwig in her lecture notes on adic spaces \cite{Ludwig20}, \S8. It was this discussion of the relative adic spectrum in loc.~cit.~that provided the initial motivation for the conditions on the ``sheaf of integral elements'' \begin{equation*}\mathcal{A}^{+}\subset \mathcal{A}\end{equation*}in our definition of a quasi-coherent sheaf of (topological) adic Banach $(\mathcal{O}_{S}, \mathcal{O}_{S}^{+})$-algebras (Definition~\ref{Sheaf of adic Banach algebras}). 
\item In complex-analytic geometry, there is a notion of quasi-coherent analytic sheaves (on a complex-analytic space) which was introduced by Ramis and Ruget in \cite{Ramis-Ruget74} and used in the works \cite{Putinar86}, \cite{Eschmeier-Putinar84}, \cite{Eschmeier-Putinar96} of Putinar and Eschmeier-Putinar on spectral theory. These are sheaves $\mathcal{F}$ of complete topological modules on a complex-analytic space $X$ with $\mathcal{F}(U)$ a nuclear Fréchet space for every open Stein subset $U\subseteq X$ which satisfy the same quasi-coherence condition as in this paper along with an additional transversality condition. The transversality condition ensures that every such quasi-coherent sheaf $\mathcal{F}$ on a Stein space $X$ is of the form \begin{equation*}\widetilde{M}=M\widehat{\otimes}_{\mathcal{O}_{X}(X)}\mathcal{O}_{X}\end{equation*}for some Fréchet module $M$ over $\mathcal{O}_{X}(X)$ (see \cite{Ramis-Ruget74}, Proposition~9). 
\item An analogous transversality condition for (ind-)Banach modules also appears in the work of Ben-Bassat and Kremnizer \cite{Ben-Bassat-Kremnizer17}, where the authors define RR-quasi-coherent modules over a Banach ring by means of the condition. In this context, we also note that the proof of \cite{Ben-Bassat-Kremnizer17}, Lemma~5.34 (respectively, of \cite{Xia22}, Theorem~3.17) applies also to Tate adic spaces and shows that for every affinoid Tate adic space $X$ and every Banach $\mathcal{O}_{X}(X)$-module $M$ which is transversal to all rational localizations (that is, for every rational subset $U\subseteq X$ the canonical morphism \begin{equation*}M\widehat{\otimes}_{\mathcal{O}_{X}(X)}^{\mathbb{L}}\mathcal{O}_{X}(U)\to M\widehat{\otimes}_{\mathcal{O}_{X}(X)}\mathcal{O}_{X}(U)\end{equation*}is an isomorphism), the pre-sheaf of (topological) Banach modules $M^{\vartriangle}$ on $X$ is a (quasi-coherent) sheaf of (topological) Banach modules on $X$. 
\item In the work of Xia~\cite{Xia22}, several notions related to the present paper are considered in the setting of Berkovich-analytic spaces. In particular, Xia uses the same transversality condition as in the work of Ben-Bassat and Kremnizer~\cite{Ben-Bassat-Kremnizer17} to define a notion of quasi-coherent sheaves on a certain kind of Berkovich $k$-analytic spaces over a nonarchimedean field $k$, called Liu $k$-analytic spaces (see, in particular, Definition~3.9 in loc.~cit.). As is the case with the RR-quasi-coherent sheaves of \cite{Ben-Bassat-Kremnizer17}, Xia's notion becomes a special case of our notion of quasi-coherent sheaves when translated into the setting of Tate adic spaces. Moreover, in \cite{Xia22}, Definition~4.2, Xia singles out a class of quasi-coherent sheaves of Banach algebras on a separated Berkovich $k$-analytic space, called quasi-coherent sheaves of Liu $k$-algebras, for which he constructs an analog of the relative $\Spec$ functor, denoted $\underline{\textrm{Sp}}_{X}$, and proves an anti-equivalence of categories analogous to Theorem~\ref{Main theorem 1}, see \cite{Xia22}, Theorem~1.2. The morphisms of Berkovich spaces characterized by this anti-equivalence are the Liu morphisms in the title of \cite{Xia22}. 
\item In this context, let us also mention that the analog of our locally affinoid morphisms of adic spaces in Berkovich geometry appears in the work of Daylies~\cite{Daylies21} under the name ``almost affinoid morphisms'' (``morphismes presque affinoïdes'' in the French original). 
\end{itemize}
\item Finally, as we already mentioned in our discussion following that theorem, the existence statement for Huber's universal compactifications in Theorem~\ref{Main theorem 3} can be seen as a special case of the very general compactification result proved by Solodov~\cite{Solodov25}. On the other hand, the explicit description of the compactification in terms of a quasi-coherent algebra given in Theorem~\ref{Main theorem 3} is original to our work. It would be interesting to study more closely the relationship between the topological space $\Spa(X, S)$ appearing in Solodov's work (see, for example, \cite{Solodov25}, Theorem~1.2) and our relative $\Spa$ construction for quasi-coherent sheaves of adic Banach algebras. 

\end{enumerate}

\subsection{Acknowledgements}The author would like to express his heartfelt gratitude to his advisor, Kiran Kedlaya, for numerous very helpful discussions, valuable feedback on various drafts of this paper as well as for his continuous encouragement, advice and support. 

The author is grateful to Jack J Garzella, Ryo Ishizuka, Jack Kelly, Kobi Kremnizer, Peter Scholze, Kazuma Shimomoto and Eva Viehmann for their comments on previous versions of this paper. In particular, the author would like to thank Peter Scholze for pointing out an important mistake in an earlier version. 

The author is also grateful to Johannes Anschütz, Lucas Mann and Nathan Wenger for very helpful discussions on notions of quasi-coherent sheaves in analytic geometry. 

This paper was finished when the author was supported by Jane Street under the Jane Street Graduate Research Fellowship 2026. 

\subsection{AI disclosure}The contents of this paper predate the AI-in-mathematics era. All results and proofs in the paper are the author's own. While the author was editing a pre-final draft of the paper, Claude Fable~5.1 was used to conduct several literature searches as well as to proofread the paper and point out some mistakes, which were subsequently corrected by the author. 

\section{Definition and examples of quasi-coherent sheaves}

\subsection{Tate adic spaces}

Recall that an adic space is called analytic if it can be covered by affinoid open subspaces which are adic spectra of Tate Huber pairs. An adic space is called locally Noetherian if it admits an open cover by affinoid subspaces $\Spa(A, A^{+})$, where $A$ is a strongly rigid-Noetherian Huber ring (see \cite{Zavyalov22}, Definition~2.8; when $A$ is Tate these are just the strongly Noetherian Tate rings of Huber). Finally, an adic space is called uniform (respectively, sousperfectoid, respectively, diamantine) if it admits a cover by affinoid open subspaces $\Spa(A, A^{+})$, where $A$ is a stably uniform Tate ring (respectively, a sousperfectoid Tate ring, respectively, a diamantine Tate ring in the sense of \cite{Hansen-Kedlaya}). Note that any of the above properties of affinoid adic spaces is inherited by rational subsets, so an adic space $X$ is analytic (resp., locally Noetherian, resp., uniform, resp., sousperfectoid, resp., diamantine) if and only if there exists a basis $\mathcal{B}$ for its topology whose members are affinoid adic spaces $\Spa(A, A^{+})$ with $A$ a Tate ring (resp., strongly rigid-Noetherian, resp., a stably uniform Tate ring, resp., a sousperfectoid Tate ring, resp., a diamantine Tate ring). 
\begin{mydef}[Tate adic spaces]\label{Tate adic spaces}A global pseudo-uniformizer of an adic space $X$ is an element $\varpi_{X}\in \mathcal{O}_{X}^{+}(X)$ such that for any affinoid open subspace $U$ of $X$ the image of $\varpi_{X}$ in $\mathcal{O}_{X}(U)$ is a topologically nilpotent unit of $\mathcal{O}_{X}(U)$. A Tate adic space is a pair $(X, \varpi_{X})$, where $X$ is a (necessarily analytic) adic space and $\varpi_{X}\in \mathcal{O}_{X}^{+}(X)$ is a global pseudo-uniformizer of $X$. A Tate adic space $(X, \varpi_{X})$ is said to be locally Noetherian (respectively, uniform, respectively, perfectoid, respectively, sousperfectoid, respectively, diamantine) if the underlying adic space $X$ has this property. 

A morphism $f: (Y, \varpi_{Y})\to (X, \varpi_{X})$ of Tate adic spaces is a morphism of adic spaces $f: Y\to X$ such that the image of $\varpi_{X}$ under the map $\mathcal{O}_{X}^{+}(X)\to \mathcal{O}_{Y}^{+}(Y)$ is equal to $\varpi_{Y}$.\end{mydef}
\begin{rmk}[Tate adic spaces as relative adic spaces]\label{Tate adic spaces as relative adic spaces}For every sheafy Tate Huber pair $(A, A^{+})$ with pseudo-uniformizer $\varpi\in A^{+}$ of $A$ the assignment $X\mapsto (X, \varpi)$ defines a faithful functor from the category of adic spaces over $S=\Spa(A, A^{+})$ to the category of Tate adic spaces. In particular, the category of rigid-analytic varieties over any nonarchimedean field $K$ admits a natural faithful functor to the category of Tate adic spaces. Conversely, the category of Tate adic spaces is the same as the category of adic spaces over the affinoid adic space\begin{equation*}S=\Spa(\mathbb{Z}((T)), \mathbb{Z}[[T]]),\end{equation*}where $\mathbb{Z}[[T]]$ is endowed with the $T$-adic topology, where $\mathbb{Z}((T))=\mathbb{Z}[[T]][T^{-1}]$ and where for a Tate adic space $(X, \varpi)$ and $U\subseteq X$ affinoid open the morphism of Huber pairs\begin{equation*}(\mathbb{Z}((T)), \mathbb{Z}[[T]])\to (\mathcal{O}_{X}(U), \mathcal{O}_{X}^{+}(U))\end{equation*}is the map sending the variable $T$ to $\varpi_{X}$.\end{rmk}We often denote a Tate adic space $(X, \varpi_{X})$ just by the letter $X$. 
\begin{mydef}If $\mathcal{P}$ is a property of adic spaces, then a Tate adic space $(X, \varpi_{X})$ is said to satisfy property $\mathcal{P}$ if its underlying adic space $X$ has $\mathcal{P}$. Similarly, if $\mathcal{P}$ is a property of morphisms of adic spaces, then a morphism of Tate adic spaces $f: (Y, \varpi_{Y})\to(X, \varpi_{X})$ is said to satisfy $\mathcal{P}$ if the underlying morphism of adic spaces $f: Y\to X$ satisfies $\mathcal{P}$.\end{mydef}
In particular, we can apply the above definition to properties of adic spaces such as being uniform, locally Noetherian, perfectoid, sousperfectoid or diamantine, and we can apply it to properties of morphisms such as being adic, quasi-compact or locally affinoid (Definition~\ref{Locally affinoid morphisms} below).    

\subsection{Definition~of quasi-coherent sheaves and adapted open subsets}
 
Recall that, just as for maps between Banach spaces over $\mathbb{R}$ and $\mathbb{C}$, a linear map between Banach spaces over a nonarchimedean field $K$ is bounded if and only if it is continuous. By contrast, a continuous map between general Tate Banach rings may fail to be bounded. For example, if $K$ is a nonarchimedean field with absolute value $\vert\cdot\vert$ and $s\in (0, 1)$, then the absolute value $\vert\cdot\vert^{s}$ on $K$ defines the same topology as $\vert\cdot\vert$, but is not bounded with respect to $\vert\cdot\vert$. Note that, in the definition~of a Banach $A$-module $M$ over a Banach ring $A$, we require that the map $A\times M\to M$ actually be bounded and not just continuous. To ensure that at least all continuous $A$-module maps between Banach $A$-modules $M\to N$ are bounded, it suffices to assume that $A$ has a pseudo-uniformizer $\varpi$ which is multiplicative with respect to the norm $\lVert\cdot\rVert_{A}$ on $A$, i.e., \begin{equation*}\lVert \varpi f\rVert=\lVert\varpi\rVert\lVert f\rVert\end{equation*}for all $f\in A$ (see \cite{Dine25-3}, Corollary~2.8, for sufficiency). To formalize this requirement, we introduce the notion of a topological Banach module over a complete Tate $\mathbb{Z}((T))$-algebra.

Note that for any complete Tate ring $A$ and any topologically nilpotent unit $\varpi$ in $A$ the topology on $A$ can be defined by a norm $\lVert\cdot\rVert_{A}$ such that $\varpi$ is multiplicative with respect to $\lVert\cdot\rVert_{A}$: Indeed, by \cite{Huber0}, Proposition~2.2.7(iii), we can choose a subring $A_{0}$ of $A$ containing $\varpi$ such that $(A_{0}, (\varpi))$ is a pair of Definition~and we can then let $\lVert\cdot\rVert_{A}$ be the canonical extension $\lVert\cdot\rVert_{A_{0},\varpi}$ of the $\varpi$-adic norm on $A_{0}$. Note that a complete Tate ring $A$ equipped with a fixed choice of a topologically nilpotent unit $\varpi$ is the same thing as a complete Huber ring which is also a topological $\mathbb{Z}((T))$-algebra (where $\varpi$ is the image of $T$ in $A$). We call such a topological $\mathbb{Z}((T))$-algebra a complete Tate $\mathbb{Z}((T))$-algebra. 
\begin{mydef}[Topological Banach module over a complete Tate $\mathbb{Z}((T))$-algebra]\label{Topological Banach module}Let $A$ be a complete Tate $\mathbb{Z}((T))$-algebra. A topological Banach $A$-module $M$ is a complete topological $A$-module such that there exists a norm $\lVert\cdot\rVert_{A}$ defining the topology on $A$ such that the image $\varpi$ of $T$ in $A$ is multiplicative with respect to $\lVert\cdot\rVert_{A}$ and an $(A, \lVert\cdot\rVert_{A})$-module norm $\lVert\cdot\rVert_{M}$ on $M$ defining the topology on $M$. A morphism of topological Banach $A$-modules is a continuous homomorphism of topological $A$-modules.\end{mydef}
Note that by \cite{Dine25-3}, Corollary~2.8, any two norms $\lVert\cdot\rVert_{A}$ and $\lVert\cdot\rVert_{A}'$ on $A$ satisfying the requirement of Definition~\ref{Topological Banach module} are bounded-equivalent, as are any two $(A, \lVert\cdot\rVert_{A})$-module norms $\lVert\cdot\rVert_{M}$, $\lVert\cdot\rVert_{M}'$ defining the topology on $M$. We also have an analogous notion of a topological Banach algebra $B$ over a complete Tate $\mathbb{Z}((T))$-algebra $A$. 
\begin{mydef}[Topological Banach algebra over a complete Tate $\mathbb{Z}((T))$-algebra]\label{Topological Banach algebra}Let $A$ be a complete Tate $\mathbb{Z}((T))$-algebra. A topological Banach $A$-algebra $B$ is a complete topological $A$-algebra such that there exists a norm $\lVert\cdot\rVert_{A}$ defining the topology on $A$ such that the image $\varpi$ of $T$ in $A$ is multiplicative with respect to $\lVert\cdot\rVert_{A}$ and an $(A, \lVert\cdot\rVert_{A})$-module norm $\lVert\cdot\rVert_{B}$ on $B$ defining the topology on $B$ which is also a ring norm.\end{mydef}
\begin{rmk}Let $A$ be a complete Tate $\mathbb{Z}((T))$-algebra (with $\varpi$ denoting the image of $T$ in $A$) and let $\lVert\cdot\rVert$ be any norm defining the topology on $A$. Let $(M, \lVert\cdot\rVert_{M})$ be any Banach module over the Banach ring $(A, \lVert\cdot\rVert)$. Then the underlying complete topological $A$-module of $M$ is always a topological Banach $A$-module: Indeed, if $A_{\leq1}$ and $M_{\leq1}$ are the closed unit balls of $A$ and $M$, respectively, for their chosen norms, then the topologies on $A$ and $M$ can also be defined by the canonical extension of the $\varpi$-adic norm on $A_{\leq1}$ and the canonical extension of the $\varpi$-adic norm on $M_{\leq1}$, respectively, and the element $\varpi$ is multiplicative with respect to these norms. We use the term ``topological Banach module'' to distinguish this notion from that of a Banach module $(M, \lVert\cdot\rVert_{M})$ with a fixed norm $\lVert\cdot\rVert_{M}$, the point being again that for an arbitrary fixed norm, continuous maps can fail to be bounded. The same remark applies to (topological) Banach algebras.  \end{rmk}
Note that any complete Tate $\mathbb{Z}((T))$-algebra is naturally a topological Banach $\mathbb{Z}((T))$-algebra and that any topological Banach $A$-algebra over a complete topological $\mathbb{Z}((T))$-algebra $A$ is a complete Tate $\mathbb{Z}((T))$-algebra via the composite map $\mathbb{Z}((T))\to A\to B$.
\begin{mydef}[Complete tensor product of topological Banach modules]\label{Complete tensor product of topological Banach modules}Let $A$ be a complete Tate $\mathbb{Z}((T))$-algebra and let $M, N$ be topological Banach $A$-modules (respectively, topological Banach $A$-algebras). Let $\lVert\cdot\rVert_{A}$ be a norm defining the topology on $A$ such that the image of $T$ in $A$ is multiplicative with respect to $\lVert\cdot\rVert_{A}$ and let $\lVert\cdot\rVert_{M}, \lVert\cdot\rVert_{N}$ be $(A, \lVert\cdot\rVert_{A})$-module norms on $M, N$ defining their respective topologies. We endow the tensor product $M\otimes_{A}N$ with the tensor product seminorm resulting from the $(A, \lVert\cdot\rVert_{A})$-module norms $\lVert\cdot\rVert_{M}$ and $\lVert\cdot\rVert_{N}$ on $M$ and $N$. By the definition~of the tensor product seminorm, if we replace $\lVert\cdot\rVert_{M}$ and $\lVert\cdot\rVert_{N}$ with bounded-equivalent norms, then the resulting tensor product seminorms are bounded-equivalent. Since any two $(A, \lVert\cdot\rVert_{A})$-module norms defining the same topology on $M$ (respectively, on $N$) are bounded-equivalent by \cite{Dine25-3}, Corollary~2.8, this means in particular that the topology on $M\otimes_{A}N$ defined by the tensor product seminorm does not depend on any particular choice of $\lVert\cdot\rVert_{A}$, $\lVert\cdot\rVert_{M}$ and $\lVert\cdot\rVert_{N}$ and we define the complete tensor product \begin{equation*}M\widehat{\otimes}_{A}N\end{equation*}of the topological Banach $A$-modules (respectively, topological Banach $A$-algebras) $M$ and $N$ as the topological Banach $A$-module (respectively, topological Banach $A$-algebra) $M\widehat{\otimes}_{A}N$ obtained by completing $M\otimes_{A}N$ with respect to the said $A$-module or $A$-algebra topology.\end{mydef}
For a fixed norm $\lVert\cdot\rVert_{A}$ on $A$ such that the image of $T$ in $A$ is norm-multiplicative, there is a canonical forgetful functor from Banach $(A, \lVert\cdot\rVert_{A})$-modules (respectively, Banach $(A, \lVert\cdot\rVert_{A})$-algebras) to topological Banach $A$-modules (respectively, topological Banach $A$-algebras) and this functor respects complete tensor products, by definition. Recall from Remark~\ref{Tate adic spaces as relative adic spaces} that for a Tate adic space $(X, \varpi_{X})$ the global pseudo-uniformizer $\varpi_{X}$ of $X$ corresponds to the datum of a morphism of adic spaces $X\to \Spa(\mathbb{Z}((T)), \mathbb{Z}[[T]])$, so for every affinoid open subspace $U$ of $X$ the complete topological ring $\mathcal{O}_{X}(U)$ is naturally a complete Tate $\mathbb{Z}((T))$-algebra. 
\begin{mydef}[Pre-sheaves and sheaves of topological Banach modules]\label{Pre-sheaves and sheaves}Let $(X, \varpi_{X})$ be a Tate adic space. A (pre-)sheaf of topological Banach modules (respectively, of topological Banach algebras) $\mathcal{F}$ on $(X, \varpi_{X})$ (or on $X$) is a (pre-)sheaf of complete topological $\mathcal{O}_{X}$-modules (respectively, of complete topological algebras) such that for every $U\in\Affd_{X}$ the complete topological $\mathcal{O}_{X}(U)$-module (respectively, the complete topological $\mathcal{O}_{X}(U)$-algebra) $\mathcal{F}(U)$ is a topological Banach module (respectively, a topological Banach algebra) over the complete Tate $\mathbb{Z}((T))$-algebra $\mathcal{O}_{X}(U)$. \end{mydef}
We also need the analogous notions of pre-sheaves and sheaves on subcategories of the category of affinoid open subspaces of $X$. 
\begin{mydef}Let $(X, \varpi_{X})$ be a Tate adic space and let $\mathcal{C}$ be a subcategory of the category of affinoid open subspaces of $X$. A (pre-)sheaf of topological Banach modules (algebras) on $\mathcal{C}$ is a (pre-)sheaf of complete topological $\mathcal{O}_{X}\vert_{\mathcal{C}}$-modules (respectively, complete topological $\mathcal{O}_{X}\vert_{\mathcal{C}}$-algebras) such that, for every $U\in \mathcal{C}$, the complete topological module $\mathcal{F}(U)$ is a topological Banach module (respectively, topological Banach algebra) over the complete Tate $\mathbb{Z}((T))$-algebra $\mathcal{O}_{X}(U)$. \end{mydef}
In the following, we denote by $\Aff_{X}$ the category of all affinoid open subspaces of an adic space $X$ and, for every $U\in \Aff_{X}$, we denote by $\Rat_{U}$ the subcategory of $\Aff_{U}$ consisting of rational subsets. 
\begin{example}\label{Defining example}Let $X=\Spa(A, A^{+})$ be an affinoid adic space such that $A=\mathcal{O}_{X}(X)$ is a complete Tate ring. Fix a pseudo-uniformizer $\varpi$ of $A$ and consider the Tate adic space $(X, \varpi)$. Let $M$ be a topological Banach $A$-module over the complete Tate $\mathbb{Z}((T))$-algebra $A$. Then \begin{equation*}\widetilde{M}(U)=M\widehat{\otimes}_{A}\mathcal{O}_{X}(U)\end{equation*}for $U\in \Rat_{X}$ defines a pre-sheaf of topological Banach $\mathcal{O}_{X}\vert_{\Rat_{X}}$-modules on $\Rat_{X}$ (here we use the complete tensor product of topological Banach modules from Definition~\ref{Complete tensor product of topological Banach modules}). \end{example} 
Our notion of quasi-coherent sheaves will be based on the following notion of $\mathcal{F}$-adapted affinoid open subspace.  
\begin{mydef}[$\mathcal{F}$-adapted open subset]\label{Adapted open subsets, definition}Let $(X, \varpi_{X})$ be a Tate adic space and let $\mathcal{F}$ be a pre-sheaf of topological Banach modules (respectively, topological Banach algebras) on $X$. An affinoid open subspace \begin{equation*}U\subseteq X\end{equation*}is called $\mathcal{F}$-adapted if the canonical morphism of pre-sheaves of complete topological $\mathcal{O}_{U}\vert_{\Rat_{U}}$-modules on $\Rat_{U}$\begin{equation*}\widetilde{\mathcal{F}(U)}\to \mathcal{F}\vert_{\Rat_{U}}\end{equation*}is an isomorphism. In other words, for every rational open subspace $V$ of $U$ the canonical morphism \begin{equation*}\can_{\mathcal{F},V\to U}: \mathcal{F}(U)\widehat{\otimes}_{\mathcal{O}_{X}(U)}\mathcal{O}_{X}(V)\to \mathcal{F}(V)\end{equation*}induced by the restriction map $\mathcal{F}(U)\to \mathcal{F}(V)$ is an isomorphism of topological Banach $\mathcal{O}_{X}(V)$-modules.\end{mydef}
We are now ready to state the definition~of quasi-coherent sheaves. 
\begin{mydef}[Quasi-coherent sheaves of topological Banach modules]\label{Quasi-coherent sheaves, definition}A (pre-)sheaf of topological Banach modules (respectively, of topological Banach algebras) $\mathcal{F}$ on $X$ is called a quasi-coherent (pre-)sheaf of topological Banach modules (respectively, of topological Banach algebras) on $X$ if there exists a basis $\mathcal{B}\subseteq \Affd_{X}$ for the topology on $X$ whose members are $\mathcal{F}$-adapted affinoid open subspaces.\end{mydef}
By abuse of terminology, we often just say ``Banach algebras and Banach modules'' when referring to topological Banach algebras and topological Banach modules.
\begin{rmk}Note that, in the situation of the above definition, the restriction map $\mathcal{F}(U)\to \mathcal{F}(V)$ is bounded for any $U, V\in \Affd_{X}$ and any choice of $\mathcal{O}_{X}(U)$-module norms defining the topologies on $\mathcal{F}(U)$ and $\mathcal{F}(V)$ whenever $\mathcal{O}_{X}(U)$ is equipped with a norm defining its topology such that $\varpi_{X}$ is multiplicative with respect to this norm, by \cite{Dine25-3}, Corollary~2.8. \end{rmk}
The following is the prototypical class of examples of quasi-coherent (pre-)sheaves. 
\begin{example}\label{Delta sheaves}Let $X$ be an affinoid Tate adic space and let $M$ be a topological Banach module over $\mathcal{O}_{X}(X)$. We define the quasi-coherent pre-sheaf of topological Banach modules $M^{\vartriangle}$ on $X$ to be the right Kan extension of the pre-sheaf of topological Banach modules $\widetilde{M}$ on $\Rat_{X}$. In other words, if $U\subseteq X$ is a rational subset of $X$, then \begin{equation*}M^{\vartriangle}(U)=M\widehat{\otimes}_{\mathcal{O}_{X}(X)}\mathcal{O}_{X}(U)\end{equation*}and for general open subsets $U$, \begin{equation*}M^{\vartriangle}(U)=\varprojlim_{V\subseteq U}M^{\vartriangle}(V),\end{equation*}where $V$ ranges over the rational subsets of $X$ contained in $U$. 

If we replace the topological Banach $\mathcal{O}_{X}(X)$-module $M$ with a topological Banach $\mathcal{O}_{X}(X)$-algebra $B$, we get a quasi-coherent pre-sheaf of topological Banach $\mathcal{O}_{X}(X)$-algebras $B^{\triangle}$.  \end{example}
\begin{mydef}[Globally quasi-coherent sheaf]\label{Globally quasi-coherent}We call a pre-sheaf of topological Banach modules (respectively, of topological Banach algebras) on an affinoid Tate adic space $X$ globally quasi-coherent if it is of the form $M^{\vartriangle}$ for some topological Banach $\mathcal{O}_{X}(X)$-module $M$ (respectively, of the form $B^{\vartriangle}$ for some topological Banach $\mathcal{O}_{X}(X)$-algebra $B$). \end{mydef}
The following lemma is very useful for working with $\mathcal{F}$-adapted affinoid open subsets. 
\begin{lemma}\label{Subsets of adapted sets are adapted}Let $\mathcal{F}$ be a (pre-)sheaf of topological Banach modules (respectively, of topological Banach algebras) on a Tate adic space $X$ and let $U\subseteq X$ be an $\mathcal{F}$-adapted affinoid open subspace. Then every rational subspace $V$ of $U$ is $\mathcal{F}$-adapted.\end{lemma}
\begin{proof}Consider an inclusion of rational subsets $W\subseteq V$. Since $U$ is $\mathcal{F}$-adapted, we know that the canonical morphisms $\can_{\mathcal{F},W\to U}$ and $\can_{\mathcal{F},V\to U}$ induced by the restriction morphisms $\res_{\mathcal{F}, W\to U}: \mathcal{F}(U)\to \mathcal{F}(W)$ and $\res_{\mathcal{F}, V\to U}: \mathcal{F}(U)\to \mathcal{F}(V)$ are isomorphisms. We have to prove that this entails that $\can_{\mathcal{F}, W\to V}$ is also an isomorphism. To begin with, for any affinoid open subspace $U^{\prime}$ of $U$, we write $\act_{U^{\prime}}$ for the canonical continuous $\mathcal{O}_{X}(U^{\prime})$-module map ($\mathcal{O}_{X}(U^{\prime})$-action) \begin{equation*}\mathcal{O}_{X}(U^{\prime})\widehat{\otimes}_{\mathcal{O}_{X}(U)}\mathcal{F}(U^{\prime})\to \mathcal{F}(U^{\prime}).\end{equation*}Recall that the continuous map \begin{equation*}\mathcal{O}_{X}(W)\widehat{\otimes}_{\mathcal{O}_{X}(U)}\mathcal{O}_{X}(V)\to \mathcal{O}_{X}(W)\end{equation*}is an isomorphism by basic properties of rational localizations. The isomorphism $\can_{\mathcal{F},W\to U}=\act_{W}\circ (\id\widehat{\otimes}\res_{\mathcal{F},W\to U})$ factors as \begin{equation*}\can_{\mathcal{F},W\to U}=\act_{W}\circ (\id\widehat{\otimes}\res_{\mathcal{F},W\to V})\circ (\id\widehat{\otimes}\res_{\mathcal{F},V\to U})=\can_{\mathcal{F},W\to V}\circ (\id\widehat{\otimes}\res_{\mathcal{F},V\to U}).\end{equation*}Therefore, to prove that $\can_{\mathcal{F},W\to V}$ is an isomorphism, it suffices to prove that \begin{equation*}\id\widehat{\otimes}\res_{\mathcal{F},V\to U}: \mathcal{O}(W)\widehat{\otimes}_{\mathcal{O}(U)}\mathcal{F}(U)\to \mathcal{O}(W)\widehat{\otimes}_{\mathcal{O}(U)}\mathcal{F}(V)\end{equation*}is an isomorphism. But $\id\widehat{\otimes}\res_{\mathcal{F},V\to U}$ fits into the canonical commutative diagram \begin{center}\begin{tikzcd} \mathcal{O}(W)\widehat{\otimes}_{\mathcal{O}(U)}\mathcal{F}(U) \arrow{r}{\id\widehat{\otimes}\res_{\mathcal{F},V\to U}} \arrow{dd}{\cong} & \mathcal{O}(W)\widehat{\otimes}_{\mathcal{O}(U)}\mathcal{F}(V) \arrow{d}{\cong} \\ & \mathcal{O}(W)\widehat{\otimes}_{\mathcal{O}(V)}\mathcal{O}(V)\widehat{\otimes}_{\mathcal{O}(U)}\mathcal{F}(V) \arrow{d}{\cong} \\ \mathcal{O}(W)\widehat{\otimes}_{\mathcal{O}(V)}(\mathcal{O}(V)\widehat{\otimes}_{\mathcal{O}(U)}\mathcal{F}(U)) \arrow{r}{\cong} & \mathcal{O}(W)\widehat{\otimes}_{\mathcal{O}(V)}\mathcal{F}(V), \end{tikzcd} \end{center}where the lower horizontal arrow is given by base extension along $\mathcal{O}(V)\to \mathcal{O}(W)$ of the isomorphism $\can_{\mathcal{F},V\to U}$. This shows that $\id\widehat{\otimes}\res_{\mathcal{F},V\to U}$ must be an isomorphism.\end{proof}
\begin{cor}\label{Adapted sets and finite families of sheaves}Let $X$ be a Tate adic space and let $\mathcal{F}_{1},\dots, \mathcal{F}_{n}$ be finitely many quasi-coherent (pre-)sheaves of topological Banach modules (respectively, of topological Banach algebras) on $X$. Then there exists a basis for the topology on $X$ consisting of affinoid open subspaces which are $\mathcal{F}_{i}$-adapted for every $i=1,\dots, n$. If $X$ is affinoid, we can choose $\mathcal{B}$ to consist only of rational subsets.\end{cor}
\begin{proof}We proceed by induction on $n$. Assume that the assertion holds for some $n\geq1$ and let $\mathcal{F}_{1},\dots,\mathcal{F}_{n+1}$ be $n+1$ quasi-coherent pre-sheaves of topological Banach modules on $X$. We have to show that every affinoid open subspace of $X$ can be covered by affinoid open subspaces of $X$ which are $\mathcal{F}_{i}$-adapted for all $i=1,\dots, n+1$ and that these affinoid open subspaces can be chosen to be rational if $X$ is affinoid. Let $\mathcal{B}_{n+1}\subseteq\Affd_{X}$ be a basis consisting of $\mathcal{F}_{n+1}$-adapted affinoid open subspaces. By the induction hypothesis, there exists a basis $\mathcal{B}_{1,\dots,n}\subseteq \Affd_{X}$ such that every $U\in\mathcal{B}_{1,\dots,n}$ is $\mathcal{F}_{i}$-adapted for all $i=1,\dots, n$. It suffices to prove that every $U\in\mathcal{B}_{1,\dots,n}$ can be covered by affinoid open subspaces (respectively, by rational subsets of $X$, if $X$ is affinoid) which are $\mathcal{F}_{i}$-adapted for all $i=1,\dots, n+1$. Given $U\in\mathcal{B}_{1,\dots,n}$, choose an open cover $(U_{j})_{j}$ of $U$ such that $U_{j}\in\mathcal{B}_{n+1}$ for all $j$ and, for every $j$, choose an open cover $(W_{jk})_{k}$ of $U_{j}$ by rational subsets of $U$ (respective, by rational subsets of $X$, if $X$ is affinoid). Then each $W_{jk}$ is also a rational subset of $U_{j}$ (and of $U$). But by virtue of Lemma~\ref{Subsets of adapted sets are adapted} this implies that each $W_{jk}$ is $\mathcal{F}_{i}$-adapted for all $i=1,\dots, n+1$.\end{proof}
\begin{cor}\label{Terminal objects and quasi-coherence}Let $X$ be a Tate adic space and let $\mathcal{F}$ be a sheaf of topological Banach modules on $X$. If $X$ has an open cover $(U_{i})_{i}$ by $\mathcal{F}$-adapted affinoid open subspaces, then $\mathcal{F}$ is a quasi-coherent sheaf of topological Banach modules.\end{cor}
\begin{proof}Let $V\subseteq X$ be an arbitrary affinoid open subspace. We have to prove that $V$ admits a covering by $\mathcal{F}$-adapted affinoid open subspaces. For every $i$, let $(W_{ij})_{j}$ be an open cover of $U_{i}\cap V$ by rational subsets of $U_{i}$. Since each $W_{ij}$ is a rational subset of $U_{i}$, Lemma~\ref{Subsets of adapted sets are adapted} ensures that $W_{ij}$ is $\mathcal{F}$-adapted for all $i, j$.\end{proof}

\subsection{Examples of quasi-coherent sheaves}

With this in mind, let us discuss some classes of examples of quasi-coherent sheaves on adic spaces, some of which already appeared in the literature. Recall that the natural topology on a finitely generated $A$-module $M$ over a topological ring $A$ is the quotient topology induced by some (equivalently, any) $A$-linear surjection $A^{n}\to M$, for some $n\geq1$. This is the finest $A$-module topology that exists on $M$ and if $A$ is a complete Tate ring, then this is the only complete metrizable $A$-module topology on $M$, by Banach's open mapping Theorem~over complete Tate rings.  
\begin{example}[Banach spaces]\label{Banach spaces}For a nonarchimedean field $K$ the underlying topological space of $X=\Spa(K, K^{\circ})$ is a point, so all sheaves of (topological) Banach modules on $X$ (where we choose a pseudo-uniformizer $\varpi\in K$ to turn $X$ into a Tate adic space) are trivially quasi-coherent and the category of quasi-coherent sheaves of (topological) Banach modules on $X$ can be identified with the category of (topological) Banach spaces over $K$. Similarly, the category of (quasi-coherent) sheaves of Banach algebras (respectively, of topological Banach algebras) on $X$ is simply the category of Banach $K$-algebras (respectively, of topological Banach $K$-algebras).\end{example}  
\begin{example}[Vector bundles]\label{Vector bundles}If $(X, \varpi_{X})$ is a Tate adic space and $\mathcal{F}$ is a vector bundle on $X$, then $\mathcal{F}$ is naturally a quasi-coherent sheaf of topological Banach $\mathcal{O}_{X}$-modules. Indeed, there is an open affinoid cover $(U_{i})_{i}$ of $X$ and for every $i$ a finite projective $\mathcal{O}_{X}(U_{i})$-module $M_{i}$ such that $\mathcal{F}(V)=M_{i}\otimes_{\mathcal{O}_{X}(U_{i})}\mathcal{O}_{X}(V)$ for every rational subset $V\subseteq U_{i}$. By considering a split surjection $F\to M_{i}$ from a finite free $\mathcal{O}_{X}(U_{i})$-module, we see that the natural topology on $M_{i}$ can be defined by a complete $\mathcal{O}_{X}(U_{i})$-module norm, with respect to any norm defining the topology on $\mathcal{O}_{X}(U_{i})$ which turns $\varpi_{X}$ into a norm-multiplicative element. For any rational subset $V$ the corresponding tensor product seminorm on $M_{i}\otimes_{\mathcal{O}_{X}(U_{i})}\mathcal{O}_{X}(V)$ is again a complete norm since it is induced by the split surjection $F\otimes_{\mathcal{O}_{X}(U_{i})}\mathcal{O}_{X}(V)\to M_{i}\otimes_{\mathcal{O}_{X}(U_{i})}\mathcal{O}_{X}(V)$. This shows that, for any rational subset $V$ of $U_{i}$, we have \begin{equation*}M_{i}\otimes_{\mathcal{O}_{X}(U_{i})}\mathcal{O}_{X}(V)=M_{i}\widehat{\otimes}_{\mathcal{O}_{X}(U_{i})}\mathcal{O}_{X}(V)\end{equation*}for this choice of norm on $M_{i}$. Consequently, the sheaf of complete topological modules $\mathcal{F}$ is a sheaf of Banach modules and the affinoid open subsets $U_{i}$ are $\mathcal{F}$-adapted. We conclude by Corollary~\ref{Terminal objects and quasi-coherence} that $\mathcal{F}$ is a quasi-coherent sheaf of topological Banach modules.\end{example}
\begin{example}[Banach bundles]\label{Banach bundles}In an online talk at Columbia University in April 2021, Kedlaya introduced the notion of a Banach bundle on an analytic adic space $X$ generalizing the notion of a vector bundle on $X$. In our terminology, Banach bundles on a Tate adic space $(X, \varpi_{X})$ are precisely the quasi-coherent sheaves of topological Banach modules $\mathcal{F}$ on $X$ such that, for some affinoid open cover $(U_{i})_{i}$ of $X$, the restriction of $\mathcal{F}$ to each $U_{i}$ is the sheaf of topological Banach modules $M_{i}^{\vartriangle}$ (see Example~\ref{Delta sheaves}) with $M_{i}$ a projective object in the category of topological Banach $\mathcal{O}_{X}(U_{i})$-modules (equivalently, a direct summand of a topologically free topological Banach $\mathcal{O}_{X}(U_{i})$-module). Kedlaya proved that every Banach bundle on an affinoid Tate adic space $X$ has vanishing higher cohomology and that the category of Banach bundles on an affinoid Tate adic space $X$ is equivalent via the global sections functor to the category of projective (topological) Banach $\mathcal{O}_{X}(X)$-modules. In our terminology, this implies that every affinoid Tate adic space $X$ is $\mathcal{F}$-adapted for every Banach bundle $\mathcal{F}$ on $X$.\end{example}
\begin{example}[Pseudocoherent and coherent sheaves]\label{Pseudocoherent}Let $(X, \varpi_{X})$ be a Tate adic space and let $\mathcal{F}$ be a pseudocoherent sheaf on the analytic adic space $X$ in the sense of Kedlaya and Liu \cite{Kedlaya-Liu2}, \S2.5. Then there exists an open affinoid cover $(U_{i})_{i}$ of $X$ such that, for every $V\in\Rat_{U}$, \begin{equation*}\mathcal{F}(V)=M_{i}\otimes_{\mathcal{O}_{X}(U_{i})}\mathcal{O}_{X}(V),\end{equation*}where $M_{i}$ is a stably pseudocoherent $\mathcal{O}_{X}(U_{i})$-module, i.e., $M_{i}$ is a pseudocoherent $\mathcal{O}_{X}(U_{i})$-module and $M_{i}\otimes_{\mathcal{O}_{X}(U_{i})}\mathcal{O}_{X}(V)$ is complete for its natural topology of finitely generated $\mathcal{O}_{X}(V)$-module for every rational subset $V$ of $U_{i}$. If we choose a norm on $\mathcal{O}_{X}(X)$ defining the topology on $\mathcal{O}_{X}(X)$ and making the element $\varpi_{X}\in \mathcal{O}_{X}(X)$ norm-multiplicative and choose an $\mathcal{O}_{X}(X)$-module norm on $M_{i}$ defining its natural topology, then the corresponding tensor product norm defines the natural topology on $M_{i}\otimes_{\mathcal{O}_{X}(U_{i})}\mathcal{O}_{X}(V)$ for all rational subsets $V$; in particular, with this choice of norm we have $M_{i}\otimes_{\mathcal{O}_{X}(U_{i})}\mathcal{O}_{X}(V)=M_{i}\widehat{\otimes}_{\mathcal{O}_{X}(U_{i})}\mathcal{O}_{X}(V)$ for all rational subsets $V\subseteq U_{i}$. This, together with Corollary~\ref{Terminal objects and quasi-coherence}, shows that $\mathcal{F}$ is a quasi-coherent sheaf of topological Banach modules, generalizing Example \ref{Vector bundles}. As a special case, we see that coherent sheaves (in the sense of \cite{Kedlaya-Liu2}, \S2.3) on locally Noetherian adic spaces over $\Spa(\mathbb{Z}((T)), \mathbb{Z}[[T]])$ form a natural class of examples of quasi-coherent sheaves of topological Banach modules. By \cite{Kedlaya-Liu2}, Theorem~2.5.5, the category of pseudocoherent sheaves on a Tate affinoid adic space $X$ is equivalent via the global sections functor to the category of stably pseudocoherent $\mathcal{O}_{X}(X)$-modules. Thus every affinoid open subspace of a Tate adic space $X$ is $\mathcal{F}$-adapted for every pseudocoherent sheaf $\mathcal{F}$ on the adic space $X$.\end{example}
\begin{example}[Sheaf of convergent power series]\label{Quasi-coherent sheaf of convergent power series}For a Tate adic space $(S, \varpi_{S})$ and a positive integer $n$ define a pre-sheaf $\mathcal{O}_{S}\langle T_1,\dots, T_n\rangle$ of topological Banach algebras on $S$ by setting \begin{equation*}\mathcal{O}_{S}\langle T_1,\dots, T_n\rangle(U)=\mathcal{O}_{S}(U)\langle T_1,\dots, T_n\rangle\end{equation*}for every affinoid open subspace $U$ of $S$, where the convergent power series algebra is equipped with the Gauss norm extended from some norm defining the topology on $\mathcal{O}_{S}(U)$ and making $\varpi_{S}\in\mathcal{O}_{S}(U)$ a norm-multiplicative element. For $U\subseteq S$ open affinoid and $V=U(\frac{f_1,\dots,f_m}{g})$ a rational subset we have\begin{align*}\mathcal{O}_{S}(V)\langle T_1,\dots, T_n\rangle=(\mathcal{O}_{S}(U)\langle X_1,\dots, X_m\rangle/\overline{(f_1-gX_1,\dots, f_{m}-gX_{m})})\langle T_1,\dots, T_n\rangle \\ =(\mathcal{O}_{S}(U)\langle X_1,\dots, X_m\rangle/\overline{(f_1-gX_1,\dots, f_{m}-gX_{m})})\widehat{\otimes}_{\mathcal{O}_{S}(U)}\mathcal{O}_{S}(U)\langle T_1,\dots, T_n\rangle \\=\mathcal{O}_{S}(V)\widehat{\otimes}_{\mathcal{O}_{S}(U)}\mathcal{O}_{S}(U)\langle T_1,\dots, T_n\rangle.\end{align*}By \cite{Bambozzi-Kremnizer20}, Proposition~3.14, $\mathcal{O}_{S}(U)\langle T_1,\dots, T_n\rangle$ is a strictly flat Banach $\mathcal{O}_{S}(U)$-module in the sense that the complete base extension functor \begin{equation*}-\widehat{\otimes}_{\mathcal{O}_{S}(U)}\mathcal{O}_{S}(U)\langle T_1,\dots, T_n\rangle\end{equation*}preserves kernels of all morphisms of Banach $\mathcal{O}_{S}(U)$-modules. In particular, we see that for every finite rational cover $(U_{i})_{i}$ of $U$ the equalizer diagram of (topological) Banach algebras\begin{equation*}0\to \mathcal{O}_{S}(U)\to \bigoplus_{i}\mathcal{O}_{S}(U_{i})\rightrightarrows \bigoplus_{i,j}\mathcal{O}_{S}(U_{i}\cap U_{j})\end{equation*}gives rise to an equalizer diagram of (topological) Banach algebras \begin{equation*}0\to \mathcal{O}_{S}(U)\langle T_1,\dots, T_n\rangle \to \bigoplus_{i}\mathcal{O}_{S}(U_{i})\langle T_1,\dots, T_n\rangle \rightrightarrows \bigoplus_{i,j}\mathcal{O}_{S}(U_{i}\cap U_{j})\langle T_1,\dots, T_n\rangle. \end{equation*}This shows that $\mathcal{O}_{S}\langle T_1,\dots, T_n\rangle$ is a quasi-coherent sheaf of topological Banach $\mathcal{O}_{S}$-algebras on $S$ and every affinoid open subspace of $S$ is $\mathcal{O}_{S}\langle T_1,\dots, T_n\rangle$-adapted. In fact, this example is a special case of Example \ref{Banach bundles}: The proof of \cite{Bambozzi-Kremnizer20}, Proposition~3.14, shows that, for any $U$ as above, $\mathcal{O}_{S}(U)\langle T_1,\dots, T_n\rangle$ is a projective (topological) Banach $\mathcal{O}_{S}(U)$-module, so $\mathcal{O}_{S}\langle T_1,\dots, T_n\rangle$ is an example of a Banach bundle in the sense of Kedlaya.\end{example}
\begin{prop}\label{Sheafy complete Tate rings and delta sheaves}Let $X=\Spa(A, A^{+})$ be a Tate affinoid adic space (with $A=\mathcal{O}_{X}(X)$) and let $B$ be any topological Banach $A$-algebra which is sheafy as a Tate ring. Then $B^{\vartriangle}$ (notation as in Example~\ref{Delta sheaves}) is a globally quasi-coherent sheaf of (topological) Banach algebras on $X$.\end{prop}
\begin{proof}It suffices to prove that the restriction $\widetilde{B}$ of $B^{\vartriangle}$ to the basis $\Rat_{X}$ of $X$ consisting of rational subsets of $X$ is a sheaf of topological Banach algebras. For a rational subset $V=\{\, x\in X\mid \vert f_{i}(x)\vert\leq \vert g(x)\vert\neq 0\,\}$ of $X$ defined by elements $f_1,\dots, f_n, g\in A$ generating the unit ideal of $A$ we have \begin{align*}B\widehat{\otimes}_{A}\mathcal{O}_{X}(V)=B\widehat{\otimes}_{A}A\langle \frac{f_1,\dots, f_n}{g}\rangle \\ =B\widehat{\otimes}_{A}A\langle T_1,\dots, T_n\rangle/\overline{(f_1-gT,\dots, f_n-gT_n)_{A}}=B\langle T_1,\dots, T_n\rangle/\overline{(f_1-gT_1,\dots, f_n-gT_n)_{B}} \\ =B\langle\frac{f_1,\dots, f_n}{g}\rangle,\end{align*}where, by abuse of notation, we identified the elements $f_1,\dots, f_n, g\in A$ with their images in $B$. Hence the assertion follows from the assumption that the Tate ring $B$ is sheafy.\end{proof}  
\begin{rmk}\label{A non-adapted affinoid adic space}By analogy with the situation in algebraic geometry one could try to show that every affinoid Tate adic space $X$ is $\mathcal{F}$-adapted for every quasi-coherent sheaf of (topological) Banach modules $\mathcal{F}$. However, there is the following counterexample due to Scholze and Weinstein (\cite{Berkeley}, Example~9.1.2). Let $K$ be some nonarchimedean field, let $V$ be the quasi-compact open subset of the $2$-dimensional closed unit ball $\mathbb{B}_{K}^{2} = \Spa(K\langle x, y\rangle, K^{\circ}\langle x, y\rangle)$ defined by $V=\{\, \vert x\vert=1 \ \text{or} \ \vert y\vert=1\,\}$ and let $i: V \hookrightarrow \mathbb{B}_{K}^{2}$ be the inclusion map. We consider these as Tate adic spaces by fixing a pseudo-uniformizer $\varpi$ of $K$. We see that $V$ is a non-affinoid open subspace of $\mathbb{B}_{K}^{2}$ covered by the two rational subsets $\{\, \vert x\vert=1\,\}$ and $\{\,\vert y\vert=1\,\}$. The ring of global sections $\Gamma(V, \mathcal{O}_{V})$ is given by \begin{align*} \Gamma(V, \mathcal{O}_{V})& = \ker(\mathcal{O}_{V}(\{\, \vert x\vert=1\,\})\oplus \mathcal{O}_{V}(\{\, \vert y\vert=1\,\}) \to \mathcal{O}_{V}(\{\, \vert x\vert=\vert y\vert=1\,\})) \\ & = \ker(K\langle x,y\rangle\langle x^{-1}\rangle \oplus K\langle x,y\rangle\langle y^{-1}\rangle \to K\langle x,y\rangle\langle x^{-1}, y^{-1}\rangle)=K\langle x,y\rangle. \end{align*}Then the sheaf of (topological) Banach $\mathcal{O}_{\mathbb{B}_{K}^{2}}$-algebras $i_{\ast}\mathcal{O}_{V}$ is quasi-coherent (see Proposition~\ref{Direct images} below) but for the rational subset $U_{0}=\{\, \vert x\vert, \vert y\vert\leq \vert\varpi\vert \,\}$, we have \begin{equation*} (i_{\ast}\mathcal{O}_{V})(U_{0})=0 \neq K\langle x,y\rangle\langle \frac{x,y}{\varpi}\rangle=\mathcal{O}_{\mathbb{B}_{K}^{2}}(U_{0}). \end{equation*}It follows that the affinoid adic space $\mathbb{B}_{K}^{2}$ (considered as an open subset of itself) is not $i_{\ast}\mathcal{O}_{V}$-adapted. For this reason we are always forced to work, given some quasi-coherent (topological) Banach module $\mathcal{F}$ on a Tate adic space $X$, with $\mathcal{F}$-adapted affinoid open covers of $X$ instead of arbitrary affinoid open covers.\end{rmk}

\section{The relative $\Spa$ functor and locally affinoid morphisms}

\subsection{Construction of the relative adic spectrum}

We now construct the relative adic spectrum of a quasi-coherent Banach algebra on an adic space $S$. To this end, we follow the usual construction in the case of schemes (cf. \cite{Stacks}, Tags~01LG, 01LL). We note that the analogous construction in the case of coherent sheaves (i.e., quasi-coherent sheaves of finite type) on locally Noetherian adic spaces is already known and due to Huber (\cite{Huber0}, Proposition~3.12.11, in the case of adic spaces whose rings of local sections have Noetherian rings of definition, see Ludwig's lecture notes \cite{Ludwig20}, Ch.~I, \S8, for the general locally Noetherian case). First we need a suitable notion of sheaf of algebras on an adic space.
\begin{mydef}[Quasi-coherent sheaf of sheafy Banach algebras]\label{Sheaf of sheafy Banach algebras}Let $S$ be a Tate adic space. A quasi-coherent sheaf of topological Banach $\mathcal{O}_{S}$-algebras $\mathcal{A}$ on $S$ is called a quasi-coherent sheaf of (topological) sheafy (resp., stably uniform, resp., strongly rigid-Noetherian, resp., perfectoid, resp., Tate) Banach $\mathcal{O}_{S}$-algebras if there exists a basis $\mathcal{B}$ for the topology on $S$ consisting of $\mathcal{A}$-adapted affinoid open subsets with the property that for every $U\in \mathcal{B}$ the Banach $\mathcal{O}_{S}(U)$-algebra $\mathcal{A}(U)$ is a sheafy (resp., stably uniform, resp., strongly rigid-Noetherian, resp., perfectoid) Huber ring.\end{mydef} 
\begin{mydef}\label{Strongly adapted subsets}Given a quasi-coherent Banach $\mathcal{O}_{S}$-algebra or a quasi-coherent topological Banach $\mathcal{O}_{S}$-algebra, an open affinoid subset $U$ of $S$ is called strongly $\mathcal{A}$-adapted if $U$ is $\mathcal{A}$-adapted and $\mathcal{A}(U)$ is sheafy.\end{mydef}
\begin{mydef}[Quasi-coherent adic Banach $(\mathcal{O}_{S}, \mathcal{O}_{S}^{+})$-algebra]\label{Sheaf of adic Banach algebras}Let $S$ be a Tate adic space. A quasi-coherent sheaf of topological adic Banach algebras on $S$, or, for short, a quasi-coherent (topological) adic Banach $(\mathcal{O}_{S}, \mathcal{O}_{S}^{+})$-algebra, is a pair \begin{equation*}(\mathcal{A}, \mathcal{A}^{+})\end{equation*}where $\mathcal{A}$ is a quasi-coherent sheaf of sheafy (topological) Banach $\mathcal{O}_{S}$-algebras on $S$ and where $\mathcal{A}^{+} \subset\mathcal{A}$ is a strict monomorphism of sheaves of complete topological rings with the following properties: \begin{enumerate}[(1)]\item For every $\mathcal{A}$-adapted affinoid open subspace $U\subseteq S$ the ring $\mathcal{A}^{+}(U)$ is a subring of integral elements of the topological Banach $\mathcal{O}_{S}(U)$-algebra $\mathcal{A}(U)$ and contains the image of $\mathcal{O}_{S}^{+}(U)$ under the structure morphism $\mathcal{O}_{S}(U)\to \mathcal{A}(U)$. \item There exists a basis $\mathcal{B}$ for the topology on $S$ consisting of strongly $\mathcal{A}$-adapted affinoid open subsets such that for every inclusion $V\subset U$, where $U\in \mathcal{B}$ and where $V=U(\frac{f_1, \dots, f_{n}}{g})$ is a rational subset of $U$, the restriction map $\mathcal{A}(U) \to \mathcal{A}(V)$ induces a ring isomorphism between rings of integral elements \begin{equation*} \mathcal{A}(U)\langle \frac{\varphi(f_1), \dots, \varphi(f_n)}{\varphi(g)}\rangle^{+} \cong \mathcal{A}^{+}(V), \end{equation*}where $\varphi: \mathcal{O}_{S}(U) \to \mathcal{A}(U)$ is the map of Banach rings induced by the structure morphism $\mathcal{O}_{S}\to \mathcal{A}$.\end{enumerate}A morphism $(\mathcal{A}, \mathcal{A}^{+})\to (\mathcal{A}', \mathcal{A}'^{+})$ of quasi-coherent sheaves of adic Banach algebras (respectively, of quasi-coherent sheaves of topological adic Banach algebras) is a morphism of quasi-coherent sheaves of (topological) Banach algebras $\mathcal{A}\to \mathcal{A}'$ such that for every $\mathcal{A}$-adapted and $\mathcal{A}'$-adapted affinoid open subspace $U\subset S$ the image of $\mathcal{A}^{+}(U)$ under $\mathcal{A}(U)\to \mathcal{A}'(U)$ is contained in $\mathcal{A}'^{+}(U)$.\end{mydef}
We denote by \begin{equation*}\Ban_{(\mathcal{O}_{S}, \mathcal{O}_{S}^{+})}\end{equation*}the category of quasi-coherent (topological) adic Banach $(\mathcal{O}_{S}, \mathcal{O}_{S}^{+})$-algebras on $S$, respectively; we do not try to define the notion of a (not necessarily quasi-coherent) topological adic Banach $(\mathcal{O}_{S}, \mathcal{O}_{S}^{+})$-algebra.
  
Let us give a name to the affinoid open subsets appearing in condition (2) of the above definition.
\begin{mydef}\label{Adic adapted subsets}Given a quasi-coherent sheaf $(\mathcal{A},\mathcal{A}^{+})$ of adic or topological adic Banach algebras on a Tate adic space $S$, an affinoid open subset $U\subseteq S$ is called $(\mathcal{A},\mathcal{A}^{+})$-adapted if $U$ is strongly $\mathcal{A}$-adapted and for every rational subset $V=U(\frac{f_1, \dots, f_{n}}{g})\subseteq U$ the restriction map $\mathcal{A}(U)\to \mathcal{A}(V)$ induces an isomorphism of rings of integral elements \begin{equation*} \mathcal{A}(U)\langle \frac{\varphi(f_1), \dots, \varphi(f_n)}{\varphi(g)}\rangle^{+} \cong \mathcal{A}^{+}(V), \end{equation*}where $\varphi: \mathcal{O}_{S}(U) \to \mathcal{A}(U)$ is the map induced by the structure morphism $\mathcal{O}_{S} \to \mathcal{A}$.\end{mydef}For a quasi-coherent adic or topological adic Banach algebra $(\mathcal{A},\mathcal{A}^{+})$ on $S$ the family of $(\mathcal{A},\mathcal{A}^{+})$-adapted affinoid open subsets enjoys some properties analogous to the properties of $\mathcal{A}$-adapted affinoid open subsets.
\begin{lemma}\label{Adic adapted covers}Let $S$ be a Tate adic space and let $(\mathcal{A},\mathcal{A}^{+})$ be a quasi-coherent sheaf of (topological) adic Banach algebras on $S$. Every rational subset of an $(\mathcal{A},\mathcal{A}^{+})$-adapted affinoid open subset of $S$ is again $(\mathcal{A},\mathcal{A}^{+})$-adapted.\end{lemma}
\begin{proof}Let $U$ be an $(\mathcal{A},\mathcal{A}^{+})$-adapted affinoid open subset, $V\subseteq U$ a rational subset of $U$ and $V'\subseteq V$ a rational subset of $V$. Then $V'$ is also a rational subset of $U$, so we can write \begin{equation*}V'=U(\frac{f_1,\dots,f_n}{g})\end{equation*}for some $f_1,\dots,f_n, g\in \mathcal{O}_{S}(U)$ generating the unit ideal of $\mathcal{O}_{S}(U)$. By abuse of notation, we identify the elements $f_1,\dots, f_n,g$ with their images in the Banach $\mathcal{O}_{S}(U)$-algebra $\mathcal{A}(U)$ or its rational localizations. Then the restriction map $\mathcal{A}(U)\to \mathcal{A}(V')$ induces an isomorphism \begin{equation*}\mathcal{A}(U)\langle\frac{f_1,\dots,f_n}{g}\rangle^{+}\cong \mathcal{A}^{+}(V').\end{equation*}But $V'$ being contained in $V$ implies that the restriction map $\mathcal{A}(U)\to \mathcal{A}(V')$ factors through $\mathcal{A}(V)$ and therefore the map \begin{equation*}\mathcal{A}(U)\langle\frac{f_1,\dots,f_n}{g}\rangle^{+}\to \mathcal{A}^{+}(V'),\end{equation*}which we know to be an isomorphism of rings, must  factor through $\mathcal{A}(V)\langle\frac{f_1,\dots,f_n}{g}\rangle^{+}$. It follows that \begin{equation*}\mathcal{A}(V)\langle\frac{f_1,\dots,f_n}{g}\rangle^{+} \to \mathcal{A}^{+}(V')\end{equation*}is surjective. On the other hand, the assumption that $U$ is $(\mathcal{A}, \mathcal{A}^{+})$-adapted also implies that the ring map \begin{equation*}\mathcal{A}(U)\langle \frac{f_1,\dots, f_n}{g}\rangle^{+}\to \mathcal{A}(V)\langle \frac{f_1,\dots, f_n}{g}\rangle^{+}\end{equation*}is an isomorphism. Hence we see that \begin{equation*}\mathcal{A}(V)\langle\frac{f_1,\dots, f_n}{g}\rangle^{+}\to \mathcal{A}^{+}(V)\end{equation*}is also injective.\end{proof} 
\begin{lemma}\label{Adic adapted covers 2}Let $S$ be a Tate adic space and let $(\mathcal{A}_{1},\mathcal{A}_{1}^{+}), \dots, (\mathcal{A}_{n},\mathcal{A}_{n}^{+})$ be a finite family of quasi-coherent sheaves of adic Banach $(\mathcal{O}_{S},\mathcal{O}_{S}^{+})$-algebras (respectively, of topological adic Banach $(\mathcal{O}_{S}, \mathcal{O}_{S}^{+})$-algebras). Then the affinoid open subsets of $S$ which are $(\mathcal{A}_{i},\mathcal{A}_{i}^{+})$-adapted for all $i=1,\dots, n$ form a basis for the topology on $S$.\end{lemma}
\begin{proof}By induction, we may assume that this holds for some number $n\geq 1$ and we have to prove that the assertion also holds for $n+1$ quasi-coherent sheaves of (topological) Banach $(\mathcal{O}_{S}, \mathcal{O}_{S}^{+})$-algebras $(\mathcal{A}_{1},\mathcal{A}_{1}^{+}), \dots, (\mathcal{A}_{n},\mathcal{A}_{n}^{+}), (\mathcal{A}_{n+1}, \mathcal{A}_{n+1}^{+})$. By the induction hypothesis, we can cover any open subspace $U\subseteq S$ by a family of affinoid open subsets $(U_{j})_{j}$ such that each $U_{j}$ is $(\mathcal{A}_{i}, \mathcal{A}_{i}^{+})$-adapted for all $i=1,\dots, n$. Cover each $U_{j}$ by $(\mathcal{A}_{n+1}, \mathcal{A}_{n+1}^{+})$-adapted affinoid open subsets $V_{jk}$. Then choose an open cover $(W_{jkl})_{l}$ of each $V_{jk}$ such that for $j, k, l$ the open subset $W_{jkl}$ is a rational subset of both $U_{j}$ and $V_{jk}$. By Lemma~\ref{Adic adapted covers} this implies that each $W_{jkl}$ is $(\mathcal{A}_{i}, \mathcal{A}_{i}^{+})$-adapted for all $i=1,\dots, n+1$.\end{proof}
Let us give a few examples of quasi-coherent adic Banach algebras.
\begin{example}For every Tate adic space $S$ we have a canonical quasi-coherent sheaf of adic Banach algebras given by the pair $(\mathcal{O}_{S}, \mathcal{O}_{S}^{+})$, where $\mathcal{O}_{S}^{+}(U)$ is defined in the usual way for every affinoid open $U\subseteq S$.\end{example}
\begin{example}If $S=\Spa(K, K^{\circ})$, where $K$ is a nonarchimedean field and $S$ is considered as a Tate adic space by choosing some pseudo-uniformizer $\varpi$ of $K$, a quasi-coherent topological adic Banach algebra over $S$ is the same thing as a Huber pair $(A, A^{+})$, where $A$ is a sheafy topological Banach $K$-algebra and $A^{+}$ is a ring of integral elements containing $K^{\circ}$.\end{example}
For the next Proposition~assume that $S$ is a Tate adic space such that for every affinoid open subspace $U$ of $S$ the Huber ring $\mathcal{O}_{S}(U)$ is strongly sheafy in the sense of Hansen and Kedlaya \cite{Hansen-Kedlaya}, Definition~4.1, i.e., both $\mathcal{O}_{S}(U)$ and $\mathcal{O}_{S}(U)\langle T_1,\dots, T_n\rangle$ are sheafy, for all integers $n\geq1$. Note that this assumption is satisfied, for example, if $S$ is locally Noetherian (by \cite{Huber2}, Theorem~2.5), if $S$ is a sousperfectoid space (for example, a perfectoid space) by \cite{Hansen-Kedlaya}, Corollary~7.4, or if $S$ is a diamantine space, by \cite{Hansen-Kedlaya}, Corollary~6.26(a). In fact, the assumption is equivalent to $S$ being covered by affinoid open subspaces that enjoy the above property, as a consequence of \cite{Hansen-Kedlaya}, Corollary~4.5. In this situation, the quasi-coherent sheaf of topological Banach $\mathcal{O}_{S}$-algebras $\mathcal{O}_{S}\langle T_1,\dots, T_n\rangle$ from Example~\ref{Quasi-coherent sheaf of convergent power series}~is a quasi-coherent sheaf of sheafy topological Banach algebras (respectively, of sheafy Banach algebras). Define a subsheaf $\mathcal{O}_{S}^{+}\langle T_1,\dots, T_n\rangle$ of $\mathcal{O}_{S}\langle T_1,\dots, T_n\rangle$ by setting \begin{equation*}\mathcal{O}_{S}^{+}\langle T_1,\dots, T_n\rangle(U)=\mathcal{O}_{S}^{+}(U)\langle T_1,\dots, T_n\rangle\end{equation*}for every affinoid open subset $U\subseteq S$.  
\begin{prop}\label{Quasi-coherent adic Banach algebra of convergent power series}The pair \begin{equation*}(\mathcal{O}_{S}\langle\underline{T}\rangle, \mathcal{O}_{S}^{+}\langle\underline{T}\rangle)=(\mathcal{O}_{S}\langle T_1,\dots, T_n\rangle, \mathcal{O}_{S}^{+}\langle T_1,\dots, T_n\rangle)\end{equation*}is a quasi-coherent topological adic Banach $(\mathcal{O}_{S}, \mathcal{O}_{S}^{+})$-algebra (respectively, a quasi-coherent adic Banach $(\mathcal{O}_{S}, \mathcal{O}_{S}^{+})$-algebra) and every affinoid open subspace $U\subseteq S$ is $(\mathcal{O}_{S}\langle\underline{T}\rangle, \mathcal{O}_{S}^{+}\langle\underline{T}\rangle)$-adapted.\end{prop}
\begin{proof}Let $U\subseteq S$ be an affinoid open and let $V=U(\frac{f_1,\dots, f_m}{g})$ be a rational subset. The subring \begin{equation*}\mathcal{O}_{S}(U)\langle\underline{T}\rangle\langle\frac{f_1,\dots,f_m}{g}\rangle^{+}\subseteq \mathcal{O}_{S}(U)\langle\underline{T}\rangle\langle\frac{f_1,\dots,f_m}{g}\rangle=\mathcal{O}_{S}(V)\langle\underline{T}\rangle\end{equation*}is the integral closure inside $\mathcal{O}_{S}(V)\langle\underline{T}\rangle$ of the open subring $\mathcal{O}_{S}^{+}(U)\langle\underline{T}\rangle\langle\frac{f_1,\dots,f_m}{g}\rangle$. The subring \begin{equation*}\mathcal{O}_{S}^{+}(U)[\frac{f_1,\dots,f_m}{g}][\underline{T}]\subseteq \mathcal{O}_{S}(U)[\frac{1}{g}][\underline{T}]=\mathcal{O}_{S}(U)[\underline{T}][\frac{f_1,\dots,f_m}{g}]\end{equation*}coincides with $\mathcal{O}_{S}^{+}(U)[\underline{T}][\frac{f_1,\dots,f_m}{g}]$, so their topological closures $\mathcal{O}_{S}^{+}(U)\langle\frac{f_1,\dots,f_m}{g}\rangle\langle\underline{T}\rangle$ and $\mathcal{O}_{S}^{+}(U)\langle\underline{T}\rangle\langle\frac{f_1,\dots,f_m}{g}\rangle$ inside $\mathcal{O}_{S}(V)\langle\underline{T}\rangle=\mathcal{O}_{S}(U)\langle\frac{f_1,\dots,f_m}{g}\rangle\langle\underline{T}\rangle=\mathcal{O}_{S}(U)\langle\underline{T}\rangle\langle\frac{f_1,\dots,f_m}{g}\rangle$ are equal. In particular, the integral closures inside $\mathcal{O}_{S}(V)\langle\underline{T}\rangle$ of those topological closures are also equal. Hence to prove the desired equality \begin{equation*}\mathcal{O}_{S}(U)\langle\underline{T}\rangle\langle\frac{f_1,\dots,f_m}{g}\rangle^{+}=\mathcal{O}_{S}^{+}(V)\langle\underline{T}\rangle\end{equation*}it suffices to show that the right hand side is the integral closure of $\mathcal{O}_{S}^{+}(U)\langle\frac{f_1,\dots,f_m}{g}\rangle\langle\underline{T}\rangle$ inside $\mathcal{O}_{S}(V)\langle\underline{T}\rangle=\mathcal{O}_{S}(U)\langle\frac{f_1,\dots,f_m}{g}\rangle\langle\underline{T}\rangle$. But since $\mathcal{O}_{S}^{+}(V)$ is by Definition~the integral closure of $\mathcal{O}_{S}^{+}(U)\langle\frac{f_1,\dots,f_m}{g}\rangle$ inside $\mathcal{O}_{S}(V)$, the polynomial ring $\mathcal{O}_{S}^{+}(V)[\underline{T}]$ is the integral closure of $\mathcal{O}_{S}^{+}(U)\langle\frac{f_1,\dots,f_m}{g}\rangle[\underline{T}]$ inside $\mathcal{O}_{S}(V)[\underline{T}]$. It follows from this that $\mathcal{O}_{S}^{+}(V)\langle\underline{T}\rangle$ is the integral closure of the open subring $\mathcal{O}_{S}^{+}(U)\langle\frac{f_1,\dots,f_m}{g}\rangle\langle\underline{T}\rangle$ inside $\mathcal{O}_{S}(V)\langle\underline{T}\rangle$ (see, for example, \cite{Huber0}, Lemma~2.4.3, or \cite{Huebner21}, Lemma~11.4).\end{proof}
\begin{mydef}\label{Definition of uniform quasi-coherent adic Banach algebra}For a Tate adic space $S$, a quasi-coherent adic Banach algebra or topological adic Banach algebra $(\mathcal{A}, \mathcal{A}^{+})$ on $S$ is called uniform (respectively, strongly rigid-Noetherian, resp., (sous)perfectoid) if $S$ admits a basis $\mathcal{B}$ consisting of $(\mathcal{A}, \mathcal{A}^{+})$-adapted subsets such that $\mathcal{A}(U)$ is stably uniform (respectively, strongly rigid-Noetherian, resp., (sous)perfectoid) for all $U\in\mathcal{B}$.\end{mydef}
\begin{lemma}\label{Morphisms of quasi-coherent adic Banach algebras}Let $S$ be a Tate affinoid adic space. Let $(\mathcal{A}, \mathcal{A}^{+})$ and $(\mathcal{A}', \mathcal{A}'^{+})$ be two quasi-coherent topological adic Banach $(\mathcal{O}_{S}, \mathcal{O}_{S}^{+})$-algebras such that $S$ is both $(\mathcal{A},\mathcal{A}^{+})$-adapted and $(\mathcal{A}', \mathcal{A}'^{+})$-adapted. Then the natural map \begin{equation*} \Hom_{\Ban_{(\mathcal{O}_{S}, \mathcal{O}_{S}^{+})}}((\mathcal{A}, \mathcal{A}^{+}), (\mathcal{A}', \mathcal{A}'^{+})) \to \Hom_{(\mathcal{O}_{S}(S), \mathcal{O}_{S}^{+}(S))}((\mathcal{A}(S), \mathcal{A}^{+}(S)), (\mathcal{A}'(S), \mathcal{A}'^{+}(S)),\end{equation*}where the set on right hand side is the set of morphisms of Huber pairs over $(\mathcal{O}_{S}(S), \mathcal{O}^{+}(S))$, is a bijection.\end{lemma}
\begin{proof} Let $V=S(\frac{f_1,\dots, f_n}{g}) \subseteq S$ be a rational subset. Denote by $\varphi: \mathcal{O}_{S}(S) \to \mathcal{A}(S)$ and $\varphi': \mathcal{O}_{S}(S) \to \mathcal{A}'(S)$ the structure maps. Given a morphism of Huber pairs $(\mathcal{A}(S), \mathcal{A}^{+}(S)) \to (\mathcal{A}'(S), \mathcal{A}'^{+}(S))$ over $(\mathcal{O}_{S}(S), \mathcal{O}_{S}^{+}(S))$, we can canonically extend it to a morphism of the rational localizations \begin{align*}(\mathcal{A}(S)\langle\frac{\varphi(f_1),\dots, \varphi(f_n)}{\varphi(g)}\rangle, \mathcal{A}(S)\langle\frac{\varphi(f_1),\dots, \varphi(f_n)}{\varphi(g)}\rangle^{+}) \\ \to (\mathcal{A}'(S)\langle\frac{\varphi'(f_1),\dots, \varphi'(f_n)}{\varphi'(g)}\rangle, \mathcal{A}'(S)\langle\frac{\varphi'(f_1),\dots, \varphi'(f_n)}{\varphi'(g)}\rangle^{+}).\end{align*}Composing this with the canonical isomorphisms of Huber pairs supplied by the definition~of $(\mathcal{A},\mathcal{A}^{+})$-adapted subsets yields a morphism of Huber pairs \begin{equation*}(\mathcal{A}(V), \mathcal{A}^{+}(V)) \to (\mathcal{A}'(V), \mathcal{A}'^{+}(V)).\end{equation*}Repeating this construction for all rational subsets $V$ of $S$ gives us a unique morphism $(\mathcal{A}, \mathcal{A}^{+}) \to (\mathcal{A}', \mathcal{A}'^{+})$ of quasi-coherent topological adic Banach algebras which recovers the given morphism $(\mathcal{A}(S), \mathcal{A}^{+}(S)) \to (\mathcal{A}'(S), \mathcal{A}'^{+}(S))$ upon taking global sections.\end{proof}   
\begin{lemma}\label{Uniform quasi-coherent algebras}Let $S$ be a Tate adic space and let $(\mathcal{A}, \mathcal{A}^{+})$ be a quasi-coherent sheaf of adic Banach algebras over $S$. Let $V\subset U$ be an inclusion of $(\mathcal{A},\mathcal{A}^{+})$-adapted affinoid open subspaces of $S$. Then the diagram of morphisms of affinoid adic spaces \begin{center} \begin{tikzcd} \Spa(\mathcal{A}(V), \mathcal{A}^{+}(V)) \arrow{r} \arrow{d}{f_{V}} & \Spa(\mathcal{A}(U), \mathcal{A}^{+}(U)) \arrow{d}{f_{U}} \\ V \arrow[hook]{r} & U \end{tikzcd} \end{center}gives rise to an isomorphism \begin{equation*} \rho_{V}^{U}: \Spa(\mathcal{A}(V), \mathcal{A}^{+}(V)) \to f_{U}^{-1}(V).\end{equation*}\end{lemma}
\begin{proof} Suppose first that $V$ is a rational subset of $U$. Write $V = U(\frac{f_1, \dots, f_n}{g})$ for some elements $f_1,\dots, f_n, g \in \mathcal{O}_{S}(U)$ generating the unit ideal of $\mathcal{O}_{S}(U)$. By the assumption that $U$ is $(\mathcal{A},\mathcal{A}^{+})$-adapted, the restriction map $\mathcal{A}(U)\to\mathcal{A}(V)$ induces an isomorphism of Huber pairs\begin{equation*} (\mathcal{A}(U)\langle\frac{\varphi(f_1), \dots,\varphi(f_n)}{\varphi(g)}\rangle, \mathcal{A}(U)\langle\frac{\varphi(f_1),\dots,\varphi(f_n)}{\varphi(g)}\rangle^{+}) \cong (\mathcal{A}(V), \mathcal{A}^{+}(V)). \end{equation*}This induces the desired isomorphism of adic spaces \begin{equation*} \rho_{V}^{U}: \Spa(\mathcal{A}(V), \mathcal{A}^{+}(V)) \to \Spa(\mathcal{A}(U), \mathcal{A}^{+}(U))(\frac{\varphi(f_1),\dots, \varphi(f_n)}{\varphi(g)}) = f_{U}^{-1}(V) \end{equation*}in the case that $V$ is rational. 

In the general case, cover $V$ by rational subsets $V_{i}$ of $U$. Then each $V_{i}$ is also a rational subset of $V$. We have to show that the image of $\Spa(\mathcal{A}(V), \mathcal{A}^{+}(V))$ in $\Spa(\mathcal{A}(U), \mathcal{A}^{+}(U))$ coincides with the pre-image of $V$ under $f_{U}$ and that the morphism $\Spa(\mathcal{A}(V), \mathcal{A}^{+}(V)) \to \Spa(\mathcal{A}(U), \mathcal{A}^{+}(U))$ is an isomorphism onto its image. But this follows from the analogous assertions for the maps $\Spa(\mathcal{A}(V_{i}), \mathcal{A}^{+}(V_{i})) \to \Spa(\mathcal{A}(U), \mathcal{A}^{+}(U))$ and $\Spa(\mathcal{A}(V_{i}), \mathcal{A}^{+}(V_{i})) \to \Spa(\mathcal{A}(V), \mathcal{A}^{+}(V))$.\end{proof}   
\begin{prop}\label{Relative adic spectrum}Let $(S, \varpi_{S})$ be either a Tate adic space. Let $(\mathcal{A}, \mathcal{A}^{+})$ be a quasi-coherent sheaf of adic Banach algebras on $S$. There exists a morphism of adic spaces\begin{equation*} \pi: \underline{\Spa}_{S}(\mathcal{A}, \mathcal{A}^{+}) \to S \end{equation*}with the following properties \begin{enumerate}[(1)]\item For every open affinoid adic subspace $U \subset S$ which is $(\mathcal{A}, \mathcal{A}^{+})$-adapted there exists an isomorphism of adic spaces $i_{U}: \pi^{-1}(U) \to \Spa(\mathcal{A}(U), \mathcal{A}^{+}(U))$. \item For $V \subset U\subset S$, where $U$ and $V$ are $(\mathcal{A}, \mathcal{A}^{+})$-adapted affinoid open subspaces of $S$, the composition \begin{center} \begin{tikzcd} \Spa(\mathcal{A}(V), \mathcal{A}^{+}(V)) \arrow{r}{i_{V}^{-1}} & \pi^{-1}(V) \arrow[hook]{r} & \pi^{-1}(U) \arrow{r}{i_{U}} & \Spa(\mathcal{A}(U), \mathcal{A}^{+}(U)) \end{tikzcd} \end{center}is the map from Lemma~\ref{Uniform quasi-coherent algebras}.\end{enumerate}If the quasi-coherent adic or topological adic Banach algebra $(\mathcal{A}, \mathcal{A}^{+})$ is strongly rigid-Noetherian (respectively, uniform, respectively, (sous)perfectoid), then the adic space $\underline{\Spa}_{S}(\mathcal{A}, \mathcal{A}^{+})$ is locally Noetherian (respectively, uniform, respectively, (sous)perfectoid).\end{prop}
\begin{proof}By the definition~of a quasi-coherent sheaf of sheafy Banach algebras (respectively, of a quasi-coherent sheaf of sheafy topological Banach algebras), the $(\mathcal{A},\mathcal{A}^{+})$-adapted affinoid open subsets of $S$ form a basis for the topology on $S$. Hence the main part of the proposition~follows from Lemma~\ref{Uniform quasi-coherent algebras} and \cite{Dine25-2}, Lemma~3.7, by choosing the basis $\mathcal{B}$ in loc.~cit.~to be the family of all $(\mathcal{A}, \mathcal{A}^{+})$-adapted affinoid open subspaces of $S$. 

It remains to prove the additional statement concerning the properties of the adic space $\underline{\Spa}_{S}(\mathcal{A}, \mathcal{A}^{+})$. Note that, by construction, the adic space $\underline{\Spa}_{S}(\mathcal{A}, \mathcal{A}^{+})$ is covered by the open subspaces \begin{equation*}\pi^{-1}(U)\cong\Spa(\mathcal{A}(U), \mathcal{A}^{+}(U)),\end{equation*}where $U$ ranges over the $(\mathcal{A}, \mathcal{A}^{+})$-adapted affinoid open subspaces of $S$. Thus $\underline{\Spa}_{S}(\mathcal{A}, \mathcal{A}^{+})$ is uniform (resp., locally Noetherian, resp., (sous)perfectoid) if and only if $\pi^{-1}(U)$, for every $(\mathcal{A}, \mathcal{A}^{+})$-adapted affinoid open subspace $U$ of $S$, can be covered by affinoid open subspaces of the form $\Spa(A, A^{+})$, where $A$ is a stably uniform Tate ring (resp., a strongly rigid-Noetherian Huber ring, resp., a (sous)perfectoid Tate ring). Now, if the quasi-coherent adic Banach $(\mathcal{O}_{S}, \mathcal{O}_{S}^{+})$-algebra $(\mathcal{A}, \mathcal{A}^{+})$ is strongly rigid-Noetherian (respectively, uniform, respectively, (sous)perfectoid), this observation, together with Lemma~\ref{Adic adapted covers}, allows us to choose a basis $\mathcal{B}'$ for the topology on $S$ consisting of $(\mathcal{A}, \mathcal{A}^{+})$-adapted affinoid open subspaces $V$ of $S$ with the property that $\mathcal{A}(V)$ is strongly rigid-Noetherian (respectively, a stably uniform Tate ring, respectively, a (sous)perfectoid Tate ring). For any $(\mathcal{A}, \mathcal{A}^{+})$-adapted affinoid open subspace $U\subseteq S$, choose a cover $(V_{j})_{j}$ of $U$ with $V_{j}\in \mathcal{B}'$. Then $(\pi^{-1}(V_{j}))_{j}$ is an open cover of $\pi^{-1}(U)$ and $\pi^{-1}(V_{j})\cong \Spa(\mathcal{A}(V_{j}), \mathcal{A}^{+}(V_{j}))$ for all $j$. This shows that $\underline{\Spa}_{S}(\mathcal{A}, \mathcal{A}^{+})$ is a locally Noetherian (resp.~uniform, resp.~(sous)perfectoid) adic space.\end{proof}
\begin{mydef}We call the adic space over $S$ \begin{equation*}\underline{\Spa}_{S}(\mathcal{A}, \mathcal{A}^{+})\end{equation*}from the above Proposition~the \textit{relative adic spectrum}~of the quasi-coherent adic Banach algebra or quasi-coherent topological adic Banach algebra $(\mathcal{A}, \mathcal{A}^{+})$.\end{mydef}
\begin{example}[Relative closed unit ball]For $S$ a Tate adic space all of whose affinoid open subspaces are adic spectra of strongly sheafy Huber pairs (see our discussion preceding Proposition~\ref{Quasi-coherent adic Banach algebra of convergent power series}) we define, using Proposition~\ref{Quasi-coherent adic Banach algebra of convergent power series}, \begin{equation*}\mathbb{B}_{S}^{n}=\underline{\Spa}_{S}(\mathcal{O}_{S}\langle T_1,\dots,T_n\rangle, \mathcal{O}_{S}^{+}\langle T_1,\dots, T_n\rangle)\end{equation*}and call it the $n$-dimensional (relative) closed unit ball over $S$. Note that this Definition~in particular makes sense if $S$ is locally Noetherian, sousperfectoid or diamantine.\end{example}
\begin{rmk}Unlike in algebraic geometry, the relative adic spectrum $\underline{\Spa}_{S}(\mathcal{A}, \mathcal{A}^{+})$ need not be an affinoid adic space even if $S$ is affinoid. However, $\underline{\Spa}_{S}(\mathcal{A}, \mathcal{A}^{+})$ is affinoid by construction if $S$ is affinoid and $(\mathcal{A},\mathcal{A}^{+})$-adapted.\end{rmk}

\subsection{Characterization of locally affinoid morphisms}
 
In this subsection, we study locally affinoid morphisms of Tate adic spaces, which can be thought of as analogs of affine morphisms of schemes. We want to use the relative $\Spa$ functor from last subsection to describe these morphisms in the same way the relative $\Spec$ functor in algebraic geometry is used to describe affine morphisms of schemes. To this end, we need to study direct images of quasi-coherent sheaves of adic Banach algebras. 
\begin{lemma}\label{Direct images, affinoid case}Let $f: X \to S$ be a morphism of affinoid Tate adic spaces. Let $\mathcal{F}$ be a globally quasi-coherent sheaf of topological Banach $\mathcal{O}_{X}$-modules. Then \begin{equation*}f_{\ast}\mathcal{F}=\widetilde{\mathcal{F}(X)_{\mathcal{O}_{S}(S)}}, \end{equation*}where $\mathcal{F}(X)_{\mathcal{O}_{S}(S)}$ is the restriction of scalars of the Banach $\mathcal{O}_{X}(X)$-module $\mathcal{F}(X)$. In other words, the direct image $f_{\ast}\mathcal{F}$ is a globally quasi-coherent sheaf of (topological) Banach $\mathcal{O}_{S}$-modules. The analogous assertion is also true for globally quasi-coherent sheaves of (topological) Banach $\mathcal{O}_{X}$-algebras.\end{lemma}
\begin{proof} We have to show that $\mathcal{F}(f^{-1}(V)) = \mathcal{F}(X)\widehat{\otimes}_{\mathcal{O}_{S}(S)}\mathcal{O}_{S}(V)$ for every rational subset $V\subseteq S$.

To this end, let $V\subseteq S$ be an arbitrary rational subset and write $V$ in the form \begin{equation*}V=\{\, x \in S \mid \vert f_{i}(x)\vert\leq \vert g(x)\vert \neq 0 \ \text{for all} \ i=1, \dots, n\,\} \end{equation*}where $f_1, \dots, f_{n}, g \in \mathcal{O}_{S}(S)$ are some elements generating the unit ideal of $\mathcal{O}_{S}(S)$. Let $\varphi: \mathcal{O}_{S}(S) \to \mathcal{O}_{X}(X)$ be the homomorphism of Tate rings corresponding to $f$. Then the images of $f_1, \dots, f_n, g$ under $\varphi$ generate the unit ideal of $\mathcal{O}_{X}(X)$ and $f^{-1}(V)$ is the rational subset of $X$ given by \begin{equation*}f^{-1}(V) = \{\, x \in X \mid \vert\varphi(f_{i})(x)\vert\leq \vert\varphi(g)(x)\vert \neq 0 \ \text{for all} i=1, \dots, n\,\}.\end{equation*}Consequently, we have $\mathcal{O}_{X}(f^{-1}(V))=\mathcal{O}_{X}(X)\langle\frac{\varphi(f_1),\dots, \varphi(f_n)}{\varphi(g)}\rangle$ and $\mathcal{O}_{S}(V)=\mathcal{O}_{S}(S)\langle\frac{f_1, \dots, f_n}{g}\rangle$, and thus \begin{equation*} \mathcal{O}_{X}(f^{-1}(V))=\mathcal{O}_{X}(X)\widehat{\otimes}_{\mathcal{O}_{S}(S)}\mathcal{O}_{S}(V).\end{equation*}Since $X$ is $\mathcal{F}$-adapted and $f^{-1}(V)$ is rational, we have \begin{equation*}\mathcal{F}(f^{-1}(V))=\mathcal{F}(X)\widehat{\otimes}_{\mathcal{O}_{X}(X)}\mathcal{O}_{X}(f^{-1}(V))=\mathcal{F}(X)\widehat{\otimes}_{\mathcal{O}_{X}(X)}(\mathcal{O}_{X}(X)\widehat{\otimes}_{\mathcal{O}_{S}(S)}\mathcal{O}_{S}(V)).\end{equation*}It follows that $\mathcal{F}(f^{-1}(V))=\mathcal{F}(X)\widehat{\otimes}_{\mathcal{O}_{S}(S)}\mathcal{O}_{S}(V)$, as desired.\end{proof}
\begin{lemma}\label{Direct images of adic Banach algebras, affinoid case}Let $f: X\to S$ be a morphism of affinoid Tate adic spaces or a morphism of affinoid Tate adic spaces. Let $(\mathcal{A}, \mathcal{A}^{+})$ be a quasi-coherent sheaf of adic Banach $(\mathcal{O}_{X}, \mathcal{O}_{X}^{+})$-algebras such that $X$ is $(\mathcal{A}, \mathcal{A}^{+})$-adapted. Then the pair $(f_{\ast}\mathcal{A}, f_{\ast}\mathcal{A}^{+})$ is a quasi-coherent sheaf of adic Banach $(\mathcal{O}_{S}, \mathcal{O}_{S}^{+})$-algebras and $S$ is $(f_{\ast}\mathcal{A}, f_{\ast}\mathcal{A}^{+})$-adapted.\end{lemma}
\begin{proof}By Lemma~\ref{Direct images, affinoid case}~we know that $S$ is $f_{\ast}\mathcal{A}$-adapted. Let $U$ be a rational subset of $S$ and let $V=U(\frac{f_1,\dots,f_{k}}{g})$ be a rational subset of $U$ defined by elements $f_1,\dots, f_{k}, g\in \mathcal{O}_{S}(U)$ generating the unit ideal in $\mathcal{O}_{S}(U)$. As in the proof of Lemma~\ref{Direct images, affinoid case}, $f^{-1}(U)$ is a rational subset of $X$ and $f^{-1}(V)=U(\frac{\varphi(f_1),\dots,\varphi(f_{k})}{\varphi(g)})$ where $\varphi$ is the map $\mathcal{O}_{S}(U)\to\mathcal{O}_{X}(f^{-1}(U))$. In particular, $\mathcal{A}(f^{-1}(U))$ is sheafy, $X$ being strongly $\mathcal{A}$-adapted. Let $\varphi'$ be the composite map $\mathcal{O}_{S}(U)\to\mathcal{O}_{X}(f^{-1}(U))\to\mathcal{A}(f^{-1}(U))$. Using that $X$ (and hence, by Lemma~\ref{Adic adapted covers}, also its rational subset $f^{-1}(U)$) is $(\mathcal{A}, \mathcal{A}^{+})$-adapted we see that \begin{align*}(f_{\ast}\mathcal{A}^{+})(V)=\mathcal{A}^{+}(f^{-1}(V))=\mathcal{A}(f^{-1}(U))\langle\frac{\varphi'(f_1),\dots,\varphi'(f_k)}{\varphi'(g)}\rangle^{+}\\ =(f_{\ast}\mathcal{A})(U)\langle\frac{\varphi'(f_1),\dots,\varphi'(f_k)}{\varphi'(g)}\rangle^{+}.\end{align*}Choosing $U=S$, this shows that $(f_{\ast}\mathcal{A}, f_{\ast}\mathcal{A}^{+})$ is a quasi-coherent sheaf of topological adic Banach $(\mathcal{O}_{S}, \mathcal{O}_{S}^{+})$-algebras (respectively, a quasi-coherent sheaf of adic Banach $(\mathcal{O}_{S}, \mathcal{O}_{S}^{+})$-algebras) and that $S$ is $(f_{\ast}\mathcal{A}, f_{\ast}\mathcal{A}^{+})$-adapted, as claimed.\end{proof} 
Recall that a morphism of adic spaces $f: X\to S$ is called quasi-compact if the pre-image $f^{-1}(U)$ of any quasi-compact open subspace $U$ of $S$ is quasi-compact. The following Lemma~is probably well-known, as it is a direct analog of \cite{Stacks}, Tags 01K4, with completely analogous proof. We include the Lemma~here since we could not locate an explicit reference in the literature.
\begin{lemma}\label{Quasi-compact morphisms} Let $f: X \to S$ be a morphism of adic spaces. The following are equivalent: \begin{enumerate}[(1)]\item $f$ is quasi-compact. \item The pre-image under $f$ of every affinoid open subspace of $S$ is quasi-compact. \item There exists some open cover $(U_{i})_{i}$ of $S$ by affinoid subspaces such that $f^{-1}(U_{i})$ is quasi-compact for all $i$.\end{enumerate}\end{lemma}
\begin{proof}The proof of \cite{Stacks}, Tag~01K4, applies verbatim with affinoid open subspaces and rational subsets playing the role of affine open subsets and standard open subsets, respectively. \end{proof}
A particular class of quasi-compact morphisms of adic spaces is given by the following definition.
\begin{mydef}\label{Locally affinoid morphisms}A morphism of adic spaces $f: X \to S$ is called locally affinoid if there exists an open cover $(U_{i})_{i}$ of $S$ by affinoid open subspaces such that $f^{-1}(U_{i})$ is affinoid for all indices $i$.\end{mydef}In the above definition we also view the empty set as an affinoid adic space. We can now prove that, at least in some important special cases, direct images of quasi-coherent sheaves on Tate adic spaces behave similarly to direct images of quasi-coherent sheaves on schemes.
\begin{prop}\label{Direct images}Let $f: X\to S$ be a locally affinoid morphism of Tate adic spaces. Suppose that $\mathcal{F}$ is a quasi-coherent sheaf of topological Banach $\mathcal{O}_{X}$-modules or topological Banach $\mathcal{O}_{X}$-algebras such that for some affinoid open cover $(U_{i})_{i}$ of $S$ the pre-images $f^{-1}(U_{i})$ are affinoid and $\mathcal{F}$-adapted. Then the sheaf $f_{\ast}\mathcal{F}$ on $S$ is a quasi-coherent sheaf of topological Banach $\mathcal{O}_{S}$-modules or Banach $\mathcal{O}_{S}$-algebras and all $U_{i}$ are $f_{\ast}\mathcal{F}$-adapted.\end{prop}
\begin{proof}By Corollary~\ref{Terminal objects and quasi-coherence} it suffices to show that all $U_{i}$ are $f_{\ast}\mathcal{F}$-adapted. Let $V$ be a rational subset of some $U_{i}$. Then, as in the proof of Lemma~\ref{Direct images, affinoid case}, we see that $f^{-1}(V)$ is a rational subset of the affinoid adic space $f^{-1}(U_{i})$ and thus, $f^{-1}(U_{i})$ being $\mathcal{F}$-adapted, we have \begin{align*}f_{\ast}\mathcal{F}(V) =\mathcal{F}(f^{-1}(V))=\mathcal{F}(f^{-1}(U_{i}))\widehat{\otimes}_{\mathcal{O}_{X}(f^{-1}(U_{i}))}\mathcal{O}_{X}(f^{-1}(V))\\=\mathcal{F}(f^{-1}(U_{i}))\widehat{\otimes}_{\mathcal{O}_{X}(f^{-1}(U_{i}))}(\mathcal{O}_{X}(f^{-1}(U_{i}))\widehat{\otimes}_{\mathcal{O}_{S}(U_{i})}\mathcal{O}_{S}(V))&=f_{\ast}\mathcal{F}(U_{i})\widehat{\otimes}_{\mathcal{O}_{S}(U_{i})}\mathcal{O}_{S}(V).\end{align*}This proves the proposition.\end{proof}
\begin{example}The simplest example in which the above proposition applies is that of the structure sheaf $\mathcal{F}=\mathcal{O}_{X}$. In this case the proposition tells us that the direct image $f_{\ast}\mathcal{O}_{X}$ is a quasi-coherent sheaf of topological Banach algebras on $S$ (respectively, a quasi-coherent sheaf of Banach algebras on $S$) for every locally affinoid adic morphism $f: X\to S$.\end{example}
\begin{example}\label{Direct images of pseudo-coherent sheaves}Let $\mathcal{F}$ be a Banach bundle or a pseudo-coherent sheaf (for example, a vector bundle) on a Tate adic space $X$ and let $f: X\to S$ be a locally affinoid morphism to another Tate adic space. Since all affinoid open subspaces of $X$ are $\mathcal{F}$-adapted, Proposition~\ref{Direct images}~tells us that the direct image $f_{\ast}\mathcal{F}$ is a quasi-coherent sheaf of topological Banach modules on $S$. Note that this example includes the case when $X$ is locally Noetherian and $\mathcal{F}$ is a coherent sheaf on $X$.\end{example}
\begin{rmk}\label{Locally affinoid morphisms and kernels}It is clear from Lemma~\ref{Quasi-compact morphisms}~that every locally affinoid morphism is quasi-compact. We caution the reader that, unlike in algebraic geometry, where the property of a morphism of schemes being affine is local for the Zariski topology and even for the fpqc topology, a morphism $f: X \to S$ of adic spaces can be locally affinoid but not ``globally affinoid'', i.e., there might exist affinoid open subspaces $U$ of $S$ such that $f^{-1}(U)$ is not affinoid. Indeed, for $K$ a nonarchimedean field, Scholze and Weinstein showed in \cite{Berkeley}, Example~9.1.2, that for the non-affinoid, quasi-compact open subset \begin{equation*}V=\{\,\vert x\vert=1\,\}\cup \,\{\vert y\vert=1\,\}\end{equation*}of the two-dimensional closed unit ball $\mathbb{B}_{K}^{2}=\Spa(K\langle x, y \rangle, K^{\circ}\langle x, y\rangle)$ the inclusion map $i: V \hookrightarrow \mathbb{B}_{K}^{2}$ is locally affinoid.

In this context, we can also show that the inclusion functor from the category of globally quasi-coherent sheaves of Banach modules into the category of all quasi-coherent sheaves of Banach modules can fail to preserve kernels. Indeed, let $V_{1}=\{\, \vert x\vert=1\,\}$, $V_{2}=\{\, \vert y\vert=1\,\}\subseteq \mathbb{B}_{K}^{2}$ and let $V=V_{1}\cup V_{2}$, as before. We view these as Tate adic spaces by fixing a pseudo-uniformizer $\varpi\in K^{\circ}$ of $K$. Let $i_{1}: V_{1}\hookrightarrow \mathbb{B}_{K}^{2}$, $i_{2}: V_{2}\hookrightarrow \mathbb{B}_{K}^{2}$ and $i_{12}: V_{1}\cap V_{2}\hookrightarrow \mathbb{B}_{K}^{2}$ be the inclusion maps. Since $V_{1}$, $V_{2}$ and $V_{1}\cap V_{2}$ are affinoid adic spaces (they are rational subsets of $\mathbb{B}_{K}^{2}$), the direct images $i_{1\ast}\mathcal{O}_{V_{1}}$, $i_{2\ast}\mathcal{O}_{V_{2}}$ and $i_{12\ast}\mathcal{O}_{V_{1}\cap V_{2}}$ are globally quasi-coherent sheaves of Banach modules, by Lemma~\ref{Direct images, affinoid case}. It is also readily seen that $i_{\ast}\mathcal{O}_{V}$ is the kernel, in the category of sheaves of topological Banach modules on $\mathbb{B}_{K}^{2}$, of the natural morphism \begin{equation*}i_{1\ast}\mathcal{O}_{V_{1}}\oplus i_{2\ast}\mathcal{O}_{V_{2}}\to i_{12\ast}\mathcal{O}_{V_{1}\cap V_{2}}.\end{equation*}Since $i: V\hookrightarrow \mathbb{B}_{K}^{2}$ is locally affinoid, the direct image $i_{\ast}\mathcal{O}_{V}$ is a quasi-coherent sheaf of topological Banach modules, by Proposition~\ref{Direct images}. Hence it is the kernel of the above morphism in the category of quasi-coherent sheaves of topological Banach modules on $\mathbb{B}_{K}^{2}$. But it is readily seen that the kernel of $i_{1\ast}\mathcal{O}_{V_{1}}\oplus i_{2\ast}\mathcal{O}_{V_{2}}\to i_{12\ast}\mathcal{O}_{V_{1}\cap V_{2}}$ in the category of globally quasi-coherent pre-sheaves (and hence also in the category of globally quasi-coherent sheaves) is given by \begin{equation*}\widetilde{\Gamma(\mathbb{B}_{K}^{2}, i_{\ast}\mathcal{O}_{V})}=\widetilde{\Gamma(V, \mathcal{O}_{V})}=\widetilde{K\langle x, y\rangle}=\mathcal{O}_{\mathbb{B}_{K}^{2}}.\end{equation*}In particular, the canonical morphism \begin{equation*}\widetilde{\Gamma(\mathbb{B}_{K}^{2}, i_{\ast}\mathcal{O}_{V})}\to i_{\ast}\mathcal{O}_{V}\end{equation*}becomes the zero morphism when restricted to the rational subset $\{\, \vert x\vert, \vert y\vert\leq \vert \varpi\vert\,\}$.\end{rmk}
We emphasize~that in the situation of the above remark the affinoid adic space $\mathbb{B}_{K}^{2}$ is not $i_{\ast}\mathcal{O}_{V}$-adapted. The next Lemma~shows that this is not a coincidence. 
\begin{lemma}\label{Locally affinoid}Let $f: X\to S$ be a quasi-compact morphism of Tate adic spaces and let $U\subseteq S$ be an affinoid open subspace. If the pre-image $f^{-1}(U)$ is also affinoid, then $U$ is $f_{\ast}\mathcal{O}_{X}$-adapted.\end{lemma}
\begin{proof}Since $f$ is quasi-compact, the direct image $f_{\ast}\mathcal{O}_{X}$ is a sheaf of topological Banach modules on $S$. The restriction $f\vert_{f^{-1}(U)}: f^{-1}(U) \to U$ is a morphism of affinoid Tate adic spaces, so the assertion follows from Lemma~\ref{Direct images, affinoid case}.\end{proof}
\begin{lemma}\label{Locally affinoid 2}For a quasi-compact morphism $f: X\to S$ of Tate adic spaces the following are equivalent:\begin{enumerate}[(1)]\item $f$ is locally affinoid.\item There exists an open cover $(U_{i})_{i}$ of $S$ by $f_{\ast}\mathcal{O}_{X}$-adapted affinoid subspaces such that $f^{-1}(U_{i})$ is affinoid for each $i$.\end{enumerate}\end{lemma}
\begin{proof}This is immediate from the above lemma.\end{proof}
With the notion of locally affinoid morphisms in place, we have the following variant of Proposition~\ref{Direct images}~for quasi-coherent adic Banach algebras.
\begin{prop}\label{Direct images 2}Let $f: X\to S$ be a locally affinoid morphism of Tate adic spaces. The pair $(f_{\ast}\mathcal{O}_{X}, f_{\ast}\mathcal{O}_{X}^{+})$ is a quasi-coherent sheaf of topological adic Banach $(\mathcal{O}_{S}, \mathcal{O}_{S}^{+})$-algebras on $S$.\end{prop}
\begin{proof}By Proposition~\ref{Direct images}~we know that $f_{\ast}\mathcal{O}_{X}$ is a quasi-coherent sheaf of topological Banach algebras (respectively, of Banach algebras) on $S$. Since $f$ is locally affinoid, Lemma~\ref{Locally affinoid 2}~shows that there exists a basis $\mathcal{B}$ for the topology on $S$ consisting of $f_{\ast}\mathcal{O}_{X}$-adapted affinoid open subsets such that the pre-image $f^{-1}(U)$ of every $U\in \mathcal{B}$ is an affinoid adic space (indeed, let $(U_{i})_{i}$ be an affinoid open cover of $S$ as in Lemma~\ref{Locally affinoid 2} and let $V\subseteq S$ be an arbitrary open subset; then we can cover $V$ by rational subsets of the $U_{i}$, each of which is again $f_{\ast}\mathcal{O}_{X}$-adapted by virtue of Lemma~\ref{Subsets of adapted sets are adapted}). In particular, for all $U\in \mathcal{B}$ the topological Banach $\mathcal{O}_{S}(U)$-algebra $\mathcal{O}_{X}(f^{-1}(U))$ is a sheafy Huber ring. Therefore, all members of the basis $\mathcal{B}$ are strongly $f_{\ast}\mathcal{O}_{X}$-adapted affinoid open subspaces of $S$, and thus $f_{\ast}\mathcal{O}_{X}$ is a quasi-coherent sheaf of sheafy topological Banach algebras (respectively, of sheafy Banach algebras). 

It remains to prove that $(f_{\ast}\mathcal{O}_{X}, f_{\ast}\mathcal{O}_{X}^{+})$ enjoys the properties (1) and (2) in the definition~of a quasi-coherent (topological) adic Banach algebra (Definition~\ref{Sheaf of adic Banach algebras}). The first property is trivially satisfied, so we only have to verify property (2) in the said definition. Let $\mathcal{B}$ be the basis of $S$ consisting of all $f_{\ast}\mathcal{O}_{X}$-adapted affinoid open subsets $U$ of $S$ with the property that $f^{-1}(U)$ is affinoid. Then, if $V=U(\frac{f_1,\dots,f_n}{g})$ is a rational subset of some $U\in \mathcal{B}$, we have \begin{equation*}f^{-1}(V)=f^{-1}(U)(\frac{f_1,\dots,f_n}{g})\end{equation*}(by a slight abuse of notation, we identify the elements $f_1,\dots,f_n, g \in \mathcal{O}_{S}(U)$ with their images in $\mathcal{O}_{X}(f^{-1}(U))$). In other words, the restriction map $\mathcal{O}_{X}(f^{-1}(U))\to \mathcal{O}_{X}(f^{-1}(V))$ induces an isomorphism of Huber pairs \begin{equation*}(\mathcal{O}_{X}(f^{-1}(V)), \mathcal{O}_{X}^{+}(f^{-1}(V)))\cong(\mathcal{O}_{X}(f^{-1}(U))\langle\frac{f_1,\dots,f_n}{g}\rangle, \mathcal{O}_{X}(f^{-1}(U))\langle\frac{f_1,\dots,f_n}{g}\rangle^{+}).\end{equation*}This concludes the proof.\end{proof}
\begin{lemma}\label{Locally affinoid 1.5}Let $f: X\to S$ be a locally affinoid morphism of Tate adic spaces and let $U\subseteq S$ be an affinoid open subspace. If $f^{-1}(U)$ is again affinoid, then $U$ is $(f_{\ast}\mathcal{O}_{X}, f_{\ast}\mathcal{O}_{X}^{+})$-adapted.\end{lemma}
\begin{proof}By Proposition~\ref{Locally affinoid 2} the pair $(f_{\ast}\mathcal{O}_{X}, f_{\ast}\mathcal{O}_{X}^{+})$ is a quasi-coherent topological adic Banach $(\mathcal{O}_{S}, \mathcal{O}_{S}^{+})$-algebra (respectively, a quasi-coherent adic Banach $(\mathcal{O}_{S}, \mathcal{O}_{S}^{+})$-algebra), so it makes sense to speak about $(f_{\ast}\mathcal{O}_{X}, f_{\ast}\mathcal{O}_{X}^{+})$-adapted affinoid open subspaces of $S$. The Lemma~follows from Lemma~\ref{Direct images of adic Banach algebras, affinoid case} applied to the morphism of affinoid Tate adic spaces $f\vert_{f^{-1}(U)}: f^{-1}(U)\to U$.\end{proof}
\begin{prop}\label{Locally affinoid 2.5}For a morphism $f:X\to S$ of Tate adic spaces the following are equivalent:\begin{enumerate}[(1)]\item $f$ is locally affinoid.\item The pair $(f_{\ast}\mathcal{O}_{X}, f_{\ast}\mathcal{O}_{X}^{+})$ is a quasi-coherent sheaf of topological adic Banach algebras on $S$ and there exists an open cover $(U_{i})_{i}$ of $S$ by $(f_{\ast}\mathcal{O}_{X},f_{\ast}\mathcal{O}_{X}^{+})$-adapted affinoid open subsets such that the pre-image $f^{-1}(U_{i})$ is affinoid for every $i$.\end{enumerate}\end{prop}
\begin{proof}The implication (2)$\Rightarrow$(1) is trivial and the implication (1)$\Rightarrow$(2) follows from Proposition~\ref{Direct images 2} and Lemma~\ref{Locally affinoid 1.5}.\end{proof}
We want to use what we proved about quasi-coherent sheaves to establish a universal property of the relative adic spectrum of a quasi-coherent topological adic Banach algebra which is analogous to the universal property enjoyed by the relative spectrum of a quasi-coherent algebra in algebraic geometry. 
Let $(S, \varpi_{S})$ be a Tate adic space. Consider the full subcategory of the category of adic spaces over $S$ spanned by the adic spaces $f: X\to S$ locally affinoid over $S$. Let us also note that, for every adic space $f: X\to S$ over $S$, the adic space $X$ naturally inherits a structure of a Tate adic space (cf.~Remark~\ref{Tate adic spaces as relative adic spaces}), so we can apply the previous results on quasi-coherent sheaves of topological adic Banach algebras to the pair $(\mathcal{O}_{X}, \mathcal{O}_{X}^{+})$; in particular, if $f: X\to S$ is a locally affinoid adic space over $S$, then $(f_{\ast}\mathcal{O}_{X}, f_{\ast}\mathcal{O}_{X}^{+})$ is a quasi-coherent sheaf of topological adic Banach $(\mathcal{O}_{S}, \mathcal{O}_{S}^{+})$-algebras, by Proposition~\ref{Direct images 2}. For $(\mathcal{A}, \mathcal{A}^{+})$ a quasi-coherent sheaf of topological adic Banach algebras on $S$ consider the functor \begin{equation*} \{\, \text{Adic spaces locally affinoid}/S\,\} \to \Sets \end{equation*}given on objects by \begin{equation*} (f: X \to S) \mapsto \Hom_{\Ban_{(\mathcal{O}_{S}, \mathcal{O}^{+}_{S})}}((\mathcal{A}, \mathcal{A}^{+}), (f_{\ast}\mathcal{O}_{X}, f_{\ast}\mathcal{O}^{+}_{X})).\end{equation*}The functor sends a morphism \begin{center}\begin{tikzcd} X' \arrow{r} \arrow{dr}{f'} & X\arrow{d}{f} \\ & S \end{tikzcd}\end{center}of locally affinoid adic spaces over $S$ to the map of sets\begin{equation*}\Hom_{\Ban_{(\mathcal{O}_{S}, \mathcal{O}^{+}_{S})}}((\mathcal{A}, \mathcal{A}^{+}), (f_{\ast}\mathcal{O}_{X}, f_{\ast}\mathcal{O}^{+}_{X}))\to \Hom_{\Ban_{(\mathcal{O}_{S}, \mathcal{O}^{+}_{S})}}((\mathcal{A}, \mathcal{A}^{+}), (f'_{\ast}\mathcal{O}_{X}, f'_{\ast}\mathcal{O}^{+}_{X}))\end{equation*}which takes a morphism of topological adic Banach $(\mathcal{O}_{S}, \mathcal{O}_{S}^{+})$-algebras $(\mathcal{A},\mathcal{A}^{+}) \to (f_{\ast}\mathcal{O}_{X},f_{\ast}\mathcal{O}_{X}^{+})$ to the composite morphism\begin{equation*}(\mathcal{A},\mathcal{A}^{+}) \to (f_{\ast}\mathcal{O}_{X},f_{\ast}\mathcal{O}_{X}^{+}) \to (f'_{\ast}\mathcal{O}_{X'},f'_{\ast}\mathcal{O}_{X'}^{+}).\end{equation*}We stress that the definition of the above functor implicitly uses Proposition~\ref{Direct images 2}~which ensures that $(f_{\ast}\mathcal{O}_{X}, f_{\ast}\mathcal{O}_{X}^{+})$ is indeed a quasi-coherent sheaf of topological adic Banach algebras on $S$, for every locally affinoid morphism $f: X\to S$. In the proof of our next Theorem~we repeatedly (and tacitly) use Lemma~\ref{Adic adapted covers 2}.    
\begin{thm}\label{Universal property of the relative spectrum}The above functor is representable by the morphism \begin{equation*}\pi: \underline{\Spa}_{S}(\mathcal{A}, \mathcal{A}^{+}) \to S\end{equation*}from Proposition~\ref{Relative adic spectrum}. \end{thm}
\begin{proof}For every $(\mathcal{A},\mathcal{A}^{+})$-adapted affinoid open subspace $U$ of $S$ we have \begin{equation*} \pi^{-1}(U)=\Spa(\mathcal{A}(U), \mathcal{A}^{+}(U)) \end{equation*}by construction, so the morphism $\pi: \underline{\Spa}_{S}(\mathcal{A}, \mathcal{A}^{+}) \to S$ is indeed locally affinoid.

We have to show that \begin{equation*} \Hom_{\Ban_{(\mathcal{O}_{S}, \mathcal{O}^{+}_{S})}}((\mathcal{A}, \mathcal{A}^{+}), (f_{\ast}\mathcal{O}_{X}, f_{\ast}\mathcal{O}^{+}_{X}))=\Hom_{S}(X, \underline{\Spa}_{S}(\mathcal{A}, \mathcal{A}^{+})) \end{equation*}for any locally affinoid morphism of adic spaces $f: X \to S$. Recall that for every $(\mathcal{A},\mathcal{A}^{+})$-adapted affinoid open subspace $U$ of $S$ we have \begin{equation*}\underline{\Spa}_{S}(\mathcal{A}|_{U}, \mathcal{A}^{+}|_{U})=\Spa(\mathcal{A}(U), \mathcal{A}^{+}(U)).\end{equation*}Suppose that the assertion of the Theorem~holds for every affinoid open subspace $U \subseteq S$ which is $(\mathcal{A},\mathcal{A}^{+})$-adapted and $(f_{\ast}\mathcal{O}_{X},f_{\ast}\mathcal{O}_{X}^{+})$-adapted. Then for every locally affinoid morphism of adic spaces $f: X \to S$ we have \begin{equation*} \Hom_{\Ban_{(\mathcal{O}_{U}, \mathcal{O}^{+}_{U})}}((\mathcal{A}\vert_{U}, \mathcal{A}^{+}\vert_{U}), ((f\vert_{f^{-1}(U)})_{\ast}\mathcal{O}_{f^{-1}(U)}, (f\vert_{f^{-1}(U)})_{\ast}\mathcal{O}^{+}_{f^{-1}(U)})) \end{equation*} \begin{equation*}=\Hom_{S}(f^{-1}(U), \underline{\Spa}_{S}(\mathcal{A}\vert_{U}, \mathcal{A}^{+}\vert_{U})) \end{equation*}for every $(\mathcal{A},\mathcal{A}^{+})$-adapted and $(f_{\ast}\mathcal{O}_{X}, f_{\ast}\mathcal{O}_{X}^{+})$-adapted affinoid open subspace $U\subseteq S$ (we use here that the restricted morphisms $f\vert_{f^{-1}(U)}$ are again locally affinoid). On the other hand, every morphism of sheaves of topological adic Banach $(\mathcal{O}_{S}, \mathcal{O}^{+}_{S})$-algebras $(\mathcal{A}, \mathcal{A}^{+})\to (f_{\ast}\mathcal{O}_{X}, f_{\ast}\mathcal{O}^{+}_{X})$ is determined by the collection of morphisms \begin{equation*} (\mathcal{A}\vert_{U}, \mathcal{A}^{+}\vert_{U}) \to (f_{\ast}\mathcal{O}_{X}\vert_{U}, f_{\ast}\mathcal{O}^{+}_{X}\vert_{U}) = ((f\vert_{f^{-1}(U)})_{\ast}\mathcal{O}_{X}\vert_{f^{-1}(U)}, (f\vert_{f^{-1}(U)})_{\ast}\mathcal{O}^{+}_{X}\vert_{f^{-1}(U)}) \end{equation*}where $U$ ranges over the affinoid open subspaces of $S$ which are both $(\mathcal{A},\mathcal{A}^{+})$-adapted and $(f_{\ast}\mathcal{O}_{X}, f_{\ast}\mathcal{O}_{X}^{+})$-adapted (since these form a basis for the topology). 

Let $(U_{i})_{i}$ be an open cover of $S$ by affinoids which are both $(\mathcal{A},\mathcal{A}^{+})$-adapted and $(f_{\ast}\mathcal{O}_{X},f_{\ast}\mathcal{O}_{X}^{+})$-adapted. A map $(\mathcal{A}, \mathcal{A}^{+}) \to (f_{\ast}\mathcal{O}_{X}, f_{\ast}\mathcal{O}^{+}_{X})$ of quasi-coherent topological adic Banach $(\mathcal{O}_{S}, \mathcal{O}^{+}_{S})$-algebras induces (by restriction) a collection of maps \begin{equation*}(\mathcal{A}\vert_{U_{i}}, \mathcal{A}^{+}\vert_{U_{i}}) \to (f_{\ast}\mathcal{O}_{X}\vert_{U_{i}}, f_{\ast}\mathcal{O}_{X}^{+}\vert_{U_{i}})\end{equation*}such that for any pair of indices $i, j$ the restriction of $(\mathcal{A}\vert_{U_{i}}, \mathcal{A}^{+}\vert_{U_{i}})\to (f_{\ast}\mathcal{O}_{X}\vert_{U_{i}}, f_{\ast}\mathcal{O}_{X}^{+}\vert_{U_{i}})$ to any $(\mathcal{A},\mathcal{A}^{+})$-adapted and $(f_{\ast}\mathcal{O}_{X},f_{\ast}\mathcal{O}_{X}^{+})$-adapted affinoid open subspace $V \subseteq U_{i}\cap U_{j}$ coincides with the restriction of $(\mathcal{A}\vert_{U_{j}}, \mathcal{A}^{+}\vert_{U_{j}}) \to (f_{\ast}\mathcal{O}_{X}\vert_{U_{j}}, f_{\ast}\mathcal{O}^{+}_{X}\vert_{U_{j}})$. By our assumption that the Theorem~holds for $(\mathcal{A},\mathcal{A}^{+})$-adapted and $(f_{\ast}\mathcal{O}_{X}, f_{\ast}\mathcal{O}_{X}^{+})$-adapted open affinoids, this collection corresponds to a collection of morphisms of adic spaces $f^{-1}(U_{i}) \to \Spa(\mathcal{A}(U_{i}), \mathcal{A}^{+}(U_{i}))$ such that for every $i, j$ and every $(\mathcal{A}, \mathcal{A}^{+})$-adapted and $(f_{\ast}\mathcal{O}_{X}, f_{\ast}\mathcal{O}_{X}^{+})$-adapted affinoid open subspace $V$ contained in $U_{i}\cap U_{j}$ the restriction of the morphism \begin{equation*}f^{-1}(U_{i}) \to \Spa(\mathcal{A}(U_{i}), \mathcal{A}^{+}(U_{i}))\end{equation*}to $f^{-1}(V) \subseteq f^{-1}(U_{i}\cap U_{j})=f^{-1}(U_{i})\cap f^{-1}(U_{j})$ coincides with the analogous restriction of $f^{-1}(U_{j}) \to \Spa(\mathcal{A}(U_{j}), \mathcal{A}^{+}(U_{j}))$. Covering $U_{i}\cap U_{j}$ by $(\mathcal{A},\mathcal{A}^{+})$-adapted and $(f_{\ast}\mathcal{O}_{X},f_{\ast}\mathcal{O}_{X}^{+})$-adapted affinoid opens, we see that subsets of the form $f^{-1}(V)$, for $V \subset U_{i}\cap U_{j}$ an $(\mathcal{A},\mathcal{A}^{+})$-adapted and $(f_{\ast}\mathcal{O}_{X},f_{\ast}\mathcal{O}_{X}^{+})$-adapted affinoid open, form an open cover of the subset $f^{-1}(U_{i}\cap U_{j})$ of $X$. Hence the above collection of morphisms defines a unique morphism of adic spaces $X \to \underline{\Spa}_{S}(\mathcal{A}. \mathcal{A}^{+})$, which proves the Theorem~provided that it holds for all affinoid open subspaces of $S$ which are $(\mathcal{A},\mathcal{A}^{+})$-adapted and $(f_{\ast}\mathcal{O}_{X},\mathcal{O}_{X}^{+})$-adapted.

Hence we may assume that $S$ is affinoid and that it is (regarded as its own open subspace) both $(\mathcal{A},\mathcal{A}^{+})$-adapted and $(f_{\ast}\mathcal{O}_{X},f_{\ast}\mathcal{O}_{X}^{+})$-adapted. But in this case $\underline{\Spa}_{S}(\mathcal{A}, \mathcal{A}^{+})$ is also affinoid, namely $\underline{\Spa}_{S}(\mathcal{A}, \mathcal{A}^{+})=\Spa(\mathcal{A}(S), \mathcal{A}^{+}(S))$, so we have \begin{equation*} \Hom_{S}(X, \underline{\Spa}_{S}(\mathcal{A}, \mathcal{A}^{+}))=\Hom_{(\mathcal{O}_{S}(S), \mathcal{O}_{S}^{+}(S))}((\mathcal{A}(S), \mathcal{A}^{+}(S)), (\mathcal{O}_{X}(X), \mathcal{O}_{X}^{+}(X)),\end{equation*}where the set on the right hand side is the set of morphisms of Huber pairs over $(\mathcal{O}_{S}(S), \mathcal{O}_{S}^{+}(S))$. The latter set is canonically in bijection with the set of morphisms of topological adic Banach $(\mathcal{O}_{S}, \mathcal{O}^{+}_{S})$-algebras $(\mathcal{A}, \mathcal{A}^{+}) \to (f_{\ast}\mathcal{O}_{X}, f_{\ast}\mathcal{O}^{+}_{X})$, by Lemma~\ref{Morphisms of quasi-coherent adic Banach algebras}, so we are done.\end{proof}
We can now prove our main result on locally affinoid morphisms and the relative adic spectrum. 
\begin{thm}\label{The affinoid equivalence}Let $(S, \varpi_{S})$ be a Tate adic space. The functor \begin{equation*}(\mathcal{A}, \mathcal{A}^{+}) \mapsto \underline{\Spa}_{S}(\mathcal{A}, \mathcal{A}^{+}) \end{equation*}is an anti-equivalence between the category of quasi-coherent sheaves of topological adic Banach $(\mathcal{O}_{S}, \mathcal{O}_{S}^{+})$-algebras and the full subcategory of the category of adic spaces over $S$ spanned by adic spaces locally affinoid over $S$.\end{thm}
\begin{proof} For a quasi-coherent topological adic Banach $(\mathcal{O}_{S}, \mathcal{O}_{S}^{+})$-algebra $(\mathcal{A}, \mathcal{A}^{+})$ denote by \begin{equation*}\pi: \underline{\Spa}_{S}(\mathcal{A}, \mathcal{A}^{+}) \to S\end{equation*}the canonical map as in Proposition~\ref{Relative adic spectrum}. For every $(\mathcal{A},\mathcal{A}^{+})$-adapted affinoid open subspace $U$ of $S$ we have \begin{equation*}\pi^{-1}(U)=\Spa(\mathcal{A}(U), \mathcal{A}^{+}(U)),\end{equation*}which gives rise to a canonical isomorphism\begin{equation*}(\pi_{\ast}\mathcal{O}_{\underline{\Spa}_{S}(\mathcal{A}, \mathcal{A}^{+})}, \pi_{\ast}\mathcal{O}_{\underline{\Spa}_{S}(\mathcal{A}, \mathcal{A}^{+})}^{+}) \cong (\mathcal{A}, \mathcal{A}^{+}).\end{equation*}From the latter isomorphism we see that the functor $(\mathcal{A}, \mathcal{A}^{+}) \mapsto \underline{\Spa}_{S}(\mathcal{A}, \mathcal{A}^{+})$ is fully faithful. To show essential surjectivity, let $f: X \to S$ be a locally affinoid morphism of adic spaces. By Proposition~\ref{Direct images 2}, $(f_{\ast}\mathcal{O}_{X}, f_{\ast}\mathcal{O}_{X}^{+})$ is a quasi-coherent sheaf of topological adic Banach algebras on $S$. By Theorem~\ref{Universal property of the relative spectrum}, the identity map $(f_{\ast}\mathcal{O}_{X}, f_{\ast}\mathcal{O}_{X}^{+})\to (f_{\ast}\mathcal{O}_{X}, f_{\ast}\mathcal{O}_{X}^{+})$ induces a canonical morphism \begin{equation*}X\to\underline{\Spa}_{S}(f_{\ast}\mathcal{O}_{X}, f_{\ast}\mathcal{O}_{X}^{+}) \end{equation*}of adic spaces. For any $(f_{\ast}\mathcal{O}_{X}, f_{\ast}\mathcal{O}_{X}^{+})$-adapted affinoid open $U\subseteq S$ the above morphism restricts to the canonical morphism \begin{equation*}f^{-1}(U)\to \Spa(\mathcal{O}_{X}(f^{-1}(U)), \mathcal{O}_{X}^{+}(f^{-1}(U))).\end{equation*}The latter is an isomorphism if $f^{-1}(U)$ is affinoid, see property (1) in Proposition~\ref{Relative adic spectrum}. But by Proposition~\ref{Locally affinoid 2.5}~there exists an open cover $(U_{i})_{i}$ of $S$ by $(f_{\ast}\mathcal{O}_{X}, f_{\ast}\mathcal{O}_{X}^{+})$-adapted affinoid open subsets such that $f^{-1}(U_{i})$ is affinoid for every $i$. It follows that $X\to\underline{\Spa}_{S}(f_{\ast}\mathcal{O}_{X}, f_{\ast}\mathcal{O}_{X}^{+})$ restricts to an isomorphism over each of the open subsets $f^{-1}(U_{i})$. Since these open subsets cover $X$, we see that $X \to \underline{\Spa}_{S}(f_{\ast}\mathcal{O}_{X}, f_{\ast}\mathcal{O}_{X})$ is an isomorphism.\end{proof}
\begin{cor}\label{Locally affinoid 3}For $f: X\to S$ a locally affinoid morphism of Tate adic spaces, the pre-image $f^{-1}(U)$ of an affinoid open subspace $U\subseteq S$ is affinoid if and only if $U$ is $(f_{\ast}\mathcal{O}_{X}, f_{\ast}\mathcal{O}_{X}^{+})$-adapted.
\end{cor}
\begin{proof}If $f^{-1}(U)$ is affinoid, then $U$ is $(f_{\ast}\mathcal{O}_{X}, f_{\ast}\mathcal{O}_{X}^{+})$-adapted by Lemma~\ref{Locally affinoid 1.5}. The converse follows from the construction of the relative adic spectrum and Theorem~\ref{The affinoid equivalence}.\end{proof}
\begin{cor}\label{Locally affinoid 4}Let $f: X\to S$ be a locally affinoid morphism of Tate adic spaces and let $g: X'\to X$ be a morphism such that $f\circ g$ is locally affinoid. Then $g$ is locally affinoid.\end{cor}
\begin{proof}Cover $S$ by affinoid open subspaces $U_{i}$ which are both $(f_{\ast}\mathcal{O}_{X}, f_{\ast}\mathcal{O}_{X}^{+})$-adapted and $((f\circ g)_{\ast}\mathcal{O}_{X}, (f\circ g)_{\ast}\mathcal{O}_{X}^{+})$-adapted. Then $f^{-1}(U_{i})$ and $(f\circ g)^{-1}(U_{i})=g^{-1}(f^{-1}(U_{i}))$ is affinoid for all $i$, by Corollary~\ref{Locally affinoid 3}. But the open subsets $f^{-1}(U_{i})$ cover $X$, so we see that $g$ is locally affinoid.\end{proof}
We have the following application of Corollary~\ref{Locally affinoid 3}.
\begin{prop}\label{Higher direct images}Let $f: X\to S$ be a locally affinoid morphism of Tate adic spaces. Let $\mathcal{F}$ be a pseudo-coherent sheaf or a Banach bundle $\mathcal{F}$ on $X$. Then $R^{i}f_{\ast}\mathcal{F}=0$ for all $i>0$.\end{prop}
\begin{proof}It suffices to prove that the stalk $(R^{i}f_{\ast}\mathcal{F})_{x}$ of $R^{i}f_{\ast}\mathcal{F}$ for each $x\in S$ is zero. We have \begin{equation*}(R^{i}f_{\ast}\mathcal{F})_{x} = \varinjlim_{x \in V}H^{i}(f^{-1}(V), \mathcal{F})\end{equation*}where $V$ ranges over all open neighbourhoods of $x$. Hence it suffices to prove that $H^{i}(f^{-1}(V), \mathcal{F})=0$ where $V$ ranges over some cofinal system of open neighbourhoods of $x$. By Proposition~\ref{Direct images 2}~the pair $(f_{\ast}\mathcal{O}_{X}, f_{\ast}\mathcal{O}_{X}^{+})$ is a quasi-coherent sheaf of topological adic Banach $(\mathcal{O}_{S},\mathcal{O}_{S}^{+})$-algebras. Since the $(f_{\ast}\mathcal{O}_{X}, f_{\ast}\mathcal{O}_{X}^{+})$-adapted open affinoid subsets $V\subseteq S$ form a basis for the topology on $S$, the $(f_{\ast}\mathcal{O}_{X}, f_{\ast}\mathcal{O}_{X}^{+})$-adapted open affinoid neighbourhoods of $x$ form a cofinal system. But by Corollary~\ref{Locally affinoid 3} the pre-image $f^{-1}(V)$ of any $(f_{\ast}\mathcal{O}_{X},f_{\ast}\mathcal{O}_{X}^{+})$-adapted affinoid open $V$ is an affinoid open of $X$. By \cite{Kedlaya-Liu2}, Theorem~2.5.5, when $\mathcal{F}$ is a pseudo-coherent sheaf (respectively, by the unpublished result of Kedlaya referred to in Example \ref{Banach bundles} when $\mathcal{F}$ is a Banach bundle), it follows that $H^{i}(f^{-1}(V), \mathcal{O}_{X})=0$ for all $i>0$ and for all $(f_{\ast}\mathcal{O}_{X}, f_{\ast}\mathcal{O}_{X}^{+})$-adapted affinoid open neighbourhoods $V$ of $x$.\end{proof}
\begin{cor}\label{Cohomology and locally affinoid morphisms}Let $f: X\to S$ be a locally affinoid morphism of Tate adic spaces. Let $\mathcal{F}$ be a pseudo-coherent sheaf or a Banach bundle on $X$. Then \begin{equation*}H^{i}(X, \mathcal{F}) = H^{i}(S, f_{\ast}\mathcal{F}) \end{equation*}for all $i\geq 0$.\end{cor}
\begin{example}\label{Cohomology of non-adapted affinoids}The above Corollary~also reveals another feature of the theory of quasi-coherent sheaves on Tate adic spaces which appears unusual from the point of view of algebraic geometry: The higher cohomology of a quasi-coherent sheaf of topological Banach modules $\mathcal{F}$ on an affinoid adic space $X$ may be non-zero, at least if $X$ is not $\mathcal{F}$-adapted. To see this, we again consider Example 9.1.2 from \cite{Berkeley}; that is, we let again \begin{equation*}V=\{\, \vert x\vert=1\,\}\cup \{\, \vert y\vert=1\,\} \subseteq \mathbb{B}_{K}^{2}\end{equation*}(for some nonarchimedean field $K$ with some fixed pseudo-uniformizer $\varpi\in K$). Let $i: V\hookrightarrow \mathbb{B}_{K}^{2}$ be the inclusion map. By the discussion in loc.~cit.~we know that $i$ is locally affinoid. Recall that by Proposition~\ref{Direct images} the sheaf of topological Banach modules $i_{\ast}\mathcal{O}_{V}$ on $X$ is quasi-coherent. However, by Corollary~\ref{Cohomology and locally affinoid morphisms} we have \begin{equation*}H^{1}(\mathbb{B}_{K}^{2}, i_{\ast}\mathcal{O}_{V}) =H^{1}(V, \mathcal{O}_{V})=\widehat{\bigoplus}_{m,n>0}Kx^{-n}y^{-m}\neq 0.\end{equation*}\end{example}
\begin{example}We can also give an example of two locally affinoid morphisms whose composition is not locally affinoid. Set notation as in the previous example. The inclusion morphism $i: V\hookrightarrow \mathbb{B}_{K}^{2}$ is locally affinoid by \cite{Berkeley}, Example~9.1.2. Let $f: \mathbb{B}_{K}^{2}\to \Spa(K, K^{\circ})$ be the structure morphism. As a morphism between affinoid adic spaces, $f$ is locally affinoid. But the composition $f\circ i: V\to \Spa(K, K^{\circ})$ is not locally affinoid since the underlying topological space of $\Spa(K, K^{\circ})$ is a point and $V$ is not affinoid. We also see that the statement of Corollary~\ref{Cohomology and locally affinoid morphisms}~becomes false if we let $\mathcal{F}$ be an arbitrary quasi-coherent sheaf of topological Banach modules on $X$. Indeed, every abelian sheaf on $\Spa(K, K^{\circ})$ has vanishing higher cohomology but \begin{equation*}H^{1}(\mathbb{B}_{K}^{2}, i_{\ast}\mathcal{O}_{V})\neq 0=H^{1}(\Spa(K, K^{\circ}), f_{\ast}i_{\ast}\mathcal{O}_{V}),\end{equation*}as we saw in the previous example.\end{example}

\section{Partially proper morphisms}
 
\subsection{Partially proper, separated and rationalizing morphisms} 
  
The connection between quasi-coherent sheaves of topological Banach algebras and locally affinoid morphisms becomes even tighter if we restrict our attention to those locally affinoid morphisms which are also partially proper in the sense of the following definition. 
\begin{mydef}[Partially proper morphism]\label{Partially proper morphisms}A morphism of adic spaces $f: X\to S$ is called \textit{partially proper} if it is quasi-separated and for every complete sheafy Huber pair $(R, R^{+})$ and every commutative square \begin{center}\begin{tikzcd}\Spa(R, R^{\circ}) \arrow{r} \arrow{d} & X \arrow{d}{f} \\ \Spa(R, R^{+}) \arrow{r} & S \end{tikzcd}\end{center}there exists a unique dotted arrow making the diagram \begin{center}\begin{tikzcd}\Spa(R, R^{\circ}) \arrow{r} \arrow{d} & X \arrow{d}{f} \\ \Spa(R, R^{+}) \arrow[dotted]{ur} \arrow{r} & S \end{tikzcd}\end{center}commute. We call $f$ \textit{proper} if it is partially proper and quasi-compact.\end{mydef}
The above definition is an immediate analog for adic spaces of the notion of partially proper morphisms of v-stacks on $\Perf$ introduced by Scholze in \cite{Diamonds}, Definition~18.4. Note that the definition does not contain any finite-type assumptions on $f$. For separated morphisms locally of $^{+}$weakly finite type between locally Noetherian analytic adic spaces, our definition coincides with Huber's original definition of partially proper morphisms of locally Noetherian analytic adic spaces (\cite{Huber3}, Definition~1.3.3), by \cite{Huber3}, Corollary~1.3.9 and \cite{Ludwig20}, Proposition~6.13. In general, the two definitions differ since we de not require our partially proper morphisms to be locally weakly of finite type.

In algebraic geometry, every proper morphism between affine schemes is finite. The analog of this assertion for partially proper (or even proper) morphisms of adic spaces is far from being true, even if the partially proper morphism in question is locally weakly of finite type (it does, however, become true if we assume that the morphism is of finite type, see Proposition~\ref{Proper locally affinoid morphisms of finite type} and Remark~\ref{Proper locally affinoid morphisms of finite type 2} below). We can describe partially proper morphisms between affinoid adic spaces as follows.
\begin{lemma}\label{Partial properness for affinoids}A morphism of affinoid adic spaces $f: X\to S$ is partially proper if and only if $X=\Spa(\mathcal{O}_{X}(X), \overline{\mathcal{O}_{S}^{+}(S)+\mathcal{O}_{X}(X)^{\circ\circ}})$, where the bar on the right hand side means integral closure inside $\mathcal{O}_{X}(X)$.\end{lemma}
\begin{proof}Note that $\overline{\mathcal{O}_{S}^{+}(S)+\mathcal{O}_{X}(X)^{\circ\circ}}$ is precisely the smallest ring of integral elements of $\mathcal{O}_{X}(X)$ containing the image of $\mathcal{O}_{S}^{+}(S)$. Suppose that $f$ is partially proper. Consider the commutative diagram \begin{center}\begin{tikzcd} \Spa(\mathcal{O}_{X}(X), \mathcal{O}_{X}^{+}(X)) \arrow{r}{=} \arrow{d} & X \arrow{d}{f} \\ \Spa(\mathcal{O}_{X}(X), \overline{\mathcal{O}_{S}^{+}(S)+\mathcal{O}_{X}(X)^{\circ\circ}}) \arrow{r} & S, \end{tikzcd}\end{center}where the left vertical arrow is induced by the inclusion map of complete sheafy Huber pairs \begin{equation*}(\mathcal{O}_{X}(X), \overline{\mathcal{O}_{S}^{+}(S)+\mathcal{O}_{X}(X)^{\circ\circ}})\hookrightarrow (\mathcal{O}_{X}(X), \mathcal{O}_{X}^{+}(X)).\end{equation*}By partial properness of $f$ we see that the left vertical arrow has a left-inverse. But then the above inclusion map of Huber pairs has a right-inverse. Thus the inclusion map \begin{equation*}\overline{\mathcal{O}_{S}^{+}(S)+\mathcal{O}_{X}(X)^{\circ\circ}}\hookrightarrow \mathcal{O}_{X}^{+}(X)\end{equation*}is surjective, as desired.

Conversely, assume that $\mathcal{O}_{X}^{+}(X)$ is the smallest ring of integral elements of $\mathcal{O}_{X}(X)$ containing the image of $\mathcal{O}_{S}^{+}(S)$. Consider a commutative diagram \begin{center}\begin{tikzcd}\Spa(R, R^{\circ}) \arrow{r} \arrow{d} & X \arrow{d}{f} \\ \Spa(R, R^{+}) \arrow{r} & S. \end{tikzcd}\end{center}Then the ring homomorphism $\mathcal{O}_{X}^{+}(X)\to R^{\circ}$ induced by the upper horizontal arrow is a homomorphism of $\mathcal{O}_{S}^{+}(S)$-algebras. Since $\mathcal{O}_{X}^{+}(X)$ is the smallest ring of integral elements of $\mathcal{O}_{X}(X)$ containing the image of $\mathcal{O}_{S}^{+}(S)$, the image of $\mathcal{O}_{X}^{+}(X)\to R^{\circ}$ must be contained in the smallest ring of integral elements of $R$ containing the image of $\mathcal{O}_{S}^{+}(S)$. A fortiori, the image of $\mathcal{O}_{X}^{+}(X)\to R^{\circ}$ must be contained in $R^{+}$.\end{proof}
\begin{example}\label{Example of a partially proper morphism}Let $K$ be a perfectoid field. The above lemma shows that the morphism of affinoid adic spaces \begin{equation*}\Spa(K\langle T^{1/p^{\infty}}\rangle, K^{\circ}\langle T^{1/p^{\infty}}\rangle)\to \Spa(K\langle T\rangle, K^{\circ}\langle T\rangle)\end{equation*}is proper. This gives a simple example of a morphism which is proper in our sense but not weakly of finite type. More generally, the Lemma~entails that every pro-finite morphism of affinoid adic spaces $f: X\to S$ (i.~e., a morphism which is the inverse limit of an inverse system of finite morphisms of affinoid adic spaces $f_{i}: X_{i}\to S$) is partially proper.\end{example} 
The following definition~is an extension of Huber's notion of morphisms locally of $^{+}$weakly finite type (\cite{Huber3}, Definition~1.2.1ii)).
\begin{mydef}[Rationalizing morphism]\label{Rationalizing morphism}~\begin{enumerate}[(1)] \item A morphism of Huber pairs $(A, A^{+})\to (B, B^{+})$ is called \textit{rationalizing} if there exists a finite subset $E\subseteq B^{+}$ such that $B^{+}$ is the smallest ring of integral elements of $B$ which contains $E$ and the image of $A^{+}\to B^{+}$. In other words, $B^{+}$ is the integral closure of the subring $A^{+}[E\cup B^{\circ\circ}]$ inside $B$. \item A morphism $f: X\to S$ of adic spaces is called \textit{rationalizing} if for every $x\in X$ there exist affinoid open subspaces $V\subseteq X$ and $U\subseteq S$ such that $x\in V$, $f(V)\subseteq U$ and the induced morphism of Huber pairs $(\mathcal{O}_{S}(U), \mathcal{O}_{S}^{+}(U))\to (\mathcal{O}_{X}(V), \mathcal{O}_{X}^{+}(V))$ is rationalizing.\end{enumerate}\end{mydef}
\begin{example}\label{Open immersions}If $(B, B^{+})$ is a rational localization of a Huber pair $(A, A^{+})$, then, writing $B=A\langle\frac{f_1,\dots, f_n}{g}\rangle$ for some elements $f_1,\dots, f_n, g\in A^{+}$ generating an open ideal of $A$, the finite subset $E=\{\, \frac{f_1}{g},\dots, \frac{f_n}{g}\,\}$ shows that the morphism $(A, A^{+})\to (B, B^{+})$ is rationalizing. It follows that every open immersion of adic spaces is rationalizing. In fact, the same argument shows that open immersions are locally of finite type. More generally, every morphism locally of $^{+}$-weakly finite type (and, in particular, every morphism locally of finite type) is rationalizing, by definition.\end{example}
\begin{example}\label{Partially proper morphisms are rationalizing}On the other hand, every proper morphism between affinoid adic spaces is rationalizing, by Lemma~\ref{Partial properness for affinoids}, and this includes many infinite-type examples such as Example~\ref{Example of a partially proper morphism}.\end{example} 
\begin{lemma}\label{First Lemma on rationalizing morphisms}Let $\varphi: (A, A^{+})\to (B, B^{+})$ and $\psi: (B, B^{+})\to (C, C^{+})$ be morphisms of Huber pairs. \begin{enumerate}[(1)]\item If $\varphi$ and $\psi$ are rationalizing, then so is $\psi\circ\varphi$. \item If $\psi\circ\varphi$ is rationalizing, then $\psi$ is rationalizing.\end{enumerate}\end{lemma}
\begin{proof}(1) Let $E\subseteq B^{+}$ be a finite subset such that $B^{+}$ is the smallest ring of integral elements of $B$ containing $E$ and the image of $A^{+}\to B^{+}$ and let $F\subseteq C^{+}$ be a finite subset such that $C^{+}$ is the smallest ring of integral elements of $C$ containing $F$ and the image of $B^{+}\to C^{+}$. Then $F\cup\psi(E)$ is a finite subset of $C^{+}$ such that $C^{+}$ is the smallest ring of integral elements of $C$ containing $F\cup\varphi(E)$ and the image of $A^{+}\to C^{+}$. 

(2) Suppose that $\psi\circ\varphi: (A, A^{+})\to (C, C^{+})$ is rationalizing and let $E\subseteq C^{+}$ be a finite subset of $C^{+}$ such that $C^{+}$ is the smallest ring of integral elements of $C$ containing $E$ and $\psi(\varphi(A^{+}))$. In other words, $C^{+}$ is integral over the subring generated by $E\cup C^{\circ\circ}$ and the image of $A^{+}\to C^{+}$. But then, a fortiori, $C^{+}$ is also integral over the subring generated by $E\cup C^{\circ\circ}$ and the image of $B^{+}\to C^{+}$.\end{proof} 
\begin{lemma}\label{Rationalizing morphisms and rational subsets}If $f: X\to S$ is a rationalizing morphism of adic spaces, then, for every affinoid open subspace $U\subseteq S$ and every open subspace $V\subseteq X$ with $f(V)\subseteq U$, there exists an affinoid open cover $(W_{k})_{k}$ of $V$ such that, for every $k$, the morphism of Huber pairs $(\mathcal{O}_{S}(U), \mathcal{O}_{S}^{+}(U))\to (\mathcal{O}_{X}(W_{k}), \mathcal{O}_{X}^{+}(W_{k}))$ is rationalizing.\end{lemma}
\begin{proof}Let $(V_{i})_{i}$ (respectively, $(U_{i})_{i}$) be an affinoid open cover of $X$ (respectively, a family of affinoid open subsets of $S$) such that $f(V_{i})\subseteq U_{i}$ and the morphisms of Huber pairs $(\mathcal{O}_{S}(U_{i}), \mathcal{O}_{S}^{+}(U_{i}))\to (\mathcal{O}_{X}(V_{i}), \mathcal{O}_{X}^{+}(V_{i}))$ are rationalizing. For each $i$ let $(U_{ij})_{j}$ be an open cover of $U_{i}\cap U$ by affinoid open subsets which are rational subsets both of $U$ and of $U_{i}$ and, for each $i, j$ and let $(V_{ik})_{k}$ be an affinoid open cover of $V\cap V_{i}$ by rational subsets of $V_{i}$. Then $f$ restricts to morphisms of affinoid adic spaces $V_{ik}\to U$, so the pre-images \begin{equation*}W_{ijk}=f^{-1}(U_{ij})\cap V_{ik}\end{equation*}of the rational subsets $U_{ij}$ under these morphisms are rational subsets of $V_{ik}$, and hence also rational subsets of $V_{i}$. Then $(W_{ijk})_{j, k}$ is an open cover of $V_{i}\cap V$ by rational subsets of $V_{i}$, for each $i$. In particular, the composite morphism of Huber pairs \begin{equation*}(\mathcal{O}_{S}(U_{i}), \mathcal{O}_{S}^{+}(U_{i}))\to (\mathcal{O}_{X}(V_{i}), \mathcal{O}_{X}^{+}(V_{i}))\to (\mathcal{O}_{X}(W_{ijk}), \mathcal{O}_{X}^{+}(W_{ijk}))\end{equation*}is rationalizing for all $i, j, k$, by Example~\ref{Open immersions} and Lemma~\ref{First Lemma on rationalizing morphisms}(1). But, by the universal property of rational localizations, the above morphism of Huber pairs $(\mathcal{O}_{S}(U_{i}), \mathcal{O}_{S}^{+}(U_{i}))\to (\mathcal{O}_{X}(W_{ijk}), \mathcal{O}_{X}^{+}(W_{ijk}))$ factors through $(\mathcal{O}_{S}(U_{ij}), \mathcal{O}_{S}^{+}(U_{ij}))$. By Lemma~\ref{First Lemma on rationalizing morphisms}(2), the resulting morphism $(\mathcal{O}_{S}(U_{ij}), \mathcal{O}_{S}^{+}(U_{ij}))\to (\mathcal{O}_{X}(W_{ijk}), \mathcal{O}_{X}^{+}(W_{ijk}))$ is rationalizing. Since each $U_{ij}$ is a rational subset of $U$, the morphism of Huber pairs $(\mathcal{O}_{S}(U), \mathcal{O}_{S}^{+}(U))\to (\mathcal{O}_{S}(U_{ij}), \mathcal{O}_{S}(U_{ij}))$ is rationalizing for all $i, j$, by Example~\ref{Open immersions}. Hence the composite morphism \begin{equation*}(\mathcal{O}_{S}(U), \mathcal{O}_{S}^{+}(U))\to (\mathcal{O}_{S}(U_{ij}), \mathcal{O}_{S}^{+}(U_{ij}))\to (\mathcal{O}_{X}(W_{ijk}), \mathcal{O}_{X}^{+}(W_{ijk}))\end{equation*}is rationalizing, by Lemma~\ref{First Lemma on rationalizing morphisms}(1). This shows that $(W_{ijk})_{i, j, k}$ is an affinoid open cover of $V$ with the sought-for property.\end{proof}
\begin{lemma}\label{Rationalizing morphisms are local}Let $f: X\to S$ be a morphism of adic spaces. The following are equivalent: \begin{enumerate}[(1)]\item $f$ is rationalizing. \item For every open subspace $U\subseteq S$ and every open subspace $V\subseteq X$ with $f(V)\subseteq U$, the morphism $f\vert_{V}: V\to U$ is rationalizing. \item There exist open covers $(V_{i})_{i}$ of $X$ and $(U_{i})_{i}$ of $S$ with $f(V_{i})\subseteq U_{i}$ for all $i$ and the maps $f\vert_{V_{i}}: V_{i}\to U_{i}$ are rationalizing. \end{enumerate}\end{lemma}
\begin{proof}The implication (2)$\Rightarrow$(3) is trivial and the implication (3)$\Rightarrow$(1) follows from the definition of rationalizing morphisms. It remains to prove that (1) implies (2). To this end, let $V\subseteq X$, $U\subseteq S$ be open subspaces such that $f(V)\subseteq U$. To prove that the map $f\vert_{V}: V\to U$ is rationalizing, we may shrink $U$ and assume that $U$ is affinoid. In this case, we apply Lemma~\ref{Rationalizing morphisms and rational subsets} to produce an affinoid open cover $(W_{k})_{k}$ of $V$ such that the corresponding maps of Huber pairs \begin{equation*}(\mathcal{O}_{S}(U), \mathcal{O}_{S}^{+}(U))\to (\mathcal{O}_{X}(W_{k}), \mathcal{O}_{X}^{+}(W_{k}))\end{equation*}are rationalizing. It follows that $f\vert_{V}: V\to U$ is rationalizing. \end{proof} 
\begin{lemma}\label{Composition of rationalizing morphisms}The composition of two rationalizing morphisms is rationalizing.\end{lemma}
\begin{proof}Let $f: X\to S$ and $g: S\to S'$ be rationalizing morphisms of adic spaces. It suffices to prove that for every affinoid open subspace $U'\subseteq S'$ the pre-image $f^{-1}(g^{-1}(U'))$ has an affinoid open cover $(W_{j})_{j}$ such that the morphism of Huber pairs $(\mathcal{O}_{S'}(U'), \mathcal{O}_{S'}^{+}(U'))\to (\mathcal{O}_{X}(W_{j}), \mathcal{O}_{X}^{+}(W_{j}))$ is rationalizing. By Lemma~\ref{Rationalizing morphisms and rational subsets}, there is an affinoid open cover $(U_{i})_{i}$ of $g^{-1}(U')$ such that the morphisms of Huber pairs $(\mathcal{O}_{S'}(U'), \mathcal{O}_{S'}^{+}(U'))\to (\mathcal{O}_{S}(U_{i}), \mathcal{O}_{S}^{+}(U_{i}))$ are rationalizing. Applying the same lemma to each $U_{i}$, we obtain for each $i$ an affinoid open cover $(W_{ij})_{j}$ of $f^{-1}(U_{i})$ such that $(\mathcal{O}_{S}(U_{i}), \mathcal{O}_{S}^{+}(U_{i}))\to (\mathcal{O}_{X}(W_{ij}), \mathcal{O}_{X}^{+}(W_{ij}))$ is rationalizing for all $i, j$. By Lemma~\ref{First Lemma on rationalizing morphisms}(1), $(W_{ij})_{i,j}$ is an affinoid open cover of $f^{-1}(g^{-1}(U'))$ with the desired property. \end{proof}
\begin{lemma}\label{Rationalizing morphisms and factorizations}If $f: X\to S$ and $g: S\to S'$ are morphisms of adic spaces such that $g\circ f$ is rationalizing, then $f$ is rationalizing.\end{lemma}
\begin{proof}Let $(U_{i})_{i}$ (respectively, $(U_{i}')_{i}$) be affinoid open covers of $S$ (respectively, $S'$) with $g(U_{i})\subseteq U_{i}'$. Let $(V_{i})_{i}$ be an affinoid open cover of $X$ with $f(V_{i})\subseteq U_{i}$. By Lemma~\ref{Rationalizing morphisms and rational subsets}, there exists, for every $i$, an affinoid open cover $(W_{ij})_{j}$ of $V_{i}$ such that the morphism of Huber pairs $(\mathcal{O}_{S'}(U_{i}'), \mathcal{O}_{S'}^{+}(U_{i}'))\to (\mathcal{O}_{S}(U_{i}), \mathcal{O}_{S}^{+}(U_{i}))\to (\mathcal{O}_{X}(W_{ij}), \mathcal{O}_{X}^{+}(W_{ij}))$ is rationalizing for every $j$. We conclude by Lemma~\ref{First Lemma on rationalizing morphisms}(2).\end{proof}
Using the notion of rationalizing morphisms, we obtain the following criterion for partial properness, which is a non-Noetherian analog of \cite{Ludwig20}, Proposition~6.13.
\begin{lemma}\label{Partial properness and affinoid fields}Suppose that $f: X\to S$ is a quasi-separated morphism of analytic adic spaces satisfying the following two properties: \begin{enumerate}[(1)] \item $f$ is rationalizing. \item For every Huber pair $(K, K^{+})$ where $K$ is a nonarchimedean field and every commutative square \begin{center}\begin{tikzcd}\Spa(K, K^{\circ}) \arrow{r} \arrow{d} & X \arrow{d}{f} \\ \Spa(K, K^{+}) \arrow{r} & S \end{tikzcd}\end{center}there exists a unique dotted arrow making the diagram \begin{center}\begin{tikzcd}\Spa(K, K^{\circ}) \arrow{r} \arrow{d} & X \arrow{d}{f} \\ \Spa(K, K^{+}) \arrow[dotted]{ur} \arrow{r} & S \end{tikzcd}\end{center}commute.\end{enumerate}Then $f$ is partially proper.\end{lemma}
\begin{proof}The proof of \cite{Ludwig20}, Proposition~6.13, applies verbatim, using the notion of a rationalizing morphism in place of the notion of a morphism of $^{+}$-weakly finite type.\end{proof}
We also have a notion of separated morphisms characterized by a valuative criterion analogous to property (2) in Lemma~\ref{Partial properness and affinoid fields}.
\begin{mydef}[Separated morphism]\label{Definition of separated morphisms}A morphism of adic spaces $f: X\to S$ is called \textit{separated} if it is quasi-separated and for every Huber pair $(K, K^{+})$ where $K$ is a nonarchimedean field and every commutative square \begin{center}\begin{tikzcd}\Spa(K, K^{\circ}) \arrow{r} \arrow{d} & X \arrow{d}{f} \\ \Spa(K, K^{+}) \arrow{r} & S \end{tikzcd}\end{center}there exists \textit{at most one} dotted arrow making the diagram \begin{center}\begin{tikzcd}\Spa(K, K^{\circ}) \arrow{r} \arrow{d} & X \arrow{d}{f} \\ \Spa(K, K^{+}) \arrow[dotted]{ur} \arrow{r} & S \end{tikzcd}\end{center}commute.\end{mydef}
\begin{rmk}By the same argument as in the proof of \cite{Diamonds}, Proposition~5.11, the above definition is also equivalent to the usual condition that $X\to X\times_{S}X$ be a closed immersion, whenever $S$ is analytic and the fiber product $X\times_{S}X$ exists as an adic space. Our definition of separated morphisms is more general than the definition from Huber's book \cite{Huber3} (see loc.~cit., Definition~1.3.1) since we do not assume our separated morphisms to be locally weakly of finite type.\end{rmk}
\begin{lemma}\label{Partial properness is local on the target}Let $f: X\to S$ be a morphism of adic spaces. Consider the following assertions: \begin{enumerate}[(1)]\item $f$ is partially proper. \item For every affinoid open subspace $U\subseteq S$, the morphism \begin{equation*}f\vert_{f^{-1}(U)}: f^{-1}(U)\to U\end{equation*}is partially proper. \item There exists an affinoid open cover $(U_{i})_{i}$ of $S$ such that the maps \begin{equation*}f\vert_{f^{-1}(U_{i})}: f^{-1}(U_{i})\to U_{i}\end{equation*}are partially proper. \end{enumerate}Then (1)$\Rightarrow$(2)$\Rightarrow$(3), and if $f$ is rationalizing and $S$ is analytic, then all three assertions are equivalent. \end{lemma}
\begin{proof}(2) obviously implies (3). 

(1)$\Rightarrow$(2): Since $f$ is quasi-separated, it is clear that $f\vert_{f^{-1}(U)}$ is again quasi-separated. Consider a commutative square \begin{center}\begin{tikzcd}\Spa(R, R^{\circ}) \arrow{r} \arrow{d} & f^{-1}(U) \arrow{d}{f\vert_{f^{-1}(U)}} \\ \Spa(R, R^{+}) \arrow{r} & U \end{tikzcd}\end{center}for a complete sheafy Huber pair $(R, R^{+})$. By partial properness of $f$, there exists a dotted arrow making the diagram \begin{center}\begin{tikzcd}\Spa(R, R^{\circ}) \arrow{r} \arrow{d} & X \arrow{d}{f} \\ \Spa(R, R^{+}) \arrow[dotted]{ur} \arrow{r} & S \end{tikzcd}\end{center}commute. But then the image of the dotted arrow is actually contained in $f^{-1}(U)$, the image of $\Spa(R, R^{+})\to S$ being contained in $U$. 

(3)$\Rightarrow$(1) when $f$ is rationalizing and $S$ analytic: The morphism $f$ is quasi-separated since each $f\vert_{f^{-1}(U_{i})}$ is quasi-separated. Consider a commutative square of the form \begin{center}\begin{tikzcd}\Spa(K, K^{\circ}) \arrow{r} \arrow{d} & X \arrow{d}{f} \\ \Spa(K, K^{+}) \arrow{r} & S. \end{tikzcd}\end{center}Recall that $\Spa(K, K^{+})$ is a totally ordered chain of points (ordered by specialization). Since the $U_{i}$ cover $S$, there exists some index $i_{0}$ with the property that the image of the unique closed point of $\Spa(K, K^{+})$ is contained in $U_{i_{0}}$. Since all other points are generizations of the unique closed point, the map $\Spa(K, K^{+})\to S$ is continuous and $U_{i_{0}}$ is open in $S$, we see that the whole image of $\Spa(K, K^{+})\to S$ is contained in $U_{i_{0}}$. Since \begin{equation*}f\vert_{f^{-1}(U_{i_{0}})}: f^{-1}(U_{i_{0}})\to U_{i_{0}}\end{equation*}is partially proper, we get a dotted arrow making the diagram \begin{center}\begin{tikzcd}\Spa(K, K^{\circ}) \arrow{r} \arrow{d} & f^{-1}(U_{i_{0}}) \subseteq X \arrow{d} \\ \Spa(K, K^{+}) \arrow[dotted]{ur} \arrow{r} & U_{i_{0}}\subseteq S \end{tikzcd}\end{center}commute. We conclude by Lemma~\ref{Partial properness and affinoid fields} that $f$ is partially proper. 
\end{proof}
\begin{cor}\label{Partially proper morphisms are rationalizing 2}Every proper locally affinoid morphism of adic spaces $f: X\to S$ is rationalizing.\end{cor}
\begin{proof}Let $(U_{i})_{i}$ be an affinoid open cover of $S$ such that the pre-images $f^{-1}(U_{i})$ are affinoid. By Lemma~\ref{Rationalizing morphisms are local}, it suffices to prove that each of the morphisms \begin{equation*}f\vert_{f^{-1}(U_{i})}: f^{-1}(U_{i})\to U_{i}\end{equation*}is rationalizing. By Lemma~\ref{Partial properness is local on the target}, these are proper morphisms between affinoid adic spaces, hence they are rationalizing by Example~\ref{Partially proper morphisms are rationalizing}. \end{proof}

\subsection{Finite morphisms}

A special class of proper locally affinoid morphisms are the finite morphisms of adic spaces which, in the locally Noetherian case, were introduced by Huber. 
\begin{mydef}[Finite morphism, \cite{Huber3}, (1.4.2), (1.4.4)]\label{Finite morphisms, definition}\begin{enumerate}[(1)]\item A morphism of Huber pairs \begin{equation*}(A, A^{+})\to (B, B^{+})\end{equation*}is called finite if it is topologically of finite type and the maps $A\to B$ and $A^{+}\to B^{+}$ are integral. \item  
A morphism $f: X\to S$ of adic spaces is called finite if there exists an affinoid open cover $(U_{i})_{i}$ of $S$ such that, for all $i$, the pre-image $f^{-1}(U_{i})$ is affinoid and the corresponding map of Huber pairs \begin{equation*}(\mathcal{O}_{S}(U_{i}), \mathcal{O}_{S}^{+}(U_{i}))\to (\mathcal{O}_{X}(f^{-1}(U_{i})), \mathcal{O}_{X}^{+}(f^{-1}(U_{i})))\end{equation*}is finite. \end{enumerate}\end{mydef}
The following lemma of Huber ensures that this definition is equivalent to a more obvious notion of finite morphisms. We include the proof of this lemma here for the reader's convenience. 
\begin{lemma}[cf.~\cite{Huber0}, Lemma~3.12.10]\label{Huber's lemma on finite morphisms}If a morphism of complete Huber pairs $(A, A^{+})\to (B, B^{+})$ is finite, then $A\to B$ is a finite ring map, the topology on $B$ is its natural topology of finitely generated $A$-module and $B^{+}$ is the integral closure of (the image of) $A^{+}$ in $B$. Conversely, if $(A, A^{+})\to (B, B^{+})$ is a morphism of complete Tate Huber pairs such that $A\to B$ is finite and $A^{+}\to B^{+}$ is integral, then $(A, A^{+})\to (B, B^{+})$ is finite. \end{lemma}
\begin{proof}We follow the argument in the proof of \cite{Huber0}, Lemma~3.12.10, which works without any Noetherian assumptions. We include the proof here for the reader's convenience. 

Suppose that the map of Huber pairs $(A, A^{+})\to (B, B^{+})$ is finite in the sense of Definition~\ref{Finite morphisms, definition}. Since $A^{+}\to B^{+}$ is integral by assumption, it suffices to prove that $B$ is module-finite over $A$ and that the topology on $B$ is its natural topology of finitely generated $A$-module. Since $(A, A^{+})\to (B, B^{+})$ is topologically of finite type, there exist rings of definition $A_{0}$ and $B_{0}$ of $A$ and $B$, respectively, such that $A\to B$ restricts to a map $A_{0}\to B_{0}$ which is (strictly) topologically of finite type, $B$ is finitely generated over $A\cdot B_{0}$ and $B^{+}$ is the integral closure of $B_{0}$ in $B$. Let $f_1,\dots, f_n$ be elements of $B_{0}$ such that $A_{0}[f_1,\dots, f_n]$ is dense in $B_{0}$. Since $f_1,\dots, f_n$ are integral over $A^{+}$, they are integral over some finite $A_{0}$-subalgebra $A_{1}$ of $A^{+}$. Then $A_{1}$ and $B_{1}=A_{1}\cdot B_{0}$ are rings of definition of $A$ and $B$, respectively, such that $A\to B$ restricts to a map $A_{1}\to B_{1}$ and $B_{1}$ has a dense $A_{1}$-subalgebra which is module-finite over $A_{1}$ (namely, $A_{1}[f_1,\dots, f_n]$). Consequently, if $I$ is an ideal of definition of $A_{1}$, then $A_{1}/I\to B_{1}/IB_{1}$ is a finite ring map. We conclude by topological Nakayama's Lemma that $A_{1}\to B_{1}$ is a finite ring map. Hence $A\to B$ is a finite ring map. Moreover, since $B$ is finitely generated over $A\cdot B_{0}$, it is a fortiori finitely generated over $A\cdot B_{1}$. Hence we obtain a surjective $A$-module map \begin{equation*}A^{m}=A_{1}^{m}\otimes_{A_{1}}A\to B_{1}\otimes_{A_{1}}A\to B\end{equation*}such that the image of $A_{1}^{m}$ under the map is the ring of definition $B_{1}$. This shows that the topology on $B$ is its natural topology of finitely generated $A$-module. 

For the converse implication, assume that $A$ is Tate and let $f_1,\dots, f_n\in B$ be elements generating $B$ as an $A$-module. Let $(A_{0}, \varpi)$ be a pair of definition of $A$ with $A_{0}\subseteq A^{+}$. Up to multiplying $f_1,\dots, f_n$ by a large enough power of $\varpi$, we may assume that $f_1,\dots, f_n$ are integral over $A_{0}$. Then the $A_{0}$-subalgebra $B_{0}=A_{0}[f_1,\dots, f_n]$ of $B$ is (module-)finite. In particular, it is a bounded subring of $B$. By Banach's open mapping theorem, the topology on $B$ is necessarily its natural topology of finitely generated $A$-module. Hence the image of $A_{0}^{n}$ under the surjective $A$-module map \begin{equation*}A^{n}\to B, (a_1,\dots, a_{n})\mapsto \sum_{i=1}^{n}a_{i}f_{i},\end{equation*}is open in $B$. It follows that the bounded subring $B_{0}$ is also open in $B$, thus a ring of definition of $B$. This implies that the morphism of complete Huber pairs $(A, A^{+})\to (B, B^{+})$ is topologically of finite type, and hence finite. \end{proof}
\begin{lemma}\label{Finite implies proper}Every finite morphism of analytic adic spaces $f: X\to S$ is proper. \end{lemma}
\begin{proof}Choose an open cover of $S$ by affinoid open subspaces $U_{i}\subseteq S$ such that $f^{-1}(U_{i})$ is affinoid and $f^{-1}(U_{i})\to U_{i}$ induces a finite morphism of Huber pairs for all $i$. By Lemma~\ref{Partial properness is local on the target}, it suffices to prove that each $f^{-1}(U_{i})\to U_{i}$ is proper, so we may assume that $S$ and $X$ are affinoid. Then finiteness implies that $\mathcal{O}_{X}^{+}(X)$ is integral over $\mathcal{O}_{S}^{+}(S)$. But then, a fortiori, $\mathcal{O}_{X}^{+}(X)$ is the smallest integrally open subring of $\mathcal{O}_{X}(X)$ containing the image of $\mathcal{O}_{S}^{+}(S)$. We conclude by Lemma~\ref{Partial properness for affinoids} that $f: X\to S$ is proper.  \end{proof}
Recall that every proper affine morphism of schemes is finite. Using the same approach as in the proof of Lemma~\ref{Huber's lemma on finite morphisms}, we can show that the analog of this statement is also true for finite type morphisms of uniform analytic adic spaces. 
\begin{prop}\label{Proper locally affinoid morphisms of finite type}Let $f: X\to S$ be a morphism locally of finite type between uniform analytic adic spaces. Then $f$ is finite if and only if it is locally affinoid and proper. \end{prop}
\begin{proof}One implication is a special case of Lemma~\ref{Finite implies proper}. Suppose that $f$ is locally affinoid and proper. We show that $f$ must be finite. Choose an affinoid open cover $(U_{i})_{i}$ of $S$ such that the pre-image $f^{-1}(U_{i})$ is affinoid for every $i$. By Lemma~\ref{Partial properness is local on the target}, $f\vert_{f^{-1}(U_{i})}: f^{-1}(U_{i})\to U_{i}$ is a proper morphism of affinoid adic spaces, for each $i$. Therefore, we may assume that $X$ and $S$ are affinoid. In this case, $\mathcal{O}_{X}^{+}(X)$ is the smallest integrally closed open subring of $\mathcal{O}_{X}(X)$ containing the image of $\mathcal{O}_{S}^{+}(S)$, by Lemma~\ref{Partial properness for affinoids}. By \cite{Huber2}, Proposition~3.6 (see also \cite{Ludwig20}, Theorem~2.25, which contains a proof of this result\footnote{Note that there are no locally Noetherian assumptions in the statement or proof of this result in loc.~cit.}), the morphism of Huber pairs \begin{equation*}(A, A^{+})=(\mathcal{O}_{S}(S), \mathcal{O}_{S}^{+}(S))\to (B, B^{+})=(\mathcal{O}_{X}(X), \mathcal{O}_{X}^{+}(X))\end{equation*}is topologically of finite type. Then there exist rings of definition $A_{0}\subseteq A^{+}$ and $B_{0}\subseteq B^{+}$ of $A$ and $B$, respectively, such that $A\to B$ restricts to an adic map of adic rings $A_{0}\to B_{0}$ which is topologically of finite type and $B^{+}$ is the integral closure of $B_{0}$ in $B$. As in the proof of Lemma~\ref{Huber's lemma on finite morphisms}, we choose elements $f_1,\dots, f_n\in B_{0}$ such that the finite-type $A_{0}$-subalgebra $A_{0}[f_1,\dots, f_n]$ is dense in $B_{0}$. Then $B_{0}$ is contained in the (topological) closure $B_{1}$ of $A^{+}[f_1,\dots, f_n]$ inside $B^{+}$. It follows that $B^{+}$ is integral over its topologically finitely generated $A^{+}$-subalgebra $B_{1}$. 

Since $B^{+}$ is the smallest integrally closed open $A^{+}$-subalgebra of $\mathcal{O}_{X}(X)$, we have \begin{equation*}B^{+}=\overline{A^{+}+\varpi B_{1}}\end{equation*}for any ideal of definition $I$ of $A_{0}$. It follows that the map \begin{equation*}A^{+}/IA^{+}\to B_{1}/I B_{1}\end{equation*}is integral. Since $B_{1}$ was topologically of finite type over $A^{+}$, the above map $A^{+}/IA^{+}\to B_{1}/IB_{1}$ is also of finite type. Consequently, this map is finite. Since the Huber rings $A$ and $B$ are uniform, $A^{+}$ and $B^{+}$ are $I$-adically complete. In particular, $A^{+}$ and $B_{1}$ are $I$-adically complete. By topological Nakayama's lemma, $B_{1}$ must then be finite over $A^{+}$. In particular, $B^{+}$ is integral over $A^{+}$. \end{proof}
\begin{rmk}\label{Proper locally affinoid morphisms of finite type 2}In the literature, the conclusion of Proposition~\ref{Proper locally affinoid morphisms of finite type} (for morphisms between affinoid adic spaces, to which the proposition reduces) is usually shown under the hypothesis that $S$ and $X$ be locally Noetherian (but not necessarily uniform), see \cite{Huber0}, Proposition~1.4.6, and \cite{Ludwig20}, Corollary~7.11. The content of our proposition is that, at least under the uniformity assumption, the same result is true without any Noetherian hypotheses. \end{rmk}

\subsection{Proper locally affinoid morphisms via quasi-coherent Banach algebras}

We can now relate partially proper morphisms to quasi-coherent Banach algebras. For every quasi-coherent sheaf $\mathcal{A}$ of sheafy topological Banach $\mathcal{O}_{S}$-algebras on a Tate adic space $S$ we can define a subsheaf $\overline{\mathcal{A}^{+}}$ of $\mathcal{A}$ by setting \begin{equation*}\overline{\mathcal{A}^{+}}(U)=\overline{\mathcal{O}_{S}^{+}(U)+\mathcal{A}(U)^{\circ\circ}}\end{equation*}for every $\mathcal{A}$-adapted affinoid open $U\subseteq S$, where the bar on the right denotes integral closure inside $\mathcal{A}(U)$.
\begin{prop}\label{Properness and adapted subsets} For every quasi-coherent sheaf of sheafy topological Banach $\mathcal{O}_{S}$-algebras $\mathcal{A}$ on a Tate adic space $S$, the pair $(\mathcal{A}, \overline{\mathcal{A}^{+}})$ is a quasi-coherent topological adic Banach $(\mathcal{O}_{S}, \mathcal{O}_{S}^{+})$-algebra and every strongly $\mathcal{A}$-adapted affinoid open subspace $U\subseteq S$ is $(\mathcal{A}, \overline{\mathcal{A}^{+}})$-adapted.\end{prop}
\begin{proof}Let $U\subseteq S$ be a strongly $\mathcal{A}$-adapted affinoid open subspace and let $V=U(\frac{f_1,\dots, f_n}{g})\subseteq U$ be a rational subset of $U$ defined by elements $f_1,\dots, f_n, g\in \mathcal{O}_{S}(U)$ generating the unit ideal in $\mathcal{O}_{S}(U)$. The subring \begin{equation*}\overline{\mathcal{A}^{+}}(V)=\overline{\mathcal{O}_{S}^{+}(V)+\mathcal{A}(V)^{\circ\circ}}=\overline{\mathcal{O}_{S}^{+}(V)+\mathcal{A}(U)\langle\frac{f_1,\dots, f_n}{g}\rangle^{\circ\circ}}\end{equation*}is the smallest integrally closed and open subring of $\mathcal{A}(V)=\mathcal{A}(U)\langle\frac{f_1,\dots,f_n}{g}\rangle$ containing $\mathcal{O}_{S}^{+}(V)$. On the other hand, consider the rational localization\begin{equation*}(\mathcal{A}(V), \mathcal{A}(V)^{+})=(\mathcal{A}(U)\langle\frac{f_1,\dots,f_n}{g}\rangle, \mathcal{A}(U)\langle\frac{f_1,\dots,f_n}{g}\rangle^{+})\end{equation*}of the complete Huber pair $(\mathcal{A}(U), \overline{\mathcal{A}^{+}}(U))$. The subring \begin{equation*}\mathcal{A}(U)\langle\frac{f_1,\dots, f_n}{g}\rangle^{+}=\overline{\overline{\mathcal{A}^{+}}(U)\langle\frac{f_1,\dots, f_n}{g}\rangle},\end{equation*}where the upper bar denotes integral closure inside $\mathcal{A}(V)$, is an integrally closed subring of $\mathcal{A}(V)$ which contains \begin{equation*}\mathcal{O}_{S}^{+}(V)=\overline{\mathcal{O}_{S}^{+}(U)\langle\frac{f_1,\dots,f_n}{g}\rangle}\end{equation*}and is open. Therefore we see that \begin{equation*}\overline{\mathcal{A}^{+}}(V)\subseteq \overline{\overline{\mathcal{A}^{+}}(U)\langle\frac{f_1,\dots, f_n}{g}\rangle}.\end{equation*}Conversely, the right hand side is the smallest integrally closed open subring of $\mathcal{A}(V)$ which contains the image of $\overline{\mathcal{A}^{+}}(U)$ and the elements $\frac{f_1}{g},\dots, \frac{f_n}{g}$. But $\overline{\mathcal{A}^{+}}(V)$ contains the elements $\frac{f_1}{g}, \dots, \frac{f_n}{g}$ since it contains \begin{equation*}\mathcal{O}_{S}^{+}(V)\supseteq \mathcal{O}_{S}^{+}(U)\langle\frac{f_1,\dots, f_n}{g}\rangle.\end{equation*}It follows that \begin{equation*}\overline{\mathcal{A}^{+}}(V)\supseteq \overline{\overline{\mathcal{A}^{+}}(U)\langle\frac{f_1,\dots,f_n}{g}\rangle}.\end{equation*}\end{proof}
\begin{lemma}\label{Relative adic spectrum is partially proper}For any Tate adic space $S$ and any quasi-coherent sheafy topological Banach $\mathcal{O}_{S}$-algebra $\mathcal{A}$ on $S$, the canonical morphism \begin{equation*}\pi: \underline{\Spa}_{S}(\mathcal{A}, \overline{\mathcal{A}^{+}})\to S\end{equation*}is proper.\end{lemma}
\begin{proof}Since every locally affinoid morphism is quasi-compact, it suffices to prove that the morphism in question is partially proper. For every $(\mathcal{A}, \overline{\mathcal{A}^{+}})$-adapted (equivalently, strongly $\mathcal{A}$-adapted) affinoid open subspace $U\subseteq S$, we have \begin{equation*}\pi^{-1}(U)=\Spa(\mathcal{A}(U), \overline{\mathcal{O}_{S}^{+}(U)+\mathcal{A}(U)^{\circ\circ}}).\end{equation*}Hence the assertion follows from Lemma~\ref{Partial properness for affinoids} and Lemma~\ref{Partial properness is local on the target}, the family of all $(\mathcal{A}, \overline{\mathcal{A}^{+}})$-adapted affinoid open subsets being a basis for the topology on $S$.\end{proof}
For $\mathcal{A}$ a quasi-coherent sheaf of sheafy topological Banach $\mathcal{O}_{S}$-algebras on a Tate adic space $S$ consider the functor \begin{equation*}\{\, \text{Adic spaces proper and locally affinoid}/S\,\}\to \Sets \end{equation*}given on objects by $(f: X\to S)\mapsto \Hom_{\Ban-\alg_{\mathcal{O}_{S}}}(\mathcal{A}, f_{\ast}\mathcal{O}_{X})$, where we use the fact that every adic space $X$ over $S$ inherits a natural structure of a Tate adic space, which allows us to talk about quasi-coherent sheaves of topological Banach algebras on $X$. As in Section 5, the above functor sends a morphism $X'\to X$ of proper and locally affinoid adic spaces over $S$ to the map of sets \begin{equation*}\Hom_{\Ban-\alg_{\mathcal{O}_{S}}}(\mathcal{A}, f_{\ast}\mathcal{O}_{X})\to\Hom_{\Ban-\alg_{\mathcal{O}_{S}}}(\mathcal{A}, f'_{\ast}\mathcal{O}_{X})\end{equation*}which takes a morphism of sheafy topological Banach $\mathcal{O}_{S}$-algebras $\mathcal{A}\to f_{\ast}\mathcal{O}_{X}$ to the composite map $\mathcal{A}\to f_{\ast}\mathcal{O}_{X}\to f'_{\ast}\mathcal{O}_{X}$.
\begin{thm}\label{Universal property of the relative spectrum 2}The above functor is representable by the morphism \begin{equation*}\pi: \underline{\Spa}_{S}(\mathcal{A}, \overline{\mathcal{A}^{+}})\to S.\end{equation*}\end{thm}
\begin{proof}We already know from Lemma~\ref{Relative adic spectrum is partially proper} that $\pi$ is a proper morphism of adic spaces. By Theorem~\ref{Universal property of the relative spectrum}~we know that for locally affinoid morphisms $f: X\to S$ there are identifications \begin{equation*}\Hom_{S}(X, \underline{\Spa}_{S}(\mathcal{A}, \overline{\mathcal{A}^{+}}))=\Hom_{\Ban_{(\mathcal{O}_{S}, \mathcal{O}_{S}^{+})}}((\mathcal{A}, \overline{\mathcal{A}^{+}}), (f_{\ast}\mathcal{O}_{X}, f_{\ast}\mathcal{O}_{X}^{+}))\end{equation*}that are natural in $f$. For every strongly $\mathcal{A}$-adapted affinoid open subspace $U\subseteq S$, the image of $\mathcal{A}(U)^{\circ\circ}$ under any continuous homomorphism $\mathcal{A}(U)\to \mathcal{O}_{X}(f^{-1}(U))$ is contained in $\mathcal{O}_{X}(f^{-1}(U))^{\circ\circ}$, so the image of $\overline{\mathcal{A}^{+}}(U)=\overline{\mathcal{O}_{S}^{+}(U)+\mathcal{A}(U)^{\circ\circ}}$ is contained in $\overline{\mathcal{O}_{S}^{+}(U)+\mathcal{O}_{X}(f^{-1}(U))^{\circ\circ}}$. A fortiori, the image of $\overline{\mathcal{A}^{+}}(U)$ under any morphism \begin{equation*}\mathcal{A}(U)\to \mathcal{O}_{X}(f^{-1}(U))\end{equation*}is contained in $\mathcal{O}_{X}^{+}(f^{-1}(U))$, the subring $\overline{\mathcal{O}_{S}^{+}(U)+\mathcal{O}_{X}(f^{-1}(U))^{\circ\circ}}$ of $\mathcal{O}_{X}(f^{-1}(U))$ being the smallest integrally closed open subring of $\mathcal{O}_{X}(f^{-1}(U))$ containing $\mathcal{O}_{S}^{+}(U)$. It follows that we have a natural identification \begin{equation*}\Hom_{\Ban_{(\mathcal{O}_{S}, \mathcal{O}_{S}^{+})}}((\mathcal{A}, \overline{\mathcal{A}^{+}}), (f_{\ast}\mathcal{O}_{X}, f_{\ast}\mathcal{O}_{X}^{+}))=\Hom_{\Ban-\alg_{\mathcal{O}_{S}}}(\mathcal{A}, f_{\ast}\mathcal{O}_{X}).\end{equation*}This concludes the proof.\end{proof}
We obtain the following analog of Theorem~\ref{The affinoid equivalence}.
\begin{thm}\label{The affinoid equivalence 2}Let $S$ be a Tate adic space. The functors \begin{equation*}\mathcal{A}\mapsto \underline{\Spa}_{S}(\mathcal{A}, \overline{\mathcal{A}^{+}})\end{equation*}and \begin{equation*}X\mapsto f_{\ast}\mathcal{O}_{X}\end{equation*}define quasi-inverse anti-equivalences of categories between the category of quasi-coherent sheaves of sheafy topological Banach $\mathcal{O}_{S}$-algebras and the full subcategory of adic spaces over $S$ spanned by adic spaces locally affinoid and proper over $S$.\end{thm}
\begin{proof}Since the pre-image in $\underline{\Spa}_{S}(\mathcal{A}, \overline{\mathcal{A}^{+}})$ of any strongly $\mathcal{A}$-adapted affinoid open subspace $U$ of $S$ is $\Spa(\mathcal{A}(U), \overline{\mathcal{A}^{+}}(U))$, the functor \begin{equation*}\mathcal{A}\mapsto \pi_{\ast}\mathcal{O}_{\underline{\Spa}_{S}(\mathcal{A}, \overline{\mathcal{A}^{+}})}\end{equation*}is isomorphic to the identity. To show that the functor $X\mapsto \underline{\Spa}_{S}(f_{\ast}\mathcal{O}_{X}, \overline{(f_{\ast}\mathcal{O}_{X})^{+}})$ is isomorphic to the identity, we let $f: X\to S$ be a locally affinoid and proper morphism of adic spaces. For any strongly $(f_{\ast}\mathcal{O}_{X}, f_{\ast}\mathcal{O}_{X}^{+})$-adapted affinoid open $U\subseteq S$, the pre-image $f^{-1}(U)$ is affinoid by Corollary~\ref{Locally affinoid 3} and then the restricted morphism \begin{equation*}f\vert_{f^{-1}(U)}: f^{-1}(U)\to U\end{equation*}is a proper morphism of affinoid adic spaces by Lemma~\ref{Partial properness is local on the target}. By Lemma~\ref{Partial properness for affinoids}, $f^{-1}(U)$ is of the form \begin{equation*}f^{-1}(U)=\Spa(\mathcal{O}_{X}(f^{-1}(U)), \overline{\mathcal{O}_{S}^{+}(U)+\mathcal{O}_{X}(f^{-1}(U))^{\circ\circ}}).\end{equation*}In other words, the morphism \begin{equation*}X\to \underline{\Spa}_{S}(f_{\ast}\mathcal{O}_{X}, \overline{(f_{\ast}\mathcal{O}_{X})^{+}})\end{equation*}induced by the identity map $f_{\ast}\mathcal{O}_{X}\to f_{\ast}\mathcal{O}_{X}$ via Theorem~\ref{Universal property of the relative spectrum 2}~is an isomorphism over every $(f_{\ast}\mathcal{O}_{X}, f_{\ast}\mathcal{O}_{X}^{+})$-adapted affinoid open subspace $U\subseteq S$. Since $(f_{\ast}\mathcal{O}_{X}, f_{\ast}\mathcal{O}_{X}^{+})$-adapted affinoid opens form a basis for the topology of $S$, this implies that \begin{equation*}X\to \underline{\Spa}_{S}(f_{\ast}\mathcal{O}_{X}, \overline{(f_{\ast}\mathcal{O}_{X})^{+}})\end{equation*}is an isomorphism.\end{proof}
\begin{cor}\label{Locally affinoid proper}Let $f: X\to S$ be a locally affinoid proper morphism of Tate adic spaces. Then an affinoid open subspace $U\subseteq S$ is $(f_{\ast}\mathcal{O}_{X}, f_{\ast}\mathcal{O}_{X}^{+})$-adapted if and only if it is strongly $f_{\ast}\mathcal{O}_{X}$-adapted.\end{cor}
\begin{proof}By Theorem~\ref{The affinoid equivalence 2}, we have \begin{equation*}(f_{\ast}\mathcal{O}_{X}, f_{\ast}\mathcal{O}_{X}^{+})=(f_{\ast}\mathcal{O}_{X}, \overline{f_{\ast}\mathcal{O}_{X}^{+}}),\end{equation*}so the assertion follows from Proposition~\ref{Properness and adapted subsets}. \end{proof}
This result also specializes to a result about finite morphisms which in the case of adic spaces covered by adic spectra of Huber pairs admitting Noetherian rings of definition was already observed by Huber in \cite{Huber0}, Korollar~3.12.12. In the non-Noetherian case, we cannot use coherent sheaves, so, to formulate our result, we need the following definition.
\begin{mydef}[Quasi-coherent sheaf of finite type]\label{Quasi-coherent of finite type}Let $S$ be a Tate adic space. \begin{enumerate}[(1)]\item A quasi-coherent sheaf of topological Banach modules $\mathcal{F}$ on $S$ is said to be of finite type if there exists an $\mathcal{F}$-adapted affinoid open cover $(U_{i})_{i}$ such that the $\mathcal{O}_{S}(U_{i})$-module $\mathcal{F}(U_{i})$ is finitely generated. \item A quasi-coherent sheaf of sheafy topological Banach algebras $\mathcal{A}$ on $S$ is said to be of finite type if we can choose an $\mathcal{A}$-adapted affinoid open cover as in (1) which is also strongly $\mathcal{A}$-adapted. \end{enumerate}\end{mydef}
By Example~\ref{Pseudocoherent}, this notion of quasi-coherent sheaves (of topological Banach modules) of finite type generalizes the notion of pseudocoherent sheaves on $S$, which was introduced and studied by Kedlaya-Liu in \cite{Kedlaya-Liu2}, \S2.5. In particular, when the Tate adic space $S$ is locally Noetherian, it coincides with the notion of coherent sheaves studied by Kedlaya-Liu in \cite{Kedlaya-Liu2}, \S2.3. 
\begin{thm}\label{The affinoid equivalence 3}Let $S$ be a Tate adic space. The functors from Theorem~\ref{The affinoid equivalence 2} restrict to an anti-equivalence of categories between the category of quasi-coherent sheaves of sheafy (topological) Banach algebras $\mathcal{A}$ of finite type on $S$ and the full subcategory of the category of adic spaces $f: X\to S$ over $S$ spanned by adic spaces finite over $S$. \end{thm}
\begin{proof}Let $\mathcal{A}$ be a quasi-coherent sheafy Banach algebra of finite type. Let $(U_{i})_{i}$ be a strongly $\mathcal{A}$-adapted affinoid open cover such that $\mathcal{A}(U_{i})$ is module-finite over $\mathcal{O}_{S}(U_{i})$. Let \begin{equation*}f: X=\underline{\Spa}_{S}(\mathcal{A}, \overline{\mathcal{A}^{+}})\to S\end{equation*}be the canonical morphism. Then \begin{equation*}f^{-1}(U_{i})=\Spa(\mathcal{A}(U_{i}), \overline{\mathcal{A}^{+}}(U_{i})).\end{equation*}In particular, the maps $\mathcal{O}_{S}(U_{i})\to \mathcal{A}(U_{i})=\mathcal{O}_{X}(f^{-1}(U_{i}))$ are finite. We have to prove that also the maps $\mathcal{O}_{S}^{+}(U_{i})\to \mathcal{O}_{X}^{+}(f^{-1}(U_{i}))$ are integral. By \cite{Kedlaya17}, Corollary~1.1.12, the topology on the Banach module $\mathcal{A}(U_{i})$ over the complete Tate ring $A_{i}=\mathcal{O}_{S}(U_{i})$ is its natural topology of finitely generated $A_{i}$-module. Hence, if $A_{0, i}$ is a ring of definition of $A_{i}$ contained in $A_{i}^{+}=\mathcal{O}_{S}^{+}(U_{i})$, there exists a finitely generated $A_{0, i}$-submodule of $\mathcal{A}(U_{i})$ which is also a ring of definition of $\mathcal{A}(U_{i})$. It follows that the integral closure of $A_{0, i}$ in $\mathcal{A}(U_{i})$ is open. In particular, the integral closure of $A_{i}^{+}$ in $\mathcal{A}(U_{i})$ is open. But, by definition, $\overline{\mathcal{A}^{+}}(U_{i})$ is the smallest integrally closed open subring of $\mathcal{A}(U_{i})$ containing $A_{i}^{+}$. It follows that $\overline{\mathcal{A}^{+}}(U_{i})=\mathcal{O}_{X}^{+}(f^{-1}(U_{i}))$ is the integral closure of $A_{i}^{+}$ in $\mathcal{A}(U_{i})$. That is, the map $\mathcal{O}_{S}^{+}(U_{i})\to \mathcal{O}_{X}^{+}(f^{-1}(U_{i}))$ is integral for every $i$. We conclude that the morphism $f: X\to S$ is finite. 

Conversely, let $f: X\to S$ be a finite morphism. Let $(U_{i})_{i}$ be an affinoid open cover of $S$ such that the pre-images $f^{-1}(U_{i})$ are affinoid and the corresponding maps of Huber pairs are finite. By Lemma~\ref{Locally affinoid 1.5}, the affinoid open subspaces $U_{i}$ are strongly $f_{\ast}\mathcal{O}_{X}$-adapted. Hence $f_{\ast}\mathcal{O}_{X}$ is a quasi-coherent topological Banach algebra of finite type. \end{proof}
As an application of our results, we obtain a short proof (in the case of Tate adic spaces) of a theorem of Huber which is referenced in \cite{Huber3}, (1.4.4), but whose proof has not appeared in published form.\footnote{In the case of adic spaces whose rings of sections have Noetherian rings of definition, the result follows from \cite{Huber0}, Korollar~3.12.12, which is a special case of our Theorem~\ref{The affinoid equivalence 3}. In \cite{Huber3}, (1.4.4), Huber refers to a manuscript in preparation titled ``Coherent sheaves on adic spaces''. It seems that this manuscript never appeared as a publication.}. 
\begin{cor}[Huber]\label{Huber's theorem}Let $f: X\to S$ be a finite morphism of locally Noetherian Tate adic spaces. For every Tate affinoid open subspace $U\subseteq S$, the pre-image $f^{-1}(U)$ is affinoid and the map of Huber pairs \begin{equation*}(\mathcal{O}_{S}(U), \mathcal{O}_{S}^{+}(U))\to (\mathcal{O}_{X}(f^{-1}(U)), \mathcal{O}_{X}^{+}(f^{-1}(U)))\end{equation*}is finite. \end{cor}
\begin{proof}By Theorem~\ref{The affinoid equivalence 3}, the quasi-coherent sheaf of topological Banach algebras $f_{\ast}\mathcal{O}_{X}$ is of finite type. Since the analytic adic space $S$ is locally Noetherian, this means that it is a coherent sheaf on $S$ in the sense of Kedlaya-Liu \cite{Kedlaya-Liu2}, \S2.3. By \cite{Kedlaya-Liu2}, Theorem~2.3.3, this entails that every affinoid open subspace $U\subseteq S$ is $f_{\ast}\mathcal{O}_{S}$-adapted and $\mathcal{O}_{X}(f^{-1}(U))$ is a finitely generated $\mathcal{O}_{S}(U)$-module. In particular, $\mathcal{O}_{X}(f^{-1}(U))$ is again a strongly Noetherian complete Tate ring and hence sheafy. It follows that every affinoid open subspace $U\subseteq S$ is strongly $f_{\ast}\mathcal{O}_{X}$-adapted. But the morphism $f$ is proper by Lemma~\ref{Finite implies proper}, so, by Corollary~\ref{Locally affinoid proper}, every affinoid open subset $U\subseteq S$ is $(f_{\ast}\mathcal{O}_{X}, f_{\ast}\mathcal{O}_{X}^{+})$-adapted. We conclude by Corollary~\ref{Locally affinoid 3} that the pre-image $f^{-1}(U)$ of every affinoid open subspace $U\subseteq S$ is affinoid. Moreover, we have seen that, for any affinoid open subspace $U\subseteq S$, the $\mathcal{O}_{S}(U)$-module $\mathcal{O}_{X}(f^{-1}(U))$ is finitely generated and hence its topology coincides with its natural topology of finitely generated $\mathcal{O}_{S}(U)$-module by \cite{Kedlaya17}, Corollary~1.1.12. We conclude as in the proof of Theorem~\ref{The affinoid equivalence 3} that the map of Huber pairs $(\mathcal{O}_{S}(U), \mathcal{O}_{S}^{+}(U))\to (\mathcal{O}_{X}(f^{-1}(U)), \mathcal{O}_{X}^{+}(f^{-1}(U)))$ is finite. \end{proof}

\subsection{Huber's compactifications}

As an application of our results on proper locally affinoid morphisms, we want to prove a compactification theorem (in the style of Huber's \cite{Huber3}, \S5.1) for locally affinoid morphisms of Tate adic spaces. We use the following straightforward analog of \cite{Huber3}, Definition~5.1.1 (note, however, that we do not assume our adic spaces to be locally Noetherian, which is a standing assumption in \cite{Huber3}, and that our notion of partially proper morphisms is broader than Huber's since we do not require partially proper morphisms to be locally weakly of finite type).
\begin{mydef}[Universal compactification]\label{Definition of universal compactifications}For a morphism $f: X\to S$ of analytic adic spaces a \textit{compactification} of $f$ is a commutative triangle of analytic adic spaces \begin{center}\begin{tikzcd}X \arrow[hook]{r}{j} \arrow{d}{f} & \overline{X} \arrow{dl}{\overline{f}} \\ S, \end{tikzcd}\end{center}where $j$ is a locally closed embedding and $\overline{f}$ is partially proper and rationalizing. A compactification of $f$ as above is called a universal compactification if for every commutative triangle of adic spaces \begin{center}\begin{tikzcd}X \arrow{r}{h} \arrow{d}{f} & Y \arrow{dl}{g} \\ S \end{tikzcd}\end{center}with $g$ partially proper and rationalizing there exists a unique morphism $i: \overline{X}\to Y$ making the diagram \begin{center}\begin{tikzcd}X\arrow{dd}{f} \arrow[hook]{dr}{j} \arrow{drr}{h} & & \\ & \overline{X} \arrow{r}{i} \arrow{dl}{\overline{f}} & Y \arrow{dll}{g} \\ S \end{tikzcd}\end{center}commute.\end{mydef}
It is clear that a universal compactification is unique whenever it exists. In order to prove our compactification theorem, we collect some further notions and results from \cite{Huber3} while avoiding any of the original Noetherian or finiteness assumptions. The following lemma is a non-Noetherian partial analog to \cite{Huber3}, Lemma~1.1.10iii). 
\begin{lemma}\label{Specializations and stalks}Let $X$ be an analytic adic space. For any two points $x, y\in X$ such that $y$ is a specialization of $x$: \begin{enumerate}[(1)]\item The canonical homomorphism of stalks $\mathcal{O}_{X,y}\to \mathcal{O}_{X,x}$ is a local ring homomorphism. \item  Endow the residue fields $k(x)$ and $k(y)$ of $\mathcal{O}_{X, x}$ and $\mathcal{O}_{X, y}$ with their respective valuation topologies. The induced injective homomorphism $k(y)\hookrightarrow k(x)$ is a homeomorphism onto its image, has dense image and restricts to an injective homomorphism $k(y)^{+}\hookrightarrow k(x)^{+}$.\end{enumerate}\end{lemma}
\begin{proof}Recall that the canonical homomorphism $\varphi_{y,x}: \mathcal{O}_{X,y}\to \mathcal{O}_{X,x}$ is defined as follows: For any $f_{y}\in \mathcal{O}_{X,y}$ choose an open neighbourhood $U$ of $y$ in $X$ such that $f_{y}$ is the image in $\mathcal{O}_{X,y}$ of some $f\in \mathcal{O}_{X}(U)$. Then $U$ is also an open neighbourhood of $x$ and $\varphi_{y,x}(f_y)\in\mathcal{O}_{X,x}$ is defined to be the image $f_{x}$ of $f$ under the canonical map $\mathcal{O}_{X}(U)\to \mathcal{O}_{X,x}$. We also denote by $v_{y}$ (respectively, $v_{x}$) the valuation on $\mathcal{O}_{X,y}$ (respectively, $\mathcal{O}_{X,x}$) induced by the point $y$ (respectively, $x$).  

(1) Suppose that $f_{y}\in \mathcal{O}_{X,y}$ is an element which belongs to the maximal ideal $\mathfrak{m}_{y}$ of $\mathcal{O}_{X,y}$. We want to prove that $\varphi_{y,x}(f_{y})$ belongs to the maximal ideal $\mathfrak{m}_{x}$ of $\mathcal{O}_{X,x}$. Choose an affinoid open neighbourhood $U$ of $y$ such that $f_{y}$ is the image in $\mathcal{O}_{X,y}$ of an element $f\in \mathcal{O}_{X}(U)$. Since $x$ is a generization of $y$, the open subset $U$ is also a neighbourhood of $x$. Since $f_y\in\mathfrak{m}_{y}$, we see that $f\in \ker(y)$ when $y\in U$ is regarded as a continuous valuation on $\mathcal{O}_{X}(U)$. Since $X$, and thus also $U$, is analytic, all generizations of $y\in U$ are secondary, by \cite{Huber1}, Section~3, so $\ker(y)=\ker(x)$. It follows that $\varphi_{y,x}(f_{y})=f_{x}\in\mathcal{O}_{X,x}$ belongs to $\ker(v_{x})=\mathfrak{m}_{x}$, as desired. 

(2) Let $U$ be an affinoid open neighborhood of $y$. Then $U$ is also an affinoid open neighborhood of $x$, since $x$ is a generization of $y$. The map \begin{equation*}\Frac(\mathcal{O}_{X}(U)/\ker(x))\to k(x)\end{equation*}has dense image (see, for example, \cite{MorelNotes}, Proposition~III.6.3.1(iv)) and factors through the map $k(y)\hookrightarrow k(x)$. Hence $k(y)\hookrightarrow k(x)$ has dense image. 

To prove that $k(y)\hookrightarrow k(x)$ restricts to a map $k(y)^{+}\hookrightarrow k(x)^{+}$, recall that $k(x)^{+}$ (respectively, $k(y)^{+}$) is the subring of $k(x)$ (respectively, $k(y)$), where the valuation of $k(x)$ (respectively, $k(y)$) induced by the valuation $v_{x}$ on $\mathcal{O}_{X,x}$ (respectively, $v_{y}$ on $\mathcal{O}_{X,y}$) is $\leq1$. Hence, choosing an affinoid open neighbourhood $U$ of $y$ and an element $f\in\mathcal{O}_{X}(U)$ as before, we have to prove that the condition $v_{y}(f_{y})\leq1$ implies $v_{x}(f_{x})\leq1$. If $v_{y}(f_{y})\leq1$, we have $\vert f(y)\vert\leq1$ when $y$ is viewed as a continuous valuation on $\mathcal{O}_{X}(U)$. In other words, $y\in U(\frac{f}{1})$. But then $x$ also belongs to $U(\frac{f}{1})$, being a generization of $y$. We conclude that $v_{x}(f_{x})\leq1$, as desired. 

Finally, we verify that the map $k(y)\hookrightarrow k(x)$ is a homeomorphism onto its image. Since $k(y)^{+}$ gets mapped into $k(x)^{+}$, we see that the pre-image of the valuation ring $k(x)^{+}$ of $k(x)$ is a valuation ring of the field $k(y)$ containing $k(y)^{+}$. But for any inclusion of valuation rings $V\subseteq W$ with the same field of fractions $K\neq W$, the topology on $K$ defined by the two valuations is the same. Hence the subspace topology on $k(y)$ induced from the valuation topology on $k(x)$ coincides with the valuation topology on $k(y)$. \end{proof}
\begin{mydef}[Huber \cite{Huber3}, Definition~1.3.5]A valuation ring of an analytic adic space $X$ is a pair $(x, A)$, where $x\in X$ and $A$ is a valuation ring of the residue field $k(x)$ of $X$ at $x$ with $A\subseteq k(x)^{+}$. The point $x$ is called the support of $(x, A)$ and we sometimes write $A$ instead of $(x, A)$ if the support $x$ is understood from the context. A point $y\in X$ is called a center of $(x, A)$ if $y$ is a specialization of $x$ and the pre-image of $A$ under the homomorphism $k(y)\hookrightarrow k(x)$ induced by the ring map $\mathcal{O}_{X,y}\to \mathcal{O}_{X,x}$ is equal to $k(y)^{+}$.\end{mydef}
\begin{lemma}[Huber \cite{Huber3}, Lemma~1.3.6]\label{Existence of valuation rings}Let $X$ be an analytic adic space. \begin{enumerate}[(1)] \item For any two points $x, y\in X$ such that $y$ is a specialization of $x$ there exists a valuation ring $A$ of $X$ such that $x$ is the support and $y$ is a center of $A$. \item If $X$ is affinoid, then every valuation ring $(x, A)$ of $X$ has at most one center on $X$ and a valuation ring $(x, A)$ has a center on $X$ if and only if the image of $\mathcal{O}_{X}^{+}(X)$ under the canonical map $\mathcal{O}_{X}^{+}(X)\to k(x)^{+}$ is contained in $A$.\end{enumerate}\end{lemma}
\begin{proof}(1) By Lemma~\ref{Specializations and stalks}, the map $k(y)\hookrightarrow k(x)$ has dense image, is a homeomorphism onto its image and maps $k(y)^{+}$ into $k(x)^{+}$. Then the (topological) closure of $k(y)^{+}$ inside $k(x)$ is a valuation ring of $k(x)$ whose intersection with $k(y)$ is $k(y)^{+}$.  

(2) See the first paragraph of the proof of \cite{Huber3}, Lemma~1.3.6ii).\end{proof}
Any morphism of analytic adic spaces $f: X\to S$ induces a map $f_{v}: v(X)\to v(S)$ between the set of valuation rings of $X$ and the set of valuation rings of $S$, given by \begin{equation*}f_{v}: (x, A)\mapsto (f(x), k(f(x))\cap A).\end{equation*}If $y\in X$ is a center of $(x, A)\in v(X)$ on $X$, then $f(y)$ is a center of $f_{v}((x, A))$ on $S$. Using these notions, we can rephrase the definition~of a separated (respectively, a partially proper) morphism of analytic adic spaces as follows.
\begin{lemma}\label{Valuative criterion 1}For a morphism of analytic adic spaces $f: X\to S$ the following are equivalent: \begin{enumerate}[(1)] \item $f$ is separated. \item For every valuation ring $(x, A)$ of $X$ and every center $y\in S$ of $f_{v}((x, A))$ on $S$, there exists at most one center $z\in X$ of $(x, A)$ on $X$ with $f(z)=y$.\end{enumerate}\end{lemma}
\begin{proof}Follows from the definitions.\end{proof}
\begin{lemma}\label{Valuative criterion 2}For a quasi-separated rationalizing morphism of analytic adic spaces $f: X\to S$ the following are equivalent: \begin{enumerate}[(1)] \item $f$ is partially proper. \item For every valuation ring $(x, A)$ of $X$ and every center $y\in S$ of $f_{v}((x, A))$ on $S$, there exists precisely one center $z\in X$ of $(x, A)$ on $X$ with $f(z)=y$.\end{enumerate}\end{lemma}
\begin{proof}Follows from the definitions and Lemma~\ref{Partial properness and affinoid fields}.\end{proof}The following technical lemma is an immediate generalization of \cite{Huber3}, Lemma~1.3.14(i) and (iii). The proof given in loc.~cit.~applies verbatim in our more general situation, but we include the argument for the reader's convenience.
\begin{lemma}\label{Huber's lemma}Let $S$ be an analytic adic space, let $d: \overline{X}\to S$ and $e: Y\to S$ be adic spaces over $S$ and let $X$ be an open subspace of $\overline{X}$. Let\begin{center}\begin{tikzcd}X\arrow{r}{f} \arrow{d}{d\vert_{X}} & Y \arrow{dl}{e} \\ S \end{tikzcd}\end{center}be a morphism of adic spaces over $S$.\begin{enumerate}[(1)]\item If $e$ is separated and every point of $\overline{X}$ is a specialization of a point of $X$, then there exists at most one $S$-morphism $\overline{f}: \overline{X}\to Y$ extending $f$. \item Assume that $e$ is rationalizing, $X$ is retrocompact in $\overline{X}$ and for every open subset $W$ of $\overline{X}$ the restriction map \begin{equation*}\mathcal{O}_{\overline{X}}(W)\to \mathcal{O}_{X}(W\cap X)\end{equation*}is an isomorphism of complete topological rings. Let $\overline{x}\in \overline{X}$ and $y\in Y$ and suppose that that there exists a valuation ring $(x, A)$ of $X$ such that $\overline{x}$ is a center of $(x, A)$ on $\overline{X}$, $y$ is a center of $f_{v}((x, A))$ on $Y$, $d(\overline{x})=e(y)$ in $S$ and $x$ is the closed point of the set of generizations of $\overline{x}$ in $X$. Then there exist an affinoid open neighbourhood $W_{\overline{x}}$ of $\overline{x}$ in $\overline{X}$ and an $S$-morphism \begin{equation*}\overline{f}_{W_{\overline{x}}}: W_{\overline{x}}\to Y\end{equation*}with \begin{equation*}\overline{f}_{W_{\overline{x}}}\vert_{W_{\overline{x}}\cap X}=f\vert_{W_{\overline{x}}\cap X}\end{equation*}and $\overline{f}_{W_{\overline{x}}}(\overline{x})=y$.\end{enumerate}\end{lemma}
\begin{proof}(1) Let $\overline{f}, g: \overline{X}\to Y$ be two $S$-morphisms extending $f$. It suffices to prove that $\overline{f}$ and $g$ coincide as maps of sets and that the induced homomorphisms on stalks \begin{equation*}\overline{f}_{\overline{x}}^{\ast}, g_{\overline{x}}^{\ast}: \mathcal{O}_{Y, \overline{f}(\overline{x})}=\mathcal{O}_{Y, g(\overline{x})}\to \mathcal{O}_{\overline{X}, \overline{x}}\end{equation*}are equal for all $\overline{x}\in \overline{X}$. Given $\overline{x}\in\overline{X}$, choose a generization $x$ of $\overline{x}$ in $X$. By Lemma~\ref{Existence of valuation rings}(1), there exists a valuation ring $A$ of $\overline{X}$ with support $x$ and center $\overline{x}$. Then $\overline{f}(\overline{x})$ and $g(\overline{x})$ are centers of the valuation ring\begin{equation*}f_{v}((x, A))=\overline{f}_{v}((x, A))=g_{v}((x, A))\end{equation*}on $Y$ such that $e(\overline{f}(\overline{x}))=e(g(\overline{x}))=d(\overline{x})$. But since $e$ is separated, Lemma~\ref{Valuative criterion 1} implies that there exists at most one center of $f_{v}((x, A))$ whose image under $e$ is $d(\overline{x})$. Consequently, $\overline{f}(\overline{x})=g(\overline{x})$. Since $\overline{x}\in \overline{X}$ was arbitrary, this shows that $\overline{f}$ and $g$ coincide as maps of sets $\vert\overline{X}\vert\to \vert Y\vert$. To prove that the induced homomorphisms on stalks also coincide, consider the commutative square \begin{center}\begin{tikzcd}\mathcal{O}_{\overline{X}, \overline{x}} \arrow{r} & \mathcal{O}_{\overline{X}, x} \\ \mathcal{O}_{Y, f(\overline{x})} \arrow{u}{\overline{f}_{\overline{x}}^{\ast}} \arrow{r} & \mathcal{O}_{Y, f(x)}. \arrow{u}{f_{x}^{\ast}} \end{tikzcd}\end{center}and the analogous commutative square with $g$ in place of $\overline{f}$. Since $x$ is a generization of $\overline{x}$, the horizontal arrows are injective. It follows that $\overline{f}_{\overline{x}}^{\ast}$ (respectively, $g_{\overline{x}}^{\ast}$) is the restriction of $f_{x}^{\ast}$ along the injective map $\mathcal{O}_{Y,f(\overline{x})}\to \mathcal{O}_{Y,f(x)}$. Consequently, $\overline{f}_{\overline{x}}^{\ast}=g_{\overline{x}}^{\ast}$. 

(2) Choose an affinoid open subspace $U\subseteq S$ and an affinoid open neighbourhood $V$ of $y$ in $Y$ with $e(V)\subseteq U$ such that there exists a finite subset $E\subseteq \mathcal{O}_{Y}(V)$ with the property that $\mathcal{O}_{Y}^{+}(V)$ is the smallest ring of integral elements of $\mathcal{O}_{Y}(V)$ containing $E$ and the image of $\mathcal{O}_{S}^{+}(U)$ in $\mathcal{O}_{Y}(V)$. Since $y$ is a center of $f_{v}((x, A))$, it is a specialization of the support $f(x)$ of $f_{v}((x, A))$. Since $y\in V$ and $V$ is open, this implies that $x\in f^{-1}(V)$. On the other hand, $x$ is the closed point of the set of generizations of $\overline{x}$ in $X$. Hence, letting $W_{0}$ be a fixed affinoid open neighborhood of $\overline{x}$ in $d^{-1}(U)$ and letting $V'$ be some quasi-compact open neighbourhood of $x$ in $f^{-1}(V)\cap W_{0}$, we see that every generization of $\overline{x}$ in $X$ lies in $V'$. Let $G_{\overline{x}}$ be the set of all generizations of $\overline{x}$ in $\overline{X}$. Write \begin{equation*}G_{\overline{x}}=\bigcap_{W\in \mathcal{U}_{\overline{x}}}W,\end{equation*}where $\mathcal{U}_{\overline{x}}$ is the family of all quasi-compact open neighbourhoods of $\overline{x}$ inside $W_{0}$. Since $X$ is retrocompact in $\overline{X}$ by hypothesis, $W\cap X$ is a quasi-compact open subset of $W_{0}$ for every $W\in\mathcal{U}_{\overline{x}}$. By the same token, $V'\cap X$ is a quasi-compact subset of $W_{0}$, so the complement $C$ of $V'\cap X$ in $W_{0}$ is constructible in $W_{0}$. Suppose that there exists no $W\in\mathcal{U}_{\overline{x}}$ such that $W\cap X$ is contained in $V'$. Then $W\cap X\cap C\neq \emptyset$ for all $W\in\mathcal{U}_{\overline{x}}$ and then\begin{equation*}(W\cap X\cap C)_{W\in\mathcal{U}_{\overline{x}}}\end{equation*}is a family of constructible subsets of $W_{0}$ which has the finite intersection property. By compactness of the constructible topology on the spectral space $W_{0}$, it follows that \begin{equation*}G_{\overline{x}}\cap X\cap C=(\bigcap_{W\in\mathcal{U}_{\overline{x}}}W)\cap X\cap C=\bigcap_{W\in\mathcal{U}_{\overline{x}}}(W\cap X\cap C)\end{equation*}is non-empty, which contradicts the assumption that $G_{\overline{x}}\cap X\subseteq V'$. This shows that there exists an affinoid open neighbourhood $W$ of $\overline{x}$ in $d^{-1}(U)$ such that \begin{equation*}W\cap X\subseteq V'\subseteq f^{-1}(V).\end{equation*}In this way we obtain a continuous ring map $\mathcal{O}_{Y}(V)\to \mathcal{O}_{X}(W\cap X)$ induced by the morphism $f\vert_{W\cap X}: W\cap X\to V$. Composing this map $\mathcal{O}_{Y}(V)\to \mathcal{O}_{X}(W\cap X)$ with the inverse of $\mathcal{O}_{\overline{X}}(W)\tilde{\to}\mathcal{O}_{X}(W\cap X)$ yields a continuous ring map $\varphi: \mathcal{O}_{Y}(V)\to \mathcal{O}_{\overline{X}}(W)$. Now set \begin{equation*}W_{\overline{x}}=W(\frac{\varphi(E)}{1})\end{equation*}and let $\psi$ be the composition of $\varphi$ with the restriction map $\mathcal{O}_{\overline{X}}(W)\to \mathcal{O}_{\overline{X}}(W_{\overline{x}})$. By the assumption that $\mathcal{O}_{Y}^{+}(V)$ is the smallest open and integrally closed subring of $\mathcal{O}_{Y}(V)$ which contains $E$ and the image of $\mathcal{O}_{S}^{+}(U)$, we see that \begin{equation*}\psi(\mathcal{O}_{Y}^{+}(V))\subseteq \mathcal{O}_{\overline{X}}^{+}(W_{\overline{x}}).\end{equation*}Consequently, $\psi$ induces a morphism of adic spaces over $U$ \begin{equation*}\overline{f}_{W_{\overline{x}}}: W_{\overline{x}}\to Y.\end{equation*}By construction, we have \begin{equation*}\overline{f}_{W_{\overline{x}}}\vert_{W_{\overline{x}}\cap X}=f\vert_{W_{\overline{x}}\cap X}.\end{equation*}Since $\overline{x}$ is a center of $(x, A)$ on $\overline{X}$ (and hence also on its subspace $W_{\overline{x}}$), the image $\overline{f}_{W_{\overline{x}}}(\overline{x})$ is a center of $\overline{f}_{W_{\overline{x}}v}((x, A))$ on $V$. But $\overline{f}_{W_{\overline{x}}}(x)=f(x)$ since $x\in W_{\overline{x}}\cap X$, so \begin{equation*}\overline{f}_{W_{\overline{x}}v}((x, A))=f_{v}((x, A)).\end{equation*}But, by assumption, $y$ is a center of $f_{v}((x, A))$ on $Y$ (and hence also on its subspace $V$). By Lemma~\ref{Existence of valuation rings}(2) the center of $f_{v}((x, A))$ on the affinoid analytic adic space $V$ is unique. It follows that $\overline{f}_{W_{\overline{x}}}(\overline{x})=y$. This concludes the proof.\end{proof}
Using the above lemma, we also get the following generalization of \cite{Huber3}, Lemma~5.1.7. Again, the proof follows the argument given in loc.~cit.
\begin{lemma}\label{Huber's Lemma 2}Let \begin{center}\begin{tikzcd}X\arrow[hook]{r}{j} \arrow{d}{f} & \overline{X} \arrow{dl}{\overline{f}} \\ S\end{tikzcd}\end{center}be a commutative diagram of analytic adic spaces such that $\overline{f}$ is partially proper and rationalizing and $j$ is a quasi-compact open embedding such that every point of $\overline{X}$ is a specialization of a point of $j(X)$ and $\mathcal{O}_{\overline{X}}\to j_{\ast}\mathcal{O}_{X}$ is an isomorphism of sheaves of complete topological rings. Then, for every open subset $U\subseteq S$, the corresponding diagram \begin{center}\begin{tikzcd}f^{-1}(U)\arrow[hook]{r} \arrow{d} & \overline{f}^{-1}(U) \arrow{dl} \\ U \end{tikzcd}\end{center}is a universal compactification of $f\vert_{f^{-1}(U)}: f^{-1}(U)\to U$. \end{lemma}
\begin{proof}Let $g: Y\to U$ be a partially proper rationalizing morphism of analytic adic spaces (in particular, $g$ is separated and rationalizing), and let $h: f^{-1}(U)\to Y$ be a morphism of adic spaces over $U$. We have to show that there exists a unique morphism of adic spaces over $U$ \begin{equation*}i: \overline{f}^{-1}(U)\to Y\end{equation*}extending $h$. Uniqueness is an immediate consequence of Lemma~\ref{Huber's lemma}(1), since $g$ is separated. For every $\overline{x}\in \overline{f}^{-1}(U)$ let $x\in f^{-1}(U)$ be the closed point of the set of generizations of $\overline{x}$ in $f^{-1}(U)$. By Lemma~\ref{Existence of valuation rings}(1) there exists a valuation ring $(x, A)$ of $\overline{f}^{-1}(U)$ with support $x$ and center $\overline{x}$. Then $\overline{f}(\overline{x})$ is a center of $\overline{f}_{v}((x, A))$ on $U$. Since $f$ is the restriction of $\overline{f}$, we have $\overline{f}_{v}((x, A))=f_{v}((x, A))$. As $h$ is a $U$-morphism, we see that \begin{equation*}f_{v}((x, A))=g_{v}(h_{v}((x, A))).\end{equation*}Applying Lemma~\ref{Valuative criterion 2} to the partially proper rationalizing morphism $g$, we find a center $y$ of $h_{v}((x, A))$ on $Y$ such that $g(y)$ is equal to the center $\overline{f}(\overline{x})$ of $g_{v}(h_{v}((x, A)))$ on $U$. Now, using Lemma~\ref{Huber's lemma}(2), we obtain an open neighbourhood $W_{\overline{x}}$ of $\overline{x}$ in $\overline{f}^{-1}(U)$ and a $U$-morphism \begin{equation*}i_{W_{\overline{x}}}: W_{\overline{x}}\to Y\end{equation*}which extends $h\vert_{W_{\overline{x}}\cap f^{-1}(U)}$. By the uniqueness assertion of Lemma~\ref{Huber's lemma}(1), the morphisms $i_{W_{\overline{x}_{1}}}$ and $i_{W_{\overline{x}_{2}}}$ coincide on the intersection $W_{\overline{x}_{1}}\cap W_{\overline{x}_{2}}$ for any two points $\overline{x}_{1}, \overline{x}_{2}\in \overline{f}^{-1}(U)$. Hence the $U$-morphisms $i_{W_{\overline{x}}}: W_{\overline{x}}\to Y$ for $\overline{x}\in \overline{f}^{-1}(U)$ glue to a unique $U$-morphism $i: \overline{f}^{-1}(U)\to Y$ with $i\vert_{f^{-1}(U)}=h$, as desired.\end{proof}
We are now ready to apply our results to the problem of constructing compactifications of locally affinoid morphisms. Let $f: X\to S$ be a locally affinoid morphism of Tate adic spaces. By Proposition~\ref{Direct images 2} and Proposition~\ref{Properness and adapted subsets}, $(f_{\ast}\mathcal{O}_{X}, f_{\ast}\mathcal{O}_{X}^{+})$ and $(f_{\ast}\mathcal{O}_{X}, \overline{(f_{\ast}\mathcal{O}_{X})^{+}})$ are quasi-coherent topological adic Banach algebras on $S$. Since for every $(f_{\ast}\mathcal{O}_{X}, f_{\ast}\mathcal{O}_{X}^{+})$-adapted affinoid open subset $U\subseteq S$ the subring $\overline{(f_{\ast}\mathcal{O}_{X})^{+}}(U)$ is the smallest ring of integral elements of $(f_{\ast}\mathcal{O}_{X})(U)=\mathcal{O}_{X}(f^{-1}(U))$ containing the image of $\mathcal{O}_{S}^{+}(U)$, there is a canonical inclusion map \begin{equation*}\overline{(f_{\ast}\mathcal{O}_{X})^{+}}(U)\hookrightarrow f_{\ast}\mathcal{O}_{X}^{+}(U)\end{equation*}for every such $U$. These inclusion maps define a morphism of quasi-coherent topological adic Banach $(\mathcal{O}_{S}, \mathcal{O}_{S}^{+})$-algebras \begin{equation*}(f_{\ast}\mathcal{O}_{X}, \overline{(f_{\ast}\mathcal{O}_{X})^{+}})\hookrightarrow (f_{\ast}\mathcal{O}_{X}, f_{\ast}\mathcal{O}_{X}^{+})\end{equation*}which, by Theorem~\ref{The affinoid equivalence}, gives rise to a triangle of morphisms of adic spaces \begin{center}\begin{tikzcd}X=\underline{\Spa}_{S}(f_{\ast}\mathcal{O}_{X}, f_{\ast}\mathcal{O}_{X}^{+})\arrow{d}{f} \arrow{r}{j} & \underline{\Spa}_{S}(f_{\ast}\mathcal{O}_{X}, \overline{(f_{\ast}\mathcal{O}_{X})^{+}})=\overline{X}\arrow{dl}{\overline{f}} \\ S, \end{tikzcd}\end{center}where $\overline{f}$ is the canonical morphism from the relative adic spectrum to the base $S$.        
\begin{thm}\label{Compactifications}Let $f: X\to S$ be a locally affinoid morphism of Tate adic spaces. Then in the above commutative triangle of analytic adic spaces\begin{center}\begin{tikzcd}X\arrow{r}{j} \arrow{d}{f} & \overline{X}=\underline{\Spa}_{S}(f_{\ast}\mathcal{O}_{X}, \overline{(f_{\ast}\mathcal{O}_{X})^{+}}) \arrow{dl}{\overline{f}} \\ S \end{tikzcd}\end{center}the morphism $\overline{f}$ is a locally affinoid proper morphism, $j$ is a locally affinoid injective morphism such that the induced map of complete topological rings $\mathcal{O}_{\overline{X}}\to j_{\ast}\mathcal{O}_{X}$ is an isomorphism and every point of $\overline{X}$ is a specialization of a point of $j(X)$. Moreover, the triangle $(\overline{X}, \overline{f}, j)$ enjoys the following universal property: For any commutative triangle of morphisms of analytic adic spaces \begin{center}\begin{tikzcd}X \arrow{r}{h} \arrow{d}{f} & Y \arrow{dl}{g} \\ S \end{tikzcd}\end{center}with $g$ locally affinoid and proper, there exists a unique morphism $i: \overline{X}\to Y$ making the diagram \begin{center}\begin{tikzcd}X\arrow{dd}{f} \arrow[hook]{dr}{j} \arrow{drr}{h} & & \\ & \overline{X} \arrow{r}{i} \arrow{dl}{\overline{f}} & Y \arrow{dll}{g} \\ S \end{tikzcd}\end{center}commute. The triangle $(\overline{X}, \overline{f}, j)$ is a universal compactification of $f$ if and only if $j$ is an open immersion if and only if $f$ is rationalizing.\end{thm}  
\begin{proof} By Lemma~\ref{Relative adic spectrum is partially proper}, the morphism $\overline{f}: \overline{X}\to S$ is proper. Since the morphism of quasi-coherent topological adic Banach algebras $(f_{\ast}\mathcal{O}_{X}, \overline{(f_{\ast}\mathcal{O}_{X})^{+}})\hookrightarrow (f_{\ast}\mathcal{O}_{X}, f_{\ast}\mathcal{O}_{X}^{+})$ is the identity on the underlying topological Banach $\mathcal{O}_{S}$-algebras, the morphism of complete topological rings $\mathcal{O}_{\overline{X}}\to j_{\ast}\mathcal{O}_{X}$ given by $j$ is an isomorphism. In particular, $j$ is injective. Moreover, by Corollary~\ref{Locally affinoid 4}, $j$ is necessarily locally affinoid and, in particular, quasi-compact. To prove that every point of $\overline{X}$ is a specialization of a point in $j(X)$, consider an $(f_{\ast}\mathcal{O}_{X}, f_{\ast}\mathcal{O}_{X}^{+})$-adapted affinoid open subspace $U\subseteq S$. If \begin{equation*}x\in \overline{f}^{-1}(U)=\Spa(\mathcal{O}_{X}(f^{-1}(U)), \overline{\mathcal{O}_{S}^{+}(U)+\mathcal{O}_{X}(f^{-1}(U))^{\circ\circ}})\end{equation*}is a rank $1$ valuation on $\mathcal{O}_{X}(f^{-1}(U))$, then $\vert e(x)\vert\leq1$ for all $e\in \mathcal{O}_{X}(f^{-1}(U))^{\circ}$ and, in particular, $\vert e(x)\vert\leq1$ for all $e\in \mathcal{O}_{X}^{+}(f^{-1}(U))$. In other words, every rank $1$ point of $\overline{f}^{-1}(U)$ actually belongs to $f^{-1}(U)=\Spa(\mathcal{O}_{X}(U), \mathcal{O}_{X}^{+}(U))$. But by \cite{Huber3}, Lemma~1.1.10(ii) (whose proof applies without the locally Noetherian hypothesis used in \cite{Huber3}), every point of the analytic adic space $\overline{f}^{-1}(U)$ has a rank $1$ generization. It follows that every point of $\overline{f}^{-1}(U)$ is a specialization of a point of $j(f^{-1}(U))$. Covering $S$ by $(f_{\ast}\mathcal{O}_{X}, f_{\ast}\mathcal{O}_{X}^{+})$-adapted affinoid open subspaces $U$ we conclude that every point of $\overline{X}$ is a specialization of a point of $j(X)$, as claimed. 

Now let \begin{center}\begin{tikzcd}X\arrow{d}{f} \arrow{r}{h} & Y \arrow{dl}{g} \\ S \end{tikzcd}\end{center}be a morphism of adic spaces over $S$ with $g: Y\to S$ locally affinoid and proper. Since $g$ is locally affinoid, $(g_{\ast}\mathcal{O}_{Y}, g_{\ast}\mathcal{O}_{Y}^{+})$ is a quasi-coherent topological adic Banach $(\mathcal{O}_{S}, \mathcal{O}_{S}^{+})$-algebra, by Proposition~\ref{Direct images 2}, and, by Theorem~\ref{The affinoid equivalence}, the morphism $h: X\to Y$ of adic spaces over $S$ induces a morphism of quasi-coherent topological adic Banach $(\mathcal{O}_{S}, \mathcal{O}_{S}^{+})$-algebras \begin{equation*}(g_{\ast}\mathcal{O}_{Y}, g_{\ast}\mathcal{O}_{Y}^{+})\to (f_{\ast}\mathcal{O}_{X}, f_{\ast}\mathcal{O}_{X}^{+}).\end{equation*}The underlying morphism of sheafy quasi-coherent topological Banach $\mathcal{O}_{S}$-algebras $g_{\ast}\mathcal{O}_{Y}\to f_{\ast}\mathcal{O}_{X}$ then induces a morphism of adic spaces over $S$ \begin{equation*}i: \overline{X}=\underline{\Spa}_{S}(f_{\ast}\mathcal{O}_{X}, \overline{(f_{\ast}\mathcal{O}_{X})^{+}})\to \underline{\Spa}_{S}(g_{\ast}\mathcal{O}_{Y}, \overline{(g_{\ast}\mathcal{O}_{Y})^{+}})=Y\end{equation*}via Theorem~\ref{The affinoid equivalence 2}. By construction, the diagram of quasi-coherent sheaves of topological adic Banach $(\mathcal{O}_{S}, \mathcal{O}_{S}^{+})$-algebras \begin{center}\begin{tikzcd}(g_{\ast}\mathcal{O}_{Y}, g_{\ast}\mathcal{O}_{Y}^{+})=(g_{\ast}\mathcal{O}_{Y}, \overline{(g_{\ast}\mathcal{O}_{Y}^{+})}) \arrow{r} \arrow{dr} & (f_{\ast}\mathcal{O}_{X}, \overline{(f_{\ast}\mathcal{O}_{X}^{+})}) \arrow[hook]{d} \\ & (f_{\ast}\mathcal{O}_{X}, f_{\ast}\mathcal{O}_{X}^{+}) \end{tikzcd}\end{center}is commutative. By Theorem~\ref{The affinoid equivalence 2}, this means precisely that the morphism $i: \overline{X}\to Y$ gives rise to a commutative diagram of adic spaces over $S$ as in the statement of the theorem. 

If $j$ is an open immersion, $f$ must be rationalizing by Lemma~\ref{Composition of rationalizing morphisms}. Now assume that $f$ is rationalizing. By Lemma~\ref{Huber's Lemma 2}, to prove that the triangle $(\overline{X}, \overline{f}, j)$ constructed above is a universal compactification, it suffices to show that the quasi-compact injective morphism $j$ is an open immersion. Covering $S$ by $(f_{\ast}\mathcal{O}_{X}, f_{\ast}\mathcal{O}_{X}^{+})$-adapted affinoid open subspaces, we reduce the problem to the case when $S$ is affinoid and $(f_{\ast}\mathcal{O}_{X}, f_{\ast}\mathcal{O}_{X}^{+})$-adapted. By Corollary~\ref{Locally affinoid 3}, this means that we may assume that $f: X\to S$ is a morphism of affinoid analytic adic spaces. By Lemma~\ref{Rationalizing morphisms and rational subsets}, there is an affinoid open cover $(W_{k})_{k}$ of $X$ such that the morphism of Huber pairs $(\mathcal{O}_{S}(S), \mathcal{O}_{S}^{+}(S))\to (\mathcal{O}_{X}(W_{k}), \mathcal{O}_{X}^{+}(W_{k}))$ is rationalizing, for every $k$. Possibly replacing the cover $(W_{k})_{k}$ by a refinement and using Lemma~\ref{First Lemma on rationalizing morphisms}(1) along with Example~\ref{Open immersions}, we may assume that $(W_{k})_{k}$ is a finite rational cover of the affinoid adic space $X$. For every $k$, choose a finite subset $E_{k}\subseteq \mathcal{O}_{X}^{+}(W_{k})$ such that $\mathcal{O}_{X}^{+}(W_{k})$ is integral over its $\mathcal{O}_{S}^{+}(S)$-subalgebra generated by $E_{k}\cup\mathcal{O}_{X}(W_{k})^{\circ\circ}$ and choose elements \begin{equation*}f_{1,k},\dots, f_{n_{k},k}, g_{k}\in \mathcal{O}_{X}(X)=\mathcal{O}_{\overline{X}}(\overline{X})\end{equation*}generating the unit ideal of $\mathcal{O}_{X}(X)$ such that \begin{equation*}W_{k}=X(\frac{f_{1,k},\dots,f_{n_{k},k}}{g_{k}}).\end{equation*}Consider the rational open subset \begin{equation*}V_{k}=\overline{X}(\frac{f_{1,k},\dots,f_{n_{k},k}}{g_{k}})\end{equation*}of $\overline{X}$. Since \begin{equation*}\frac{f_{i, k}}{g_{k}}\in \mathcal{O}_{X}^{+}(W_{k})\end{equation*}for all $i$, the morphism of Huber pairs \begin{equation*}(\mathcal{O}_{X}(X), \overline{\mathcal{O}_{S}^{+}(S)+\mathcal{O}_{X}(X)^{\circ\circ}})\to (\mathcal{O}_{X}(W_{k}), \mathcal{O}_{X}^{+}(W_{k}))\end{equation*}factors through \begin{equation*}(\mathcal{O}_{X}(W_{k}), \overline{\mathcal{O}_{S}^{+}(S)+\mathcal{O}_{X}(X)^{\circ\circ}}\langle \frac{f_{1, k},\dots, f_{n_{k}, k}}{g_{k}}\rangle).\end{equation*}This shows that the image $j(W_{k})$ of $W_{k}$ in $\overline{X}$ is contained in $V_{k}$. On the other hand, $\overline{\mathcal{O}_{S}^{+}(S)+\mathcal{O}_{X}(W_{k})^{\circ\circ}}$ is the smallest subring of integral elements of $\mathcal{O}_{X}(W_{k})$ containing $\mathcal{O}_{S}^{+}(S)$, so \begin{equation*}\overline{\mathcal{O}_{S}^{+}(S)+\mathcal{O}_{X}(W_{k})^{\circ\circ}}\subseteq \overline{\mathcal{O}_{S}^{+}(S)+\mathcal{O}_{X}(X)^{\circ\circ}}\langle \frac{f_{1, k},\dots, f_{n_{k}, k}}{g_{k}}\rangle.\end{equation*}Therefore, if $x\in V_{k}$ satisfies $\vert e(x)\vert\leq1$ for all $e\in E_{k}$, then $\vert f(x)\vert\leq1$ for all \begin{equation*}f\in E_{k}\cup\mathcal{O}_{X}(W_{k})^{\circ\circ}.\end{equation*}It follows that the morphism $j\vert_{W_{k}}: W_{k}\hookrightarrow V_{k}$ identifies $W_{k}$ with the rational subset \begin{equation*}V_{k}(\frac{E_{k}}{1})=\{\, x\in V_{k}\mid \vert e(x)\vert\leq 1 \ \text{for all}~e\in E_{k}\,\}\end{equation*}of $V_{k}$. In particular, $j\vert_{W_{k}}: W_{k}\hookrightarrow V_{k}$ is an open immersion for each $k$. Since $(W_{k})_{k}$ is an open cover of $X$, this entails that $j: X\hookrightarrow \overline{X}$ is an open immersion.\end{proof}
\begin{rmk}By examining the proof of \cite{Huber3}, Theorem~5.1.5, one sees that every morphism of analytic adic spaces $f: X\to S$ which is rationalizing, separated (in the sense of Definition~\ref{Definition of separated morphisms}) and taut admits a universal compactification as in Definition~\ref{Definition of universal compactifications}, without locally Noetherian assumptions. In fact, such a non-Noetherian generalization of Huber's theorem (still for morphisms of ${}^{+}$weakly finite type, but for many adic spaces which need not even be analytic) was proved recently by Solodov~\cite{Solodov25}. 

Hence the main content of Theorem~\ref{Compactifications} is the explicit description of the universal compactification (when $S$ is Tate and $f$ is locally affinoid) as the relative adic spectrum of the quasi-coherent topological adic Banach algebra $(f_{\ast}\mathcal{O}_{X}, \overline{(f_{\ast}\mathcal{O}_{X})^{+}})$. 

We note that this description is new even in the case of morphisms ${}^{+}$weakly of finite type between locally Noetherian adic spaces, the setting of Huber's original result. The proof of Theorem~\ref{Compactifications} also provides a simpler, more streamlined argument for the existence of the universal compactification in the case of locally affinoid morphisms. \end{rmk}
\begin{rmk}\label{Caution on compactifications}In the context of Theorem~\ref{Compactifications}, let us emphasize that the compactification $\overline{f}$ of a finite type locally affinoid morphism $f: X\to S$ need not be of finite type. Indeed, if $\overline{f}$ was of finite type, then, at least if $X$ and $S$ are uniform or locally Noetherian, $\overline{f}$ would have to be finite by Proposition~\ref{Proper locally affinoid morphisms of finite type} and Remark~\ref{Proper locally affinoid morphisms of finite type 2}, respectively. That this need not be the case can already be seen in the simplest examples: The universal compactification \begin{equation*}\overline{\mathbb{B}_{K}^{n}}=\Spa(K\langle T_1,\dots, T_n\rangle, K^{\circ}+K^{\circ\circ}\langle T_1,\dots, T_n\rangle)\end{equation*}of the $n$-dimensional closed unit ball \begin{equation*}\mathbb{B}_{K}^{n}=\Spa(K\langle T_1,\dots, T_n\rangle, K^{\circ}\langle T_1,\dots, T_n\rangle)\end{equation*}over a nonarchimedean field $K$ is far from being finite over its base $\Spa(K, K^{\circ})$. \end{rmk}

\section{Relation to morphisms of formal schemes}

In this section we relate the notion of quasi-coherent sheaves introduced in this paper to the notion of adically quasi-coherent sheaves on formal schemes introduced by Fujiwara and Kato in \cite{FK}. This will allow us to relate locally affinoid morphisms of Tate adic spaces to affine morphisms of formal schemes. Throughout this section, we use the terminology on formal schemes from our previous work \cite{Dine25-2}, which we recall below. 

In what follows, we fix a pair $(R, \varpi)$, where $R$ is a complete sheafy adic ring and $\varpi\in R$ is a non-zero-divisor and non-unit generating an ideal of definition of $R$. We view $R[\varpi^{-1}]$ as a Tate ring with pair of definition $(R, \varpi)$ and we denote by $\overline{R}$ the integral closure of $R$ inside $R[\varpi^{-1}]$. We restrict our attention to formal schemes $\mathfrak{X}$ equipped with an adic morphism $\mathfrak{X}\to \Spf(R)$, i.e., for any affine open subset $\mathfrak{U}\subseteq \mathfrak{X}$ the topology on $\mathcal{O}_{\mathfrak{X}}(\mathfrak{U})$ is the $\varpi$-adic topology. As in \cite{Dine25-2}, we call such formal schemes \textit{adic formal $R$-schemes}. For an adic formal $R$-scheme $\mathfrak{X}$ and an integer $k\geq0$, we denote by $\mathfrak{X}_{k}$ the scheme whose underlying topological space is $\vert\mathfrak{X}\vert$ and whose structure sheaf is $\mathcal{O}_{\mathfrak{X}}/\varpi^{k+1}\mathcal{O}_{\mathfrak{X}}$. In particular, $\mathfrak{X}_{0}$ denotes the special fiber of $\mathfrak{X}$. 

We call an adic formal $R$-scheme $\mathfrak{X}$ \textit{locally rig-sheafy} (respectively, \textit{locally stably uniform}) if it can be covered by affine open subsets $\mathfrak{U}=\Spf(A)$ (where $A=\mathcal{O}_{\mathfrak{X}}(\mathfrak{U})$, a complete adic ring) such that the Tate ring $A[\varpi^{-1}]$, with pair of definition $(A, \varpi)$, is sheafy (respectively, stably uniform)\footnote{We call such adic rings and affine open subsets \textit{rig-sheafy} (respectively, \textit{stably uniform}).}. By the famous sheafiness result of Buzzard-Verberkmoes \cite{BV} and Mihara \cite{Mihara}, the second condition implies the first. The locally stably uniform condition is equivalent to $\mathcal{O}_{\mathfrak{X}}(\mathfrak{U})$ being stably uniform for every affine open $\mathfrak{U}\subseteq \mathfrak{X}$ (\cite{Dine25-2}, Lemma~5.15). For locally rig-sheafy adic formal $R$-schemes $\mathfrak{X}$, there exists a natural analytic adic space $\mathfrak{X}_{\eta}^{\ad}$ over $\Spa(R[\varpi^{-1}], \overline{R})$, called the \textit{(adic analytic) generic fiber} of $\mathfrak{X}$, which is functorial in $\mathfrak{X}$ and is equipped with a continuous and surjective specialization map \begin{equation*}\spc_{\mathfrak{X}_{\eta}^{\ad},\mathfrak{X}}: \vert\mathfrak{X}_{\eta}^{\ad}\vert\to \vert\mathfrak{X}\vert\end{equation*}such that the pre-image of any rig-sheafy affine open subset of $\mathfrak{X}$ is affinoid (\cite{Dine25-2}, Lemma~3.8 and Definition~3.9). For a morphism of locally rig-sheafy adic formal $R$-schemes \begin{equation*}f_{0}: \mathfrak{X}\to \mathfrak{S},\end{equation*}we denote by $f_{0\eta}$ the corresponding adic space morphism between their generic fibers. In what follows, we always view $\mathfrak{X}_{\eta}^{\ad}$ as a Tate adic space by equipping it with the global pseudo-uniformizer $\varpi$, but we usually suppress $\varpi$ from our notation. 

Finally, for an adic space $X$ over $\Spa(R[\varpi^{-1}], \overline{R})$, we call a locally rig-sheafy adic formal $R$-scheme $\mathfrak{X}$ a \textit{formal $R$-model} of $X$ if $X\cong\mathfrak{X}_{\eta}^{\ad}$, and for a morphism $f$ of adic spaces over $\Spa(R[\varpi^{-1}], \overline{R})$, we call a morphism of locally rig-sheafy adic formal $R$-schemes $f_{0}$ a formal $R$-model of $f$ is $f=f_{0\eta}$. 

\subsection{Locally affinoid morphisms and affine morphisms} 

We recall the following special case of \cite{FK}, Ch.~I, Definition~3.1.3.
\begin{mydef}[Fujiwara-Kato]\label{Adically quasi-coherent sheaf}Let $\mathfrak{X}$ be an adic formal $R$-scheme. A sheaf of $\mathcal{O}_{\mathfrak{X}}$-modules $\mathcal{F}$ is called \textit{adically quasi-coherent}~if it is complete, i.e.,\begin{equation*}\mathcal{F}=\varprojlim_{k}\mathcal{F}/\varpi^{k}\mathcal{F},\end{equation*}and if the sheaf of $\mathcal{O}_{\mathfrak{X}}/\varpi^{k+1}\mathcal{O}_{\mathfrak{X}}$-modules $\mathcal{F}/\varpi^{k+1}\mathcal{F}$ is a quasi-coherent sheaf on the scheme $\mathfrak{X}_{k}$ for all integers $k\geq 0$.\end{mydef}
Our goal in this section is to relate quasi-coherent sheaves of adic Banach algebras on a uniform adic space $X$ over $\Spa(R[\varpi^{-1}], \overline{R})$ to adically quasi-coherent sheaves on formal $R$-models of $X$. 
\begin{mydef}[Rig-sheafy and uniform adically quasi-coherent algebras]An adically quasi-coherent sheaf of algebras $\mathcal{A}_{0}$ on a locally rig-sheafy $\varpi$-torsion-free adic formal $R$-scheme $\mathfrak{S}$ is called \textit{rig-sheafy} (respectively, \textit{uniform}) if for every rig-sheafy affine open subset $\mathfrak{U}$ of $\mathfrak{S}$, the complete adic ring $\mathcal{A}_{0}(\mathfrak{U})$ is again rig-sheafy (respectively, stably uniform) and $\varpi$-torsion-free. In other words, $\mathcal{A}_{0}(\mathfrak{U})$ is $\varpi$-torsion-free and the Tate ring $\mathcal{A}_{0}(\mathfrak{U})[\varpi^{-1}]$ with pair of definition $(\mathcal{A}_{0}(\mathfrak{U}), \varpi)$ is sheafy (respectively, stably uniform). \end{mydef}
Recall the notion of relative formal spectrum $\Spf(\mathcal{A}_{0})$ of an adically quasi-coherent sheaf of algebras on an adic formal $R$-scheme $\mathfrak{S}$, which was introduced by Fujiwara-Kato in \cite{FK}, Ch.~I, \S4.1. 
\begin{lemma}\label{Formal spectrum of a rig-sheafy algebra}Let $\mathfrak{S}$ be a locally rig-sheafy $\varpi$-torsion-free adic formal $R$-scheme and let $\mathcal{A}_{0}$ be an adically quasi-coherent sheaf of algebras on $\mathfrak{S}$. If $\mathcal{A}_{0}$ is rig-sheafy, then the relative formal spectrum $\Spf(\mathcal{A}_{0})$ of $\mathcal{A}_{0}$ is a locally rig-sheafy $\varpi$-torsion-free adic formal $R$-scheme. Moreover, $\mathcal{A}_{0}$ is uniform if and only if $\Spf(\mathcal{A}_{0})$ is a locally stably uniform $\varpi$-torsion-free adic formal $R$-scheme. \end{lemma}
\begin{proof}The relative formal spectrum $\Spf(\mathcal{A}_{0})$ is covered by affine open subsets of the form $\Spf(\mathcal{A}_{0}(\mathfrak{U}))$ where $\mathfrak{U}$ is an affine open subset of $\mathfrak{S}$. Hence the assertions of the lemma follow from the definitions and \cite{Dine25-2}, Lemma~5.15, which says that every affine open subset of a locally stably uniform adic formal $R$-scheme is stably uniform. \end{proof}
\begin{lemma}\label{Affine morphisms and locally affinoid morphisms}Let $\mathfrak{S}$ be a locally rig-sheafy $\varpi$-torsion-free adic formal $R$-scheme and let $\mathcal{A}_{0}$ be a rig-sheafy adically quasi-coherent sheaf of algebras on $\mathfrak{S}$. Consider the affine morphism of locally rig-sheafy adic formal $R$-schemes \begin{equation*}f_{0}: \mathfrak{X}=\Spf(\mathcal{A}_{0})\to \mathfrak{S}\end{equation*}and let \begin{equation*}f=f_{0\eta}: X\to S\end{equation*}be the induced morphism between their generic fibers. Then $f$ is locally affinoid. More precisely, for every rig-sheafy affine open subset $\mathfrak{U}\subseteq \mathfrak{S}$, the pre-image $U=\spc_{S,\mathfrak{S}}^{-1}(\mathfrak{U})$ of $\mathfrak{U}$ in $S$ is an affinoid open subspace such that $f^{-1}(U)$ is affinoid.\end{lemma}
\begin{proof}Consider the commutative diagram of topological spaces and continuous maps \begin{center}\begin{tikzcd}X\arrow{r}{f} \arrow{d}{\spc_{X,\mathfrak{X}}} & S \arrow{d}{\spc_{S,\mathfrak{S}}} \\ \mathfrak{X}\arrow{r}{f_{0}} & \mathfrak{S}.\end{tikzcd}\end{center}By \cite{Dine25-2}, Lemma~3.8, pre-images of rig-sheafy affine open subsets under the vertical arrows in this diagram are affinoid. On the other hand, since $\mathcal{A}_{0}$ is rig-sheafy, the pre-image \begin{equation*}f_{0}^{-1}(\mathfrak{U})=\Spf(\mathcal{A}_{0}(\mathfrak{U}))\end{equation*}of any rig-sheafy affine open subset $\mathfrak{U}\subseteq \mathfrak{S}$ is rig-sheafy. The assertion follows from this and commutativity of the diagram.  \end{proof}
\begin{prop}\label{Relative formal spectrum and relative adic spectrum}Let $\mathfrak{S}$ be a locally rig-sheafy $\varpi$-torsion-free adic formal $R$-scheme with generic fiber $S=\mathfrak{S}_{\eta}^{\ad}$. Let $\mathcal{A}_{0}$ be a rig-sheafy adically quasi-coherent sheaf of algebras on $\mathfrak{S}$. Then there exists a quasi-coherent sheaf of adic Banach algebras $(\mathcal{A}, \mathcal{A}^{+})$ on $S$ such that for every rig-sheafy affine open $\mathfrak{U}\subseteq \mathfrak{S}$ with generic fiber $U\subseteq S$, we have \begin{equation*}\mathcal{A}(U)=\mathcal{A}_{0}(\mathfrak{U})[\varpi^{-1}],\end{equation*}with pair of definition $(\mathcal{A}_{0}(\mathfrak{U}), \varpi)$, and \begin{equation*}\mathcal{A}^{+}(U)=\overline{\mathcal{A}_{0}(\mathfrak{U})},\end{equation*}the integral closure of $\mathcal{A}_{0}(\mathfrak{U})$ inside $\mathcal{A}(U)$. Moreover, if \begin{equation*}f_{0}: \mathfrak{X}=\Spf(\mathcal{A}_{0})\to \mathfrak{S}\end{equation*}is the canonical morphism from the relative formal spectrum of $\mathcal{A}_{0}$, then $f_{0\eta}$ is the canonical morphism \begin{equation*}X=\underline{\Spa}_{S}(\mathcal{A}, \mathcal{A}^{+})\to S\end{equation*}from the relative adic spectrum.\end{prop}
\begin{proof}Note that the second assertion (about the relative adic spectrum being the generic fiber of the relative formal spectrum) follows from the first assertion by definition of the relative adic spectrum, the relative formal spectrum and the adic analytic generic fiber. Hence it suffices to prove the first assertion (about the sheaves of algebras). 

By \cite{FK}, Ch.~I, Proposition~4.1.7, we can identify $\mathcal{A}_{0}$ with the direct image sheaf $f_{0\ast}\mathcal{O}_{\mathfrak{X}}$. By Lemma~\ref{Formal spectrum of a rig-sheafy algebra}, the adic formal $R$-scheme $\Spf(\mathcal{A}_{0})$ is locally rig-sheafy, so we can consider its adic analytic generic fiber $X$. By Lemma~\ref{Affine morphisms and locally affinoid morphisms}, the map of adic spaces $f=f_{0\eta}: X\to S$ is locally affinoid. By Proposition~\ref{Direct images 2}, this entails that the pair $(f_{\ast}\mathcal{O}_{X}, f_{\ast}\mathcal{O}_{X}^{+})$ is a quasi-coherent sheaf of adic Banach algebras on $S$. For a rig-sheafy affine open subset $\mathfrak{U}\subseteq \mathfrak{S}$, set $\mathfrak{V}=f_{0}^{-1}(\mathfrak{U})$, a rig-sheafy affine open subset of $\mathfrak{X}$. Let $U=\mathfrak{U}_{\eta}^{\ad}$ and $V=\mathfrak{V}_{\eta}^{\ad}$ and view them as affinoid open subspaces of $S$ and $X$, respectively, via the identifications $\mathfrak{U}_{\eta}^{\ad}=\spc_{S,\mathfrak{S}}^{-1}(\mathfrak{U})$ and $\mathfrak{V}_{\eta}^{\ad}=\spc_{X,\mathfrak{X}}^{-1}(\mathfrak{V})$ from \cite{Dine25-2}, Lemma~3.8. Then we have $\mathcal{A}_{0}(\mathfrak{U})[\varpi^{-1}]=\mathcal{O}_{\mathfrak{S}}(f_{0}^{-1}(\mathfrak{U}))[\varpi^{-1}]=\mathcal{O}_{X}(V)$ and \begin{equation*}\overline{\mathcal{A}_{0}(\mathfrak{U})}=\overline{\mathcal{O}_{\mathfrak{S}}(f_{0}^{-1}(\mathfrak{U}))}=\mathcal{O}_{X}^{+}(V)\end{equation*}by definition of the adic analytic generic fiber. But by the same commutative diagram as in the proof of Lemma~\ref{Affine morphisms and locally affinoid morphisms}, we have $V=f^{-1}(U)$. It follows that \begin{equation*}(\mathcal{A}, \mathcal{A}^{+})=(f_{\ast}\mathcal{O}_{X}, f_{\ast}\mathcal{O}_{X}^{+})\end{equation*}is a quasi-coherent sheaf of adic Banach algebras with the desired properties.  \end{proof} 
In the above situation, the formal model $\mathfrak{S}$ gives an abundance of $(\mathcal{A}, \mathcal{A}^{+})$-adapted affinoid open covers of $S$. 
\begin{lemma}\label{Adapted subsets and formal models}In the situation of Proposition~\ref{Relative formal spectrum and relative adic spectrum}, for every rig-sheafy affine open subset $\mathfrak{U}\subseteq \mathfrak{S}$, the affinoid open subspace \begin{equation*}U=\mathfrak{U}_{\eta}^{\ad}=\spc_{S,\mathfrak{S}}^{-1}(\mathfrak{U})\end{equation*}is $(\mathcal{A}, \mathcal{A}^{+})$-adapted. \end{lemma}
\begin{proof}By Theorem~\ref{The affinoid equivalence}, or by the proof of Proposition~\ref{Relative formal spectrum and relative adic spectrum} we may identify $(\mathcal{A}, \mathcal{A}^{+})$ with $(f_{\ast}\mathcal{O}_{X}, f_{\ast}\mathcal{O}_{X}^{+})$. By Corollary~\ref{Locally affinoid 3}, the affinoid open subspace $U$ being $(\mathcal{A}, \mathcal{A}^{+})$-adapted is equivalent to the pre-image $f^{-1}(U)$ being affinoid. But $f^{-1}(U)$ is affinoid by Lemma~\ref{Affine morphisms and locally affinoid morphisms}.  \end{proof}
\begin{prop}\label{Quasi-coherent adic Banach algebras and formal models}Let $\mathfrak{S}$ be a locally rig-sheafy $\varpi$-torsion-free quasi-compact quasi-separated adic formal $R$-scheme with adic analytic generic fiber $S=\mathfrak{S}_{\eta}^{\ad}$ and let $(\mathcal{A}, \mathcal{A}^{+})$ be a uniform quasi-coherent (topological) adic Banach $(\mathcal{O}_{S}, \mathcal{O}_{S}^{+})$-algebra. Suppose that $\mathfrak{S}$ has a rig-sheafy affine open cover $(\mathfrak{U}_{i})_{i}$ such that the pre-images \begin{equation*}U_{i}=\spc_{S,\mathfrak{S}}^{-1}(\mathfrak{U}_{i})\end{equation*}are $(\mathcal{A}, \mathcal{A}^{+})$-adapted affinoid open subspaces. Then the direct image sheaf $\spc_{S,\mathfrak{S}\ast}\mathcal{A}^{+}$ is an adically quasi-coherent sheaf of algebras on $\mathfrak{S}$ and \begin{equation*}\Spf(\spc_{S,\mathfrak{S}\ast}\mathcal{A}^{+})\to \mathfrak{S}\end{equation*}is a formal model of \begin{equation*}\underline{\Spa}_{S}(\mathcal{A}, \mathcal{A}^{+})\to S.\end{equation*}\end{prop}
\begin{proof}For every rig-sheafy affine open subset $\mathfrak{U}\subseteq \mathfrak{S}$, the pre-image $U=\spc_{S,\mathfrak{S}}^{-1}(\mathfrak{U})$ is an affinoid open subspace of $S$ (\cite{Dine25-2}, Lemma~3.8). Since the quasi-coherent adic Banach algebra $(\mathcal{A}, \mathcal{A}^{+})$ is uniform, \begin{equation*}\mathcal{A}^{+}(U)=\spc_{S,\mathfrak{S}\ast}\mathcal{A}^{+}(\mathfrak{U})\end{equation*}is a ring of definition of $\mathcal{A}(U)$. Hence, once we prove that $\spc_{S,\mathfrak{S}\ast}\mathcal{A}^{+}$ is adically quasi-coherent, it follows that \begin{equation*}\underline{\Spa}_{S}(\mathcal{A}, \mathcal{A}^{+})\to S\end{equation*}is indeed the generic fiber of \begin{equation*}\Spf(\spc_{S,\mathfrak{S}\ast}\mathcal{A}^{+})\to \mathfrak{S}.\end{equation*}

Let $(\mathfrak{U}_{i})_{i}$ be a rig-sheafy affine open cover of $\mathfrak{S}$. It suffices to prove that the sheaf\begin{equation*}(\spc_{S,\mathfrak{S}\ast}\mathcal{A}^{+})\vert_{\mathfrak{U}_{i}}=\spc_{U_{i},\mathfrak{U}_{i}\ast}(\mathcal{A}^{+}\vert_{U_{i}})\end{equation*}on $\mathfrak{U}_{i}$ (where in the above equality we tacitly used \cite{Dine25-2}, Lemma~3.8) is adically quasi-coherent for each $i$. Hence it suffices to prove the proposition in the case when $\mathfrak{S}$ is a rig-sheafy affine adic formal $R$-scheme and its generic fiber $S$ is $(\mathcal{A}, \mathcal{A}^{+})$-adapted (we use here that $\spc_{U_{i},\mathfrak{U}_{i}}^{-1}(\mathfrak{U}_{i})=(\mathfrak{U}_{i})_{\eta}^{\ad}$ for all $i$, by \cite{Dine25-2}, Lemma~3.8). For this it suffices to check that \begin{equation*}\mathcal{A}^{+}(S(\frac{1}{g}))=(\spc_{S,\mathfrak{S}\ast}\mathcal{A}^{+})(D(g))=\mathcal{A}^{+}(S)\widehat{\otimes}_{A}A\langle g^{-1}\rangle=\mathcal{A}^{+}(S)\langle g^{-1}\rangle\end{equation*}for every $g\in A$. Since $S$ is $(\mathcal{A}, \mathcal{A}^{+})$-adapted, $\mathcal{A}^{+}(S(\frac{1}{g}))=\mathcal{A}(S)\langle g^{-1}\rangle^{+}$, the integral closure of $\mathcal{A}^{+}(S)\langle g^{-1}\rangle$ inside $\mathcal{A}(S)\langle g^{-1}\rangle$. Hence it suffices to prove that $\mathcal{A}^{+}(S)\langle g^{-1}\rangle$ is integrally closed inside $\mathcal{A}(S)\langle g^{-1}\rangle$. By \cite{Huber0}, Lemma 2.4.3(iv), this follows from $\mathcal{A}^{+}(S)[g^{-1}]$ being an integrally closed open subring of $\mathcal{A}(S)[g^{-1}]$. \end{proof}
As a consequence of these results, we obtain the following generalization of \cite{Dine25-2}, Proposition~5.8. 
\begin{prop}\label{Normalization of a formal model}Let $X$ be a uniform adic space over $\Spa(R[\varpi^{-1}], \overline{R})$ and let \begin{equation*}f: X\to S\end{equation*}be a locally affinoid morphism of adic spaces over $\Spa(R[\varpi^{-1}], \overline{R})$. Let $\mathfrak{S}$ be a $\varpi$-torsion-free formal $R$-model of $S$ and suppose that there exists a rig-sheafy affine open cover $(\mathfrak{U}_{i})_{i}$ of $\mathfrak{S}$ such that, for every index $i$, the pre-image $f^{-1}(U_{i})$ of $U_{i}=\spc_{S,\mathfrak{S}}^{-1}(\mathfrak{U}_{i})$ is affinoid. Then there exists a unique affine adic morphism of locally rig-sheafy adic formal $R$-schemes $f_{0}: \mathfrak{X}\to\mathfrak{S}$ such that $\mathfrak{X}_{\eta}^{\ad}=X$, $f_{0\eta}=f$ and $\mathfrak{X}$ is $\varpi$-torsion-free and integrally closed in its generic fiber $X$.\end{prop}
\begin{proof}By Lemma~\ref{Locally affinoid 1.5}, the assumption that the pre-images $f^{-1}(U_{i})$ are affinoid implies that the affinoid open subspaces $U_{i}$ of $S$ are $(f_{\ast}\mathcal{O}_{X}, f_{\ast}\mathcal{O}_{X}^{+})$-adapted. Moreover, since $X$ is uniform, the quasi-coherent adic Banach algebra $(f_{\ast}\mathcal{O}_{X}, f_{\ast}\mathcal{O}_{X}^{+})$ is uniform. We have seen in Proposition~\ref{Quasi-coherent adic Banach algebras and formal models} that in this case $\spc_{S,\mathfrak{S}\ast}f_{\ast}\mathcal{O}_{X}^{+}$ is an adically quasi-coherent sheaf of algebras on $\mathfrak{S}$. Set \begin{equation*}\mathfrak{X}=\Spf(\spc_{S,\mathfrak{S}\ast}f_{\ast}\mathcal{O}_{X}^{+})\end{equation*}and let $f_{0}: \mathfrak{X}\to\mathfrak{S}$ be the canonical morphism, which is adic and affine. By the last assertion of Proposition~\ref{Quasi-coherent adic Banach algebras and formal models}, we have $\mathfrak{X}_{\eta}^{\ad}=X$ and $f_{0\eta}=f$.  By construction, $\mathcal{O}_{\mathfrak{X}}(f_{0}^{-1}(\mathfrak{U}_{i}))=\mathcal{O}_{X}^{+}(f^{-1}(U_{i}))$, so $\mathcal{O}_{\mathfrak{X}}(f_{0}^{-1}(\mathfrak{U}_{i}))$ is an integrally closed subring of \begin{equation*}\mathcal{O}_{\mathfrak{X}}(f_{0}^{-1}(\mathfrak{U}_{i}))[\varpi^{-1}]=\mathcal{O}_{X}(f^{-1}(U_{i}))\end{equation*}for all $i$. Since $X$ is uniform, $\mathcal{O}_{\mathfrak{X}}(f_{0}^{-1}(\mathfrak{U}_{i}))=\mathcal{O}_{X}^{+}(f^{-1}(U_{i}))$ is a ring of definition of the Tate ring $\mathcal{O}_{X}(f^{-1}(U_{i}))$ for all $i$. In particular, $\mathfrak{X}$ is an integrally closed formal $R$-model of $X$. 

It remains to prove the uniqueness part of the proposition. To this end, let $g_{0}: \mathfrak{Z}\to \mathfrak{S}$ be another affine adic morphism of $\varpi$-torsion-free adic formal $R$-schemes such that $\mathfrak{Z}_{\eta}^{\ad}=X$, $g_{0\eta}=f$ and $\mathfrak{Z}$ is integrally closed in its generic fiber. By \cite{FK}, Ch.~I, Theorem~4.1.8, $g_{0\ast}\mathcal{O}_{\mathfrak{Z}}$ is an adically quasi-coherent $\mathcal{O}_{\mathfrak{S}}$-algebra, $\mathfrak{Z}=\Spf(g_{0\ast}\mathcal{O}_{\mathfrak{Z}})$ and $g_{0}$ is the canonical morphism $\Spf(g_{0\ast}\mathcal{O}_{\mathfrak{Z}})\to \mathfrak{S}$. By loc.~cit., Proposition~4.1.10, the pre-image $g_{0}^{-1}(\mathfrak{U})$ of every affine open subset $\mathfrak{U}\subseteq\mathfrak{S}$ is affine. Let $\mathfrak{U}$ be a rig-sheafy affine open subset of $\mathfrak{S}$ and let $U=\spc_{S,\mathfrak{S}}^{-1}(\mathfrak{U})$. By \cite{Dine25-2}, Lemma~5.16, the generic fiber $X$ being uniform implies that \begin{equation*}\mathfrak{Z}=\Spf(g_{0\ast}\mathcal{O}_{\mathfrak{Z}})\end{equation*}is locally stably uniform. Thus the adically quasi-coherent algebra $g_{0\ast}\mathcal{O}_{\mathfrak{Z}}$ is rig-sheafy (and, in fact, uniform) by Lemma~\ref{Formal spectrum of a rig-sheafy algebra}. By Lemma~\ref{Affine morphisms and locally affinoid morphisms} applied to the affine morphism $g_{0}$ and its generic fiber $g_{0\eta}=f$, the pre-image $f^{-1}(U)$ is affinoid. Since $X=\mathfrak{Z}_{\eta}^{\ad}$ and $f=g_{0\eta}$, the commutative square \begin{center}\begin{tikzcd}X \arrow{r}{f} \arrow{d}{\spc_{X,\mathfrak{Z}}} & S \arrow{d}{\spc_{S,\mathfrak{S}}} \\ \mathfrak{Z} \arrow{r}{g_{0}} & \mathfrak{S} \end{tikzcd}\end{center}(together with \cite{Dine25-2}, Lemma~3.8) shows that $g_{0}^{-1}(\mathfrak{U})_{\eta}^{\ad}$ is equal to the affinoid open subspace \begin{equation*}f^{-1}(U)=\Spa(\mathcal{O}_{X}(f^{-1}(U)), \mathcal{O}_{X}^{+}(f^{-1}(U)))\end{equation*}of $X$. Therefore, we have \begin{equation*}\mathcal{O}_{\mathfrak{Z}}(g_{0}^{-1}(\mathfrak{U}))[\varpi^{-1}]=\mathcal{O}_{X}(f^{-1}(U)),\end{equation*}and the integral closure of $\mathcal{O}_{\mathfrak{Z}}(g_{0}^{-1}(\mathfrak{U}))$ inside $\mathcal{O}_{\mathfrak{Z}}(g_{0}^{-1}(\mathfrak{U}))[\varpi^{-1}]$ is equal to $\mathcal{O}_{X}^{+}(f^{-1}(U))$. Since $\mathfrak{Z}$ is integrally closed in its generic fiber, we have \begin{equation*}\mathcal{O}_{\mathfrak{Z}}(g_{0}^{-1}(\mathfrak{U}))=\mathcal{O}_{X}^{+}(f^{-1}(U)).\end{equation*}Since $\mathfrak{U}$ was an arbitrary rig-sheafy affine open subset, this means that \begin{equation*}g_{0\ast}\mathcal{O}_{\mathfrak{Z}}=\spc_{S,\mathfrak{S}\ast}f_{\ast}\mathcal{O}_{X}^{+}\end{equation*}as adically quasi-coherent algebras on $\mathfrak{S}$. By \cite{FK}, Ch.~I, Theorem 4.1.8, this implies the desired equality $g_{0}=f_{0}$.\end{proof}

\subsection{Proper locally affinoid morphisms and universally closed affine morphisms}

We can now characterize proper locally affinoid morphisms of Tate adic spaces in terms of formal models. 
\begin{prop}\label{Partially proper morphisms and formal models}Let $f_{0}: \mathfrak{X}\to \mathfrak{S}$ be an affine morphism of locally rig-sheafy quasi-compact quasi-separated $\varpi$-torsion-free adic formal $R$-schemes with $\mathfrak{X}$ locally stably uniform and let \begin{equation*}f=f_{0\eta}: X=\mathfrak{X}_{\eta}^{\ad}\to S=\mathfrak{S}_{\eta}^{\ad}\end{equation*}be the induced morphism between the adic analytic generic fibers over $\Spa(R[\varpi^{-1}], \overline{R})$. Then $f$ is proper if and only if $f_{0}$ is universally closed.\end{prop}
\begin{proof}Suppose that $f$ is proper. Let $\mathfrak{U}$ be a rig-sheafy affine open subset of $\mathfrak{S}$ with generic fiber $U=\spc_{S,\mathfrak{S}}^{-1}(\mathfrak{U})$, and consider its pre-image $f_{0}^{-1}(\mathfrak{U})$ under $f_{0}$, an affine open subset of $\mathfrak{X}$. Since $f_{0}$ is affine, $\mathfrak{X}=\Spf(f_{0\ast}\mathcal{O}_{\mathfrak{X}})$ by \cite{FK}, Ch.~I, Theorem~4.1.8. Since $\mathfrak{X}$ is locally stably uniform, Lemma~\ref{Formal spectrum of a rig-sheafy algebra} ensures that $f_{0}^{-1}(\mathfrak{U})=\Spf(\mathcal{O}_{\mathfrak{X}}(f_{0}^{-1}(\mathfrak{U})))$ is rig-sheafy. Its adic analytic generic fiber is identified with $\spc_{X,\mathfrak{X}}^{-1}(f_{0}^{-1}(\mathfrak{U}))$ via \cite{Dine25-2}, Lemma~3.8, so, by the commutative diagram \begin{center}\begin{tikzcd}X \arrow{r}{f} \arrow{d}{\spc_{X,\mathfrak{X}}} & S \arrow{d}{\spc_{S,\mathfrak{S}}} \\ \mathfrak{X} \arrow{r}{f_{0}} & \mathfrak{S}, \end{tikzcd}\end{center}it is identified with $f^{-1}(U)$. By Lemma~\ref{Partial properness is local on the target}, the restriction $f\vert_{f^{-1}(U)}: f^{-1}(U)\to U$ is proper. Since $f^{-1}(U)$ is the generic fiber of the rig-sheafy affine $f_{0}^{-1}(\mathfrak{U})$, it is affinoid, so $f\vert_{f^{-1}(U)}$ is a proper morphism between affinoid analytic adic spaces. By Lemma~\ref{Partial properness for affinoids} this means that $\mathcal{O}_{X}^{+}(f^{-1}(U))$ is the smallest open integrally closed $\mathcal{O}_{S}^{+}(U)$-subalgebra of $\mathcal{O}_{X}(f^{-1}(U))$. This entails that\begin{equation*}\mathcal{O}_{X}^{+}(f^{-1}(U))=\overline{\mathcal{O}_{S}^{+}(U)[\varpi\mathcal{O}_{\mathfrak{X}}(f_{0}^{-1}(\mathfrak{U}))]},\end{equation*}the ring $\mathcal{O}_{\mathfrak{X}}(f_{0}^{-1}(\mathfrak{U}))$ being a ring of definition of $\mathcal{O}_{X}(f^{-1}(U))$ (indeed, if $B$ is an open integrally closed subring of $\mathcal{O}_{X}(f^{-1}(U))$, then $B$ contains $\varpi^{n}\mathcal{O}_{\mathfrak{X}}(f_{0}^{-1}(\mathfrak{U}))$ for some $n$, in which case it contains $(\varpi f)^{n}$ for every element $f\in\mathcal{O}_{\mathfrak{X}}(f_{0}^{-1}(\mathfrak{U}))$; but then $\varpi f\in B$, since $B$ was integrally closed). Since $\mathcal{O}_{S}^{+}(U)=\overline{\mathcal{O}_{\mathfrak{S}}(\mathfrak{U})}$, the above displayed formula means that \begin{equation*}\mathcal{O}_{X}^{+}(f^{-1}(U))=\overline{\mathcal{O}_{\mathfrak{S}}(\mathfrak{U})[\varpi\mathcal{O}_{\mathfrak{X}}(f_{0}^{-1}(\mathfrak{U}))]}=\overline{\mathcal{O}_{\mathfrak{S}}(\mathfrak{U})+\varpi\mathcal{O}_{\mathfrak{X}}(f_{0}^{-1}(\mathfrak{U}))}.\end{equation*}In particular, we see that the map \begin{equation*}\mathcal{O}_{\mathfrak{S}}(\mathfrak{U})/\varpi\mathcal{O}_{\mathfrak{S}}(\mathfrak{U})\to\mathcal{O}_{\mathfrak{X}}(f_{0}^{-1}(\mathfrak{U}))/\varpi\mathcal{O}_{\mathfrak{X}}(f_{0}^{-1}(\mathfrak{U}))\end{equation*}is integral. In particular, the morphism of affine schemes \begin{equation*}f_{0}^{-1}(\mathfrak{U})_{0}\to \mathfrak{U}_{0}\end{equation*}is universally closed. Since $\mathfrak{U}\subseteq\mathfrak{S}$ was an arbitrary rig-sheafy affine open, we conclude that the morphism of schemes $\mathfrak{X}_{0}\to\mathfrak{S}_{0}$ induced by $f_0$ is universally closed. By \cite{FK}, Ch.~I, Proposition~4.5.9, this is equivalent to the morphism of formal schemes $f_{0}: \mathfrak{X}\to\mathfrak{S}$ being universally closed. 

Conversely, suppose that $f_{0}$ is a universally closed affine morphism of formal schemes, let $(\mathfrak{U}_{i})_{i}$ be a rig-sheafy affine open cover of $\mathfrak{S}$ and let $(U_{i})_{i}$ be the corresponding affinoid open cover of the generic fiber $S$ (i.e., $U_{i}=\spc_{S,\mathfrak{S}}^{-1}(\mathfrak{U}_{i})$). By Lemma~\ref{Partial properness is local on the target}, to prove that $f=f_{0\eta}$ is proper, it suffices to prove that $f\vert_{f^{-1}(U_{i})}: f^{-1}(U_{i})\to U_{i}$ is proper for every $i$. Hence we may assume that $\mathfrak{S}$ is rig-sheafy and affine, in which case $\mathfrak{X}$ is also rig-sheafy and affine (it is affine since $f_{0}$ is an affine morphism and it is rig-sheafy since $f_{0\ast}\mathcal{O}_{\mathfrak{X}}$ is uniform by Lemma~\ref{Formal spectrum of a rig-sheafy algebra}). In this case, the morphism of affine schemes $\mathfrak{X}_{0}\to \mathfrak{S}_{0}$ induced by $f_{0}$ is integral. In other words, the ring map \begin{equation*}\mathcal{O}_{\mathfrak{S}}(\mathfrak{S})/\varpi\mathcal{O}_{\mathfrak{S}}(\mathfrak{S})\to \mathcal{O}_{\mathfrak{X}}(\mathfrak{X})/\varpi\mathcal{O}_{\mathfrak{X}}(\mathfrak{X})\end{equation*}is integral. It follows that \begin{equation*}\mathcal{O}_{X}^{+}(X)=\overline{\mathcal{O}_{\mathfrak{X}}(\mathfrak{X})}=\overline{\mathcal{O}_{\mathfrak{S}}(\mathfrak{S})+\varpi\mathcal{O}_{\mathfrak{X}}(\mathfrak{X})}=\overline{\mathcal{O}_{S}^{+}(S)+\varpi\mathcal{O}_{\mathfrak{X}}(\mathfrak{X})},\end{equation*}the smallest integrally closed and open $\mathcal{O}_{S}^{+}(S)$-subalgebra of $\mathcal{O}_{X}(X)$. By Lemma \ref{Partial properness for affinoids}, this entails that the morphism $f$ is proper. \end{proof}

\bibliographystyle{plain} 
\bibliography{Bib}
\end{document}